\documentclass[notheoremnums]{dpreprint}

\newtheorem{theorem}{Theorem}[section]
\newtheorem*{theorem*}{Theorem}
\newtheorem*{proposition*}{Proposition}

\newtheorem{proposition}[theorem]{Proposition}

\newtheorem{lemma}[theorem]{Lemma}
\newtheorem{notation}[theorem]{Notation}
\theoremstyle{definition}
\newtheorem{definition}[theorem]{Definition}

\newtheorem*{definition*}{Definition}

\newtheorem{maintheorem}{Theorem}

\newcommand{\muB}{\lambda_{B}}

\newcommand{\len}[1]{\| #1 \|}

\usepackage{xcolor}
\usepackage{mathtools}
\DeclarePairedDelimiter{\ceiling}{\lceil}{\rceil}
\DeclarePairedDelimiter{\Bigfloor}{\Big{\lfloor}}{\Big{\rfloor}}

\begin{document}

\dtitle[Mixing subshifts of minimal word complexity]{Measure-Theoretically Mixing Subshifts of Minimal Word Complexity}
\dauthor[Darren~Creutz]{Darren Creutz}{darren.a.creutz@vanderbilt.edu}{Vanderbilt University, 1326 Stevenson Center Lane, Nashville TN 37240}{}
%\dnoauthor
\datewritten{\today}


\dabstract{
We resolve a long-standing open question on the relationship between measure-theoretic dynamical complexity and symbolic complexity by establishing the exact word complexity at which measure-theoretic strong mixing manifests:

For every superlinear $f : \mathbb{N} \to \mathbb{N}$, i.e.~$\nicefrac{f(q)}{q} \to \infty$, there exists a subshift admitting a (strongly) mixing of all orders probability measure with word complexity $p$ such that $\nicefrac{p(q)}{f(q)} \to 0$.

For a subshift with word complexity $p$ which is non-superlinear, i.e.~$\liminf \nicefrac{p(q)}{q} < \infty$, every ergodic probability measure is partially rigid.
}

\makepreprint

\dsectionstar{Introduction}

Among measure-theoretic dynamical properties of measure-preserving transformations, strong mixing of all orders is the `most complex': every finite collection of measurable sets tends asymptotically toward independence, necessarily implying a significant amount of randomness.  Despite this, `low complexity' mixing transformations exist--there are mixing transformation with zero entropy--raising the question of how deterministic a mixing transformation can be.

Word complexity, the number $p(q)$ of distinct words of length $q$ appearing in the language of the subshift, provides a more fine-grained means of quantifying complexity in the zero entropy setting, leading to the question of how low the word complexity of a mixing transformation can be.

Ferenczi \cite{ferenczichacon} initially conjectured that mixing transformations' word complexity should be superpolynomial but quickly refuted this himself \cite{ferenczi1996rank} showing that the staircase transformation, proven mixing by Adams \cite{adams1998smorodinsky}, has quadratic word complexity.
Recent  joint work of the author and R.~Pavlov and S.~Rodock \cite{CPR} exhibited subshifts admitting mixing measures with word complexity functions which are subquadratic but superlinear by more than a logarithm.  We exhibit subshifts admitting mixing measures with complexity arbitrarily close to linear:

\begin{maintheorem}
For every $f : \mathbb{N} \to \mathbb{N}$ which is superlinear, $\nicefrac{f(q)}{q} \to \infty$, there exists a subshift, admitting a strongly mixing probability measure, with word complexity $p$ such that $\nicefrac{p(q)}{f(q)} \to 0$.
\end{maintheorem}

Our examples, which we call quasi-staircase transformations, are mixing rank-one transformations hence mixing of all orders \cite{kalikow}, \cite{Ry93}.  We establish their word complexity is optimal:

\begin{maintheorem}\label{mainB}
Every subshift of non-superlinear word complexity, $\liminf \nicefrac{p(q)}{q} < \infty$, equipped with an ergodic probability measure is partially rigid hence not strongly mixing,
\end{maintheorem}

Non-superlinear complexity subshifts are conjugate to $S$-adic shifts (Donoso, Durand, Maass and Petite \cite{DDMP2}).  Named by Vershik and the subject of a well-known conjecture of Host, $S$-adic subshifts are quite structured (see e.g.~\cite{leroy2} for more information on $S$-adicity).

Our work may be viewed as saying there is a sharp divide in `measure-theoretic complexity', precisely at superlinear word complexity, between highly structured and highly complicated: as soon as the word complexity is  `large enough' to escape the $S$-adic structure and partial rigidity, there is already `enough room' for (strong) mixing of all orders.

Cyr and Kra established that superlinear complexity is the dividing line for a subshift admitting countably many ergodic measures: there exists subshifts with complexity arbitrarily close to linear which admit uncountably many ergodic measures \cite{CK4} and non-superlinear complexity implies at most countably many \cite{CK2}, \cite{boshernitzan}.  Our work implies that in the non-superlinear case, the at most countably many measures are all partially rigid (with a uniform rigidity constant).  Their result, like ours, indicates that superlinear word complexity is the line at which complicated measure-theoretic phenomena can manifest.

Beyond the structure imposed by $S$-adicity, linear complexity subshifts are known to be structured in various ways (e.g.~\cite{cassaigneetal}, \cite{CK3}, \cite{DDMP}, \cite{DOP}, \cite{PS2}, \cite{PS}).  Our work indicates there is no hope for similar phenomena in any superlinear setting.

\section{Definitions and preliminaries}

\subsection{Symbolic dynamics}

\begin{definition}
A \textbf{subshift} on the finite set $\mathcal{A}$ is any subset $X \subset \mathcal{A}^{\mathbb{Z}}$ which is closed in the product topology and shift-invariant: for all $x = (x_{n})_{n \in \mathbb{Z}} \in X$ and $k \in \mathbb{Z}$, the translate $(x_{n+k})_{n \in \mathbb{Z}}$ of $x$ by $k$ is also in $X$.
\end{definition} 

\begin{definition}
A \textbf{word} is any element of $\mathcal{A}^\ell$ for some $\ell$, the \textbf{length} of $w$, written $\len{w}$.  A word $w$ is a \textbf{subword} of a word or bi-infinite sequence $x$ if there exists $k$ so that $w_{i} = x_{i+k}$ for all $1 \leq i \leq \len{w}$.  A word $u$ is a \textbf{prefix} of $w$ when $u_{i} = w_{i}$ for $1 \leq i \leq \len{u}$ and a word $v$ is a suffix of $w$ when $v_{i} = w_{i+\len{w}-\len{v}}$ for $1 \leq i \leq \len{v}$.
\end{definition}

For words $v,w$, we denote by $vw$ their concatenation--the word obtained by following $v$ immediately by $w$. We write such concatenations with product or exponential notation, e.g.~$\prod_i w_i$ or $0^n$.

\begin{definition}
The \textbf{language} of a subshift $X$ is $\mathcal{L}(X) = \{ w : \text{$w$ is a subword of some $x \in X$} \}$.
\end{definition}

\begin{definition}
The \textbf{word complexity function} of a subshift $X$ over $\mathcal{A}$ is the function $p_X: \mathbb{N} \rightarrow \mathbb{N}$ defined by $p_X(q) = |\mathcal{L}(X) \cap \mathcal{A}^q|$, the number of words of length $q$ in the language of $X$.
\end{definition}

When $X$ is clear from context, we suppress the subscript and just write $p(n)$.

For subshifts on the alphabet $\{0,1\}$, we consider:

\begin{definition}
The set of \textbf{right-special} words is
$
\mathcal{L}^{RS}(X) = \{ w \in X : w0, w1 \in \mathcal{L}(X) \}
$.
\end{definition}

Cassaigne \cite{cassaigne} showed the well-known: $
p(q) = p(m) + \sum_{\ell = m}^{q-1} |\{ w \in \mathcal{L}^{RS} : \len{w} = \ell \}|
$
for $m < q$.

\subsection{Ergodic theory}

\begin{definition}
A \textbf{transformation} $T$ is a measurable map on a standard Borel or Lebesgue measure space $(Y,\mathcal{B},\mu)$ that is measure-preserving: $\mu(T^{-1}B) = \mu(B)$ for all $B \in \mathcal{B}$.
\end{definition}

\begin{definition}Two transformations $T$ on $(Y,\mathcal{B},\mu)$ and $T^{\prime}$ on $(Y^{\prime},\mathcal{B}^{\prime},\mu^{\prime})$ are \textbf{measure-theoretically isomorphic} when there exists 
a bijective map $\phi$ between full measure subsets $Y_0 \subset Y$ and $Y'_0 \subset Y'$ where 
$\mu(\phi^{-1} A) = \mu'(A)$ for all measurable $A \subset Y'_0$ and $(\phi \circ T)(y) = (T' \circ \phi)(y)$ for all $y \in Y_0$.
\end{definition}

\begin{definition}
A transformation $T$ is \textbf{ergodic} when $A = T^{-1} A$ implies that $\mu(A) = 0$ or $\mu(A^c) = 0$. 
\end{definition}

\begin{theorem}[Mean Ergodic Theorem]
If $T$ is ergodic and on a finite measure space and $f \in L^2(Y)$, 
\[
\lim_{n \rightarrow \infty} \int \big{|}\frac{1}{n} \sum_{i = 0}^{n-1} f \circ T^{-i} - \int f~d\mu \big{|} \ d\mu = 0
\]
\end{theorem}

\begin{definition}
A transformation $T$ is \textbf{mixing} when for all $A, B \in \mathcal{B}$, $\mu(T^{n} A \cap B) \rightarrow \mu(A) \mu(B)$.
\end{definition}

\subsection{Rank-one transformations}

A \textbf{rank-one transformation} is a transformation $T$ constructed by ``cutting and stacking".  Here $Y$ represents a (possibly infinite) interval, $\mathcal{B}$ is the induced $\sigma$-algebra from $\mathbb{R}$, and $\mu$ is Lebesgue measure. We give a brief description, referring the reader to 
\cite{fghsw21} or  \cite{silva2008invitation} for more details.

The transformation is defined inductively on larger and larger portions of the space through Rohlin towers or \textbf{columns}, denoted $C_n$. Each column $C_n$ consists of \textbf{levels} $I_{n,j}$ where $0 \leq j < h_{n}$ is the height of the level within the column. All levels $I_{n,j}$ in $C_n$ are 
intervals with the same length, $\mu(I_{ n})$, and the total number of levels in a column is the \textbf{height} of the column, denoted by $h_n$. The transformation $T$ is defined on all levels $I_{n,j}$ except the top one $I_{n, h_n - 1}$ by sending each $I_{n,j}$ to $I_{n,j+1}$ using the unique order-preserving affine map.

Start with $C_1=[0,1)$ with height $h_1=1$. To obtain $C_{n+1}$ from $C_n$, we require a \textbf{cut sequence}, $\{r_n\}$ such that $r_n \geq 1$ for all $n$. Make $r_n$ vertical cuts of $C_n$ to create $r_n+1$ \textbf{subcolumns} of equal width. Denote a \textbf{sublevel} of $C_n$ by $I_{n,j}^{[i]}$ where $0 \leq a < h_{n}$ is the height of the level within that column, and 
$i$ represents the position of the subcolumn, where $i=0$ represents the leftmost subcolumn and $i=r_n$ is the rightmost subcolumn. After cutting $C_n$ into subcolumns, add extra intervals called \textbf{spacers} on top of each subcolumn to function as levels of the next column. The \textbf{spacer sequence}, $\{s_{n,i}\}$ such that $0 \leq i \leq r_{n}$ and $s_{n,i} \geq 0$, specifies how many sublevels to add above each subcolumn. Spacers are the same width as the sublevels, act as new levels in the column $C_{n+1}$, and are taken to be the leftmost intervals in $[1,\infty)$ not in $C_{n}$.
After the spacers are added, stack the subcolumns with their spacers right on top of left, i.e.~so that $I_{n,0}^{[i+1]}$ is directly above $I_{n,h_{n}-1}^{[i]}$. This gives the next column, $C_{n+1}$.

Each column $C_n$ defines $T$ on $\bigcup_{j = 0}^{h_n - 2} I_{n,j}$ and 
the partially defined map $T$ on $C_{n+1}$ agrees with that of $C_n$, extending the definition of $T$ to a portion of the top level of $C_n$ where it was previously undefined. Continuing this process gives the sequence of columns $\{C_1, \dots, C_n, C_{n+1}, \dots\}$ and $T$ is then the limit of the partially defined maps. 

Though this construction could result in $Y$ being an infinite interval with infinite Lebesgue measure, $Y$ has finite measure if and only if 
$\sum_{n} \frac{1}{r_{n}h_{n}}\sum_{i=0}^{r_{n}} s_{n,i} < \infty$, see \cite{CreutzSilva2010}.
All rank-one transformations we define satisfy this condition, and for convenience we renormalize so that $Y = [0,1)$. Every rank-one transformation is ergodic and invertible.

The reader should be aware that we are making $r_{n}$ cuts and obtaining $r_{n}+1$ subcolumns (following Ferenczi \cite{ferenczi1996rank}), while other papers (e.g.~\cite{Creutz2021}) use $r_{n}$ as the number of subcolumns.

\subsection{Symbolic models of rank-one transformations}

For a rank-one transformation defined as above, we define a subshift $X(T)$ on the alphabet $\{ 0,1 \}$ which is measure-theoretically isomorphic to $T$:

\begin{definition}
The \textbf{symbolic model} $X(T)$ of a rank-one transformation $T$ is given by the sequence of words: $B_1 = 0$ and 
\begin{equation*}\label{rk1word}
B_{n+1}=B_n1^{s_{n,0}}B_n1^{s_{n,1}}\cdots B_{n} 1^{s_{n,r_n}}= \prod_{i=0}^{r_n} B_n1^{s_{n,i}}
\end{equation*}
and $X(T)$ is the set of all bi-infinite sequences such that every subword is a subword of some $B_n$.
\end{definition}

The words $B_n$ are a symbolic coding of the column $C_n$: $0$ represents $C_1$ and $1$ represents the spacers.  There is a natural measure associated to $X(T)$:

\begin{definition}
The \textbf{empirical measure} for a symbolic model $X(T)$ of a rank-one transformation $T$ is the measure $\nu$ defined by, for each word $w$,
\[
\nu([w]) = \lim_{n \rightarrow \infty} \frac{|\{1 \leq j \leq \len{B_{n}} - \len{w} \ : \ {B_n}_{[j,j+\len{w})} = w\}|}{\len{B_{n}} - \len{w}}
\]
\end{definition}

Danilenko \cite{danilenkopr} (combined with \cite{deljunco} and \cite{kalikow}) proved that the symbolic model $X(T)$ of a rank-one subshift, equipped with its empirical measure, is measure-theoretically isomorphic to the cut-and-stack construction (see \cite{adamsferenczipeterson17}; see \cite{fghsw21} for the full generality including odometers).

Due to this isomorphism, we move back and forth between rank-one and symbolic model terminology as needed and write $\mathcal{L}(T)$ for the language of $X(T)$.

\section{Quasi-staircase transformations}

\begin{definition}\label{qsdef} Given nondecreasing sequences of integers $\{ a_{n} \}$, $\{ b_{n} \}$ and $\{ c_{n} \}$ tending to infinity such that $c_{1} \geq 1$ and $c_{n+1} \geq c_{n} + b_{n}$, a \textbf{quasi-staircase transformation} is a rank-one transformation with cut sequence $r_{n} = a_{n}b_{n}$ and spacer sequence $s_{n,t} = c_{n} + \floor*{\frac{t}{a_{n}}}$ for $0 \leq t < r_{n}$ and $s_{n,r_{n}} = 0$.
\end{definition}

The symbolic representation of a quasi-staircase is $B_{1} = 0$ and
\[
B_{n+1} = \Big{(}\prod_{i=0}^{b_{n}-1}\big{(}B_{n}1^{c_{n}+i}\big{)}^{a_{n}}\Big{)}B_{n}
\]

The height sequence of a quasi-staircase is $h_{1} = 1$ and $h_{n+1} = (a_{n}b_{n}+1)h_{n} + a_{n}b_{n}c_{n} + \frac{1}{2}a_{n}b_{n}(b_{n}-1)$.

\subsection{Quasi-staircase right-special words}

\begin{lemma}\label{2ss}
The following hold:
\begin{enumerate}[\hspace{10pt}(i)\hspace{5pt}]
\item\label{2s1} $1^{\ell} \in \mathcal{L}^{RS}(T)$ for all $\ell$.
\item\label{2s0} If $w$ is a suffix of $1^{c_n}(B_{n}1^{c_{n}})^{a_{n}}$ then $w \in \mathcal{L}^{RS}(T)$. 
\item\label{2si} If $w$ is a suffix of $1^{c_{n}+i-1}(B_{n}1^{c_{n}+i})^{a_{n}}$ for $0 < i < b_{n}$ then $w \in \mathcal{L}^{RS}(T)$.
\item\label{2sspecial} If $w$ is a suffix of $1^{c_{n}+b_{n}-1}B_{n}1^{c_{n}}$ then $w \in \mathcal{L}^{RS}(T)$.
\end{enumerate}
\end{lemma}
\begin{proof}
Since suffixes of right-special words are right-special, it suffices to show the words $w$ is claimed to be a suffix of are right-special.

\textit{(\ref{2s1})}: For $n$ such that $\ell < c_{n}$, as the word $1^{c_{n}}B_{n}$ is a subword of $B_{n+1}$, so are $1^{\ell+1}$ and $1^{\ell}0$ since $\ell < c_{n}$ and $B_{n}$ starts with $0$.

\textit{(\ref{2s0})}:  $B_{n+2}$ has $1^{c_{n+1}}B_{n+1} = 1^{c_{n+1}-c_{n}}1^{c_{n}}B_{n+1}$ as a subword which has $1^{c_{n}}(B_{n}1^{c_{n}})^{a_{n}}B_{n}$ as a subword which gives $1^{c_{n}}(B_{n}1^{c_{n}})^{a_{n}}0 \in \mathcal{L}(T)$.  $B_{n+1}$ has $(B_{n}1^{c_{n}})^{a_{n}}B_{n}1^{c_{n}+1}$ as a prefix which has suffix $1^{c_{n}}(B_{n}1^{c_{n}})^{a_{n}-1}B_{n}1^{c_{n}+1}$ and that word is $1^{c_{n}}(B_{n}1^{c_{n}})^{a_{n}}1$ giving $1^{c_{n}-1}(B_{n}1^{c_{n}})^{a_{n}}1 \in \mathcal{L}(T)$.

\textit{(\ref{2si})}:  $B_{n+1}$ has $1^{c_{n}+i-1}(B_{n}1^{c_{n}+i})^{a_{n}}B_{n}$ as a subword which gives $1^{c_{n}+i-1}(B_{n}1^{c_{n}+i})^{a_{n}}0 \in \mathcal{L}(T)$.  When $i < b_{n} - 1$, $B_{n+1}$ has $(1^{c_{n}+i}B_{n})^{a_{n}}1^{c_{n}+i+1}$ as a subword which gives $1^{c_{n}+i-1}(B_{n}1^{c_{n}+i})^{a_{n}}1 \in \mathcal{L}(T)$; when $i = b_{n}-1$, $B_{n+2}$ has the subword $(1^{c_{n}+b_{n}-1}B_{n})^{a_{n}}1^{c_{n+1}}$ so $1^{c_{n}+b_{n}-2}(B_{n}1^{c_{n}+b_{n}-1})^{a_{n}}1^{c_{n+1}-c_n-b_n+1} \in \mathcal{L}(T)$ so $1^{c_{n}+b_{n}-2}(B_{n}1^{c_{n}+b_{n}-1})^{a_{n}}1 \in \mathcal{L}(T)$ as $c_{n+1} \geq c_{n} + b_{n}$.

\textit{(\ref{2sspecial})}: $B_{n+2}$ has $B_{n+1}1^{c_{n+1}}B_{n+1}$ as a subword which has $B_{n+1}1^{c_{n+1}}B_{n}1^{c_{n}}B_{n}$ as a prefix, and that word has $1^{c_{n}+b_{n}-1}B_{n}1^{c_{n}}0$ as a subword since $c_{n}+b_{n}-1 < c_{n+1}$.  Also $B_{n+2}$ has $B_{n+1}1^{c_{n+1}}$ as a subword which has  $1^{c_{n}+b_{n}-1}B_{n}1^{c_{n+1}}$ as a suffix which then has $1^{c_{n}+b_{n}-1}B_{n}1^{c_{n}}1$ as a subword.
\end{proof}

\begin{lemma}\label{L:010}
Let $01^{z}0 \in \mathcal{L}(T)$.  Then there are unique $n$ and $i$ with $0 \leq i < b_{n}$ such that $z = c_{n} + i$.   $01^{c_{n}+i}0$ is not a subword of $B_{m}$ for $m \leq n$ and for every $x \in X(T)$ and every occurrence of $01^{c_{n}+1}0$ in $x$, $01^{c_{n}+i}0$ occurs as a suffix of $1^{c_{n+1}}(\prod_{j=0}^{i-1}(B_{n}1^{c_{n}+j})^{a_{n}})(B_{n}1^{c_{n}+i})^{q}0$ for some $1 \leq q \leq a_{n}$ (adopting the convention that $\prod_{0}^{-1}$ is the empty word).
\end{lemma}
\begin{proof}
As every $B_{n}$ begins and ends with $0$, the only such words are of the form $01^{c_{n}+i}0$.  Since $c_{n+1} \geq c_{n} + b_{n}$, such $n$ and $i$ are unique.  This also gives that $1^{c_{n}}$ is not a subword of $B_{n}$.

The word $01^{c_{n}+i}0$ only occurs inside $B_{n+1}$ due to $c_{n+1} \geq c_{n} + b_{n}$, and only as part of the $(B_{n}1^{c_{n}+i})^{a_{n}}$ in its construction, and $B_{n+1}$ is always preceded by $1^{c_{n+1}}$
\end{proof}

\begin{proposition}\label{P:twosuccs}
If $w \in \mathcal{L}^{RS}(T)$ then at least one of the following holds:
\begin{enumerate}[\hspace{10pt}(i)\hspace{5pt}]
\item\label{Pt-1} $w = 1^{\len{w}}$
\item\label{Pt-3} $w$ is a suffix of $1^{c_{n}+i-1}\big{(}B_{n}1^{c_{n}+i}\big{)}^{a_{n}}$ for some $n$ and $0 \leq i < b_{n}$ 
\item\label{Pt-2} $w$ is a suffix of $1^{c_{n}+b_{n}-1}B_{n}1^{c_{n}}$ for some $n$ 
\item\label{Pt-weird} $w = 1^{c_{n}}(B_{n}1^{c_{n}})^{a_{n}}$
\end{enumerate}
\end{proposition}
\begin{proof}
Let $w \in \mathcal{L}^{RS}(T)$.  Since $c_{1} \geq 1$, the word $00 \notin \mathcal{L}(T)$ so $w$ does not end in $0$.  If $w = 1^{\len{w}}$ then $w$ is of form $(\ref{Pt-1})$ so from here on, assume that $w$ contains at least one $0$.

Let $z \geq 1$ such that $w$ has $01^{z}$ as a suffix.  Then $w0$ has $01^{z}0$ as a suffix so $z = c_{n} + i$ for some unique $n \geq 1$ and $0 \leq i < b_{n}$ by Lemma \ref{L:010}.  As $w0$ has $01^{c_{n}+i}0$ as a suffix, $w0$ shares a suffix with the word $1^{c_{n+1}}(\prod_{j=0}^{i-1}(B_{n}1^{c_{n}+j})^{a_{n}})(B_{n}1^{c_{n}+i})^{q}0$ for some $1 \leq q \leq a_{n}$.

First consider the case when $i > 0$.  If $w$ is a suffix of $1^{c_{n}+i-1}(B_{n}1^{c_{n}+i})^{a_{n}}$ then it is of form $(\ref{Pt-3})$ so we need only consider $w$ that have $01^{c_{n}+i-1}(B_{n}1^{c_{n}+i})^{q}$ as a suffix.  For such $w$, the word $w1$ has the suffix $01^{c_{n}+i-1}(B_{n}1^{c_{n}+i})^{q-1}B_{n}1^{c_{n}+i+1}$ but that word is only in $\mathcal{L}(T)$ if $q-1 = a_{n}$ which is impossible.

Now consider the case when $i = 0$, i.e.~$z = c_{n}$.  If $w$ is a suffix of $1^{c_{n}-1}(B_{n}1^{c_{n}})^{a_{n}}$ then it is of form $(\ref{Pt-3})$ so we may assume that $w$ has $1^{c_{n}-1}(B_{n}1^{c_{n}})^{q}$ as a strict suffix for some $1 \leq q \leq a_{n}$.  Since $B_{n}1^{c_{n}}$ is always preceded by $1^{c_{n}}$ (possibly as part of some $1^{c_{n+1}+i}$ or $1^{c_{n}+i}$), $w$ cannot have $01^{c_{n}-1}B_{n}1^{c_{n}}$ as a subword so $w$ has $1^{c_{n}}(B_{n}1^{c_{n}})^{q}$ as a suffix for some $1 \leq q \leq a_{n}$.

Take $q$ maximal so that $w$ has $1^{c_{n}}(B_{n}1^{c_{n}})^{q}$ as a suffix.

Consider first when $w$ has $1^{c_{n}}(B_{n}1^{c_{n}})^{a_{n}}$ as a suffix, i.e.~when $q = a_{n}$.  If $w = 1^{c_{n}}(B_{n}1^{c_{n}})^{a_{n}}$ then it is of form $(\ref{Pt-weird})$.  If $w$ has $01^{c_{n}}(B_{n}1^{c_{n}})^{a_{n}}$ as a suffix then $w0 \notin \mathcal{L}(T)$ as $0(1^{c_{n}}B_{n})^{a_{n}}1^{c_{n}}0 \notin \mathcal{L}(T)$.  If $w$ has $11^{c_{n}}(B_{n}1^{c_{n}})^{a_{n}}$ as a suffix then $w1$ has $1^{c_{n}+1}(B_{n}1^{c_{n}})^{a_{n}-1}B_{n}1^{c_{n}+1}$ as a suffix but that is not in $\mathcal{L}(T)$.

So we may assume $q < a_{n}$.  Since $1^{c_{n}}(B_{n}1^{c_{n}})^{q}$ is then of form $(\ref{Pt-3})$, we may assume $1^{c_{n}}(B_{n}1^{c_{n}})^{q}$ is a strict suffix of $w$.

Consider when $w$ has $01^{c_{n}}(B_{n}1^{c_{n}})^{q}$ as a suffix.  As $01^{c_{n}}(B_{n}1^{c_{n}})^{q}$ only appears as a suffix of $B_{n}1^{c_{n}}(B_{n}1^{c_{n}})^{q}$ and that word is always preceded by $1^{c_{n}}$ (possibly as part of some $1^{c_{n+1}+i}$), $w$ then shares a suffix with $1^{c_{n}}(B_{n}1^{c_{n}})^{q+1}$.  As $q$ is maximal, then $w$ is a suffix of $1^{c_{n}-1}(B_{n}1^{c_{n}})^{q+1}$ and, as $q < a_{n}$, this means $w$ is of form $(\ref{Pt-3})$.

We are left with the case when $w$ has $1^{c_{n}+1}(B_{n}1^{c_{n}})^{q}$ as a suffix for some $1 \leq q < a_{n}$.  If $q \geq 2$ then $w1$ has $1^{c_{n}+1}(B_{n}1^{c_{n}})^{q-1}B_{n}1^{c_{n}+1}$ as a suffix but that is not in $\mathcal{L}(T)$ for $q - 1 \geq 1$.  So we are left with the situation when $w$ shares a suffix with $1^{c_{n}+1}B_{n}1^{c_{n}}$.  So $w0$ shares a suffix with $1^{c_{n}+1} B_{n} 1^{c_{n}}0$ which must share a suffix with $1^{c_{n+1}}B_{n}1^{c_{n}}0$, meaning that $w$ shares a suffix with
$1^{c_{n+1}}B_{n}1^{c_{n}}$.  If $w$ is a suffix of $1^{c_{n}+b_{n}-1}B_{n}1^{c_{n}}$ then it is of form $(\ref{Pt-2})$.  If not then $w$ has the suffix
$1^{c_{n}+b_{n}}B_{n}1^{c_{n}}$ so $w1$ has suffix $1^{c_{n}+b_{n}}B_{n}1^{c_{n}+1}$ which is not in $\mathcal{L}(T)$ since $B_{n}1^{c_{n}+1}$ is always preceded by $B_{n}1^{c_{n}}$ or $B_{n}1^{c_{n}+1}$.
\end{proof}

\subsection{The \texorpdfstring{level-$n$}{level-n} complexity functions}

\begin{definition}
For a word $w$, define the \textbf{tail length} $z(w)$ such that
$w = u01^{z(w)}$
for some (possibly empty) word $u$ with the conventions that $z(1^{\len{w}}) = \infty$ and $z(u0) = 0$.
\end{definition}

\begin{definition}
For $1 \leq n < \infty$, the set of \textbf{level-$n$ generating words} is
\begin{align*}
W_{n} &= \{ w \in \mathcal{L}^{RS}(T) : c_{n} \leq z(w) < c_{n+1} \} 
\end{align*}
\end{definition}

\begin{proposition}\label{P:union}
$\mathcal{L}^{RS}(T) = \{ 1^{\ell} : \ell \in \mathbb{N} \} \sqcup \bigsqcup_{n=1}^{\infty} W_{n}$.
\end{proposition}
\begin{proof}
$\{ c_{n} \}$ is strictly increasing so the $W_{n}$ are disjoint.  Lemma \ref{2ss} (\ref{2s1}) says $1^{\ell} \in \mathcal{L}^{RS}(T)$ for all $\ell$ and as every word in $W_{n}$ has $0$ as a subword, these are disjoint from the $W_{n}$.  If $z(w) < c_{1}$ then $w0 \notin \mathcal{L}(T)$ by Lemma \ref{L:010} so all right-special words with $0$ as a subword are in some $W_{n}$.
\end{proof}

\begin{definition}
The \textbf{level-$n$ complexity} is
$
p_{n}(q) = |\{ w \in W_{n} : \len{w} < q \}|
$.
\end{definition}

By definition, $p_{n}(\ell+1) - p_{n}(\ell) = |\{ w \in W_{n} : \len{w} = \ell \}|$.

\begin{proposition}\label{P:p}
The complexity function $p$ satisfies
$p(q) = 1 + q + \sum_{n=1}^{\infty} p_{n}(q)$.
\end{proposition}
\begin{proof}
Using Proposition \ref{P:union} and that $p(\ell+1)-p(\ell) = |\{w\in\mathcal{L}^{RS}: \len{w}=\ell\}|$,
\begin{align*}
p(q) - p(1) &= \sum_{\ell=1}^{q-1} (p(\ell+1) - p(\ell)) = \sum_{\ell=1}^{q-1} |\{ w \in \mathcal{L}^{RS}(T) : \len{w} = \ell \}| \\
&= \sum_{\ell=1}^{q-1} \Big{(} \sum_{n=1}^{\infty} |\{ w \in W_{n} : \len{w} = \ell \}| + |\{ 1^{\ell} \}|\Big{)} 
= \sum_{\ell=1}^{q-1} \Big{(}\sum_{n=1}^{\infty} (p_{n}(\ell+1) - p_{n}(\ell)) + 1 \Big{)} \\
&= \sum_{n=1}^{\infty} \Big{(}\sum_{\ell=1}^{q-1} (p_{n}(\ell+1) - p_{n}(\ell))\Big{)} + q-1 
= \sum_{n=1}^{\infty} (p_{n}(q)-p_{n}(1)) + q - 1
\end{align*}
All words in $W_{n}$ have length at least $1 + c_{n} > 1$ so $p_{n}(1) = 0$.  The claim follows as $p(1) = 2$.
\end{proof}

\subsection{Counting quasi-staircase words}

\begin{lemma}\label{L:2s}
	If $w \in W_{n}$ then exactly one of the following holds:
\begin{enumerate}[\hspace{10pt}(i)\hspace{5pt}]
\item\label{Lt-2} $w$ is a suffix of $1^{c_{n}+i-1}\big{(}B_{n}1^{c_{n}+i}\big{)}^{a_{n}}$ and $\len{w} > c_{n}+i$ for some $0 \leq i < b_{n}$; 
\item\label{Lt-1} $w$ is a suffix of $1^{c_{n}+b_{n}-1}B_{n}1^{c_{n}}$ and $\len{w} > h_{n} + 2c_{n}$; or
\item\label{Lt-weird} $w = 1^{c_{n}}(B_{n}1^{c_{n}})^{a_{n}}$
\end{enumerate}
\end{lemma}
\begin{proof}
	The only words in Proposition \ref{P:twosuccs} which have $c_{n} \leq z(w) < c_{n+1}$ are of the stated forms; Lemma \ref{2ss} (\ref{2s0}), (\ref{2si}) and (\ref{2sspecial}) state that these words are in $\mathcal{L}^{RS}(T)$.  The forms do not overlap due to the restriction on $\len{w}$ in form $(\ref{Lt-1})$.
\end{proof}

\begin{lemma}\label{Mi}
Fix $0 \leq i < b_{n}$.  For $c_{n} + i < \ell < a_{n}h_{n} + (a_{n}+1)(c_{n} + i)$ there is exactly one word in $W_{n}$ of form $(\ref{Lt-2})$ for that value of $i$; for $\ell$ not in that range, there are no words in $W_{n}$ of form $(\ref{Lt-2})$ for that $i$.
\end{lemma}
\begin{proof}
For $w \in W_{n}$ of form $(\ref{Lt-2})$, $w = u1^{c_{n}+i}$ where $u$ is a nonempty suffix of $1^{c_{n}+i-1}(B_{n}1^{c_{n}+i})^{a_{n}-1}B_{n}$.  The word $u$ is unique if it exists which is exactly when $c_{n}+i = \len{1^{c_{n}+i}} < \len{w} \leq \len{1^{c_{n}+i-1}(B_{n}1^{c_{n}+i})^{a_{n}}} = a_{n}h_{n} + (a_{n}+1)(c_{n}+i) - 1$.
\end{proof}

\begin{lemma}\label{M0}
For $h_{n} + 2c_{n} < \ell < h_{n} + 2c_{n} + b_{n}$ there is exactly one word in $W_{n}$ of form $(\ref{Lt-1})$; for $\ell$ not in that range, there are no words in $W_{n}$ of form $(\ref{Lt-1})$.
\end{lemma}
\begin{proof}
To be of that form, $w = u1^{c_{n}}$ where $u$ is a nonempty suffix of $1^{c_{n}+b_{n}-1}B_{n}$ that has $1^{c_{n}+1}$ as a prefix.  The word $u$ is unique if it exists and it exists exactly when $h_{n} + 2c_{n} + 1 = \len{1^{c_{n}+1}B_{n}1^{c_{n}}} \leq \len{w} \leq \len{1^{c_{n}+b_{n}-1}B_{n}1^{c_{n}}} = h_{n} + 2c_{n} + b_{n} - 1$.
\end{proof}

\begin{lemma}\label{cf0}
	If $\ell \leq c_{n}$ then $p_{n}(\ell+1) - p_{n}(\ell) = 0$.
\end{lemma}
\begin{proof}
	Every $w \in W_{n}$ has subwords $1^{c_{n}}$ and $0$ so $\len{w} \geq c_{n} + 1$ therefore $p_{n+1}(\ell) = p_{n}(\ell) = 0$.
\end{proof}

\begin{lemma}\label{cf1}
	If $c_n < \ell < c_n + b_n$ then $p_{n}(\ell +1)-p_{n}(\ell)=\ell-c_n$.
\end{lemma}

\begin{proof}
	Lemma \ref{Mi} applies for $0 \leq i < \ell- c_{n}$ but not for $\ell - c_{n} \leq i < b_{n}$.  Lemma \ref{M0} does not apply.
\end{proof}

\begin{lemma} \label{cf2}
	If $c_n+b_n \leq \ell \leq h_n+2c_n$ then $p_{n}(\ell+1)-p_{n}(\ell)=b_n$.
\end{lemma}

\begin{proof}
	Lemma \ref{Mi} applies for all $0 \leq i < b_{n}$ and Lemma \ref{M0} does not apply.
\end{proof}

\begin{lemma} \label{cf3}
	If $h_n+2c_n < \ell < h_{n} + 2c_{n} + b_{n}$ then $p_{n}(\ell+1)-p_{n}(\ell)=b_n + 1$.
\end{lemma}

\begin{proof}
	Lemma \ref{Mi} applies for all $0 \leq i < b_{n}$ and Lemma \ref{M0} applies.
\end{proof}

\begin{lemma} \label{cf4}
	If $h_n+2c_n+b_{n} \leq \ell < a_{n}h_n+(a_{n}+1)c_n$ then $p_{n}(\ell+1)-p_{n}(\ell)=b_n$.
\end{lemma}

\begin{proof}
	Lemma \ref{Mi} applies for all $0 \leq i < b_{n}$ and Lemma \ref{M0} does not apply.
\end{proof}

\begin{lemma} \label{cf5weird}
	$p_{n}(a_{n}h_{n} + (a_{n}+1)c_{n} + 1) - p_{n}(a_{n}h_{n} + (a_{n}+1)c_{n}) = b_{n} + 1$.
\end{lemma}
\begin{proof}
	Lemma \ref{Mi} applies for all $0 \leq i < b_{n}$ and Lemma \ref{M0} does not apply.  Lemma \ref{L:2s} form $(\ref{Lt-weird})$ gives one additional word in $W_{n}$.
\end{proof}

\begin{lemma}\label{cf5}
	If $a_{n}h_n+(a_{n}+1)c_n + 1 < \ell < a_{n}h_{n} + (a_{n}+1)(c_{n} + b_{n}-1)$ then $p_{n}(\ell+1) - p_{n}(\ell) \leq b_n$.
\end{lemma}
\begin{proof}
	Lemma \ref{Mi} applies for some subset of $0 \leq i < b_{n}$ and Lemma \ref{M0} does not apply.
\end{proof}

\begin{lemma}\label{cf5prime}
	$p(a_{n}h_{n} + (a_{n}+1)c_{n} + (a_{n}+1)(b_{n}-1)) - p(a_{n}h_{n} + (a_{n}+1)c_{n}) = \frac{1}{2}(a_{n}+1)b_{n}(b_{n}-1) + 1$.
\end{lemma}
\begin{proof}
	For each $0 \leq i < b_{n}$, Lemma \ref{Mi} applies for $\ell = a_{n}h_{n} + (a_{n}+1)c_{n} + y$ exactly when $0 \leq y < (a_{n}+1)i$, therefore there are a total of $(a_{n}+1)\frac{1}{2}b_{n}(b_{n}-1)$ words in $W_{n}$ of the enclosed lengths from Lemma \ref{Mi}.  Lemma \ref{M0} does not apply and Lemma \ref{L:2s} form $(\ref{Lt-weird})$ gives one additional word.
\end{proof}

\begin{lemma}\label{cf6}
	If $a_{n}h_{n} + (a_{n}+1)(c_{n}+b_{n}-1) \leq \ell$ then $p_{n}(\ell+1)-p_{n}(\ell)=0$.
\end{lemma}
\begin{proof}
	 Neither Lemma \ref{Mi} nor \ref{M0} apply.
\end{proof}

\subsection{Bounding the complexity of quasi-staircases}

Since $p_{n}(\ell+1) - p_{n}(\ell) = 0$ for $\ell \geq  a_{n}h_{n} + (a_{n}+1)(c_{n} + b_{n} - 1)$, we define:
\begin{definition}
The \textbf{post-productive sequence} is
\[
m_{n} = a_{n}h_{n} + (a_{n}+1)(c_{n} + b_{n} - 1) \quad\quad m_{0} = 0
\]
\end{definition}

\begin{lemma}\label{pmn}
$p_{n}(m_{n}) = h_{n+1} - h_{n}$
\end{lemma}
\begin{proof}
By Lemma \ref{cf0}, $p_{n}(c_{n}) = \sum_{\ell=0}^{c_{n}-1} (p_{n}(\ell+1) - p_{n}(\ell)) = 0$.

By Lemma \ref{cf1},
$
p_{n}(c_{n} + b_{n}) - p_{n}(c_{n}) = \sum_{\ell=c_{n}}^{c_{n} + b_{n}-1} (\ell - c_{n}) = \frac{1}{2}b_{n}(b_{n} - 1)
$.

By Lemma \ref{cf2},
$
p_{n}(h_{n} + 2c_{n}+1) - p_{n}(c_{n} + b_{n}) = (h_{n} + c_{n} + 1 - b_{n})b_{n}
$.

By Lemma \ref{cf3},
$
p_{n}(h_{n} + 2c_{n} + b_{n}) - p_{n}(h_{n} + 2c_{n}+1) = (b_{n}+1)(b_{n}-1)
$.

By Lemma \ref{cf4},
$
p_{n}(a_{n}h_{n} + (a_{n}+1)c_{n}) - p_{n}(h_{n} + 2c_{n} + b_{n}) = ((a_{n}-1)h_{n} + (a_{n}-1)c_{n} - b_{n})b_{n}
$.

By Lemma \ref{cf5prime},
$
p_{n}(m_{n}) - p(a_{n}h_n+(a_{n}+1)c_n) = \frac{1}{2}(a_{n}+1)b_{n}(b_{n}-1) + 1
$.
Therefore
\begin{align*}
p_{n}(m_{n}) &= 
\frac{1}{2}b_{n}(b_{n}-1) + (h_{n}+c_{n}+1-b_{n})b_{n} + (b_{n}+1)(b_{n}-1) \\
&\quad\quad + ((a_{n}-1)h_{n} + (a_{n}-1)c_{n} - b_{n})b_{n} + \frac{1}{2}(a_{n}+1)b_{n}(b_{n}-1) + 1 \\
&= a_{n}b_{n}h_{n} + a_{n}b_{n}c_{n} + \frac{1}{2}a_{n}b_{n}(b_{n}-1) + b_{n}(b_{n}-1) + b_{n} - b_{n}^{2} + b_{n}^{2} - 1 - b_{n}^{2} + 1 
= h_{n+1} - h_{n}
\qedhere
\end{align*}
\end{proof}

\begin{definition}
For $q \in \mathbb{N}$ define
\[
\rho(q) = \max \{ n : m_{n} \leq q \}  \quad\quad\quad\text{and}\quad\quad\quad \beta(q) = \min \{ n : q < c_{n+1} \}
\]
\end{definition}

\begin{lemma} $\rho(q) \leq \beta(q)$
 \end{lemma}
\begin{proof}
If  $\beta(q) \leq \rho(q) - 1$ then $m_{\rho(q)} \leq q < c_{\beta(q)+1} \leq c_{\rho(q)-1+1} = c_{\rho(q)} < m_{\rho(q)}$ is impossible.
\end{proof}

\begin{lemma}\label{L:aa}
If $q < c_{n}$ then $p_{n}(q) = 0$.
If $c_{n} \leq q < m_{n}$ then
$
p_{n}(q) \leq (q - c_{n} + 1)b_{n}
$.
If $m_{n} \leq q$ then $p_{n}(q) = h_{n+1} - h_{n}$.
\end{lemma}
\begin{proof}
Lemma \ref{cf0} gives $p_{n}(\ell+1) - p_{n}(\ell) = 0$ for $0 \leq \ell < c_{n}$.  Lemmas \ref{cf1}, \ref{cf2}, \ref{cf3}, \ref{cf4} and \ref{cf5} all give $p_{n}(\ell+1) - p_{n}(\ell) \leq b_{n}$ for $c_{n} \leq \ell < m_{n}$ except for Lemma \ref{cf3} which gives $p_{n}(\ell+1) - p_{n}(\ell) = b_{n} + 1$ for exactly $b_{n}-1$ values of $\ell$ and Lemma \ref{cf5weird} which gives one additional word.  Then, for $c_{n} \leq q < m_{n}$,
\begin{align*}
p_{n}(q) &= \sum_{\ell=0}^{q-1} (p_{n}(\ell+1) - p_{n}(\ell)) 
= \sum_{\ell=0}^{c_{n}-1} 0 + \sum_{\ell=c_{n}}^{q-1} (p_{n}(\ell+1) - p_{n}(\ell)) 
\leq (q-c_{n})b_{n} + b_{n}
\end{align*}
Lemma \ref{cf6} says $p_{n}(\ell+1) - p_{n}(\ell) = 0$ for $\ell \geq m_{n}$ so when $q \geq m_{n}$, $p_{n}(q) = p_{n}(m_{n})$ and Lemma \ref{pmn} gives the final statement.
\end{proof}

\begin{proposition}\label{P:pbounds}
$
p(q) \leq q\Big{(} 2 + \sum_{n=\rho(q)}^{\beta(q)}b_{n}\Big{)}
$ for all $q$.
\end{proposition}
\begin{proof}
For $n$ such that $\beta(q) < n$, by Lemma \ref{cf0}, $p_{n}(q) = 0$.  Proposition \ref{P:p} and Lemma \ref{L:aa} give, using that $h_{1} = 1$ so $1 + \sum_{n=1}^{\rho(q)}(h_{n+1}-h_{n}) = h_{\rho(q)+1}$,
\begin{align*}
p(q) &= q + 1 + \sum_{n=1}^{\rho(q)} p_{n}(q) + \sum_{n=\rho(q)+1}^{\beta(q)} p_{n}(q) + \sum_{n=\beta(q)+1}^{\infty} p_{n}(q) \\
&\leq q + 1 + \sum_{n=1}^{\rho(q)} (h_{n+1} - h_{n}) + \sum_{n=\rho(q)+1}^{\beta(q)} (q - c_{n} + 1)b_{n} + 0 
\leq q + h_{\rho(q) + 1} + \sum_{n=\rho(q)+1}^{\beta(q)} qb_{n}
\end{align*}
\begin{flalign*}
\text{and} && h_{\rho(q) + 1} &= h_{\rho(q)} + b_{\rho(q)}(a_{\rho(q)}h_{\rho(q)} + a_{\rho(q)}c_{\rho(q)} + \frac{1}{2}a_{\rho(q)}(b_{\rho(q)}-1)) && \\
&&&\leq h_{\rho(q)} + b_{\rho(q)}m_{\rho(q)}
\leq m_{\rho(q)}(1 + b_{\rho(q)})
\leq q(1 + b_{\rho(q)}) && \qedhere
\end{flalign*}
\end{proof}

\section{Quasi-staircase complexity arbitrarily close to linear}

\begin{lemma}\label{P:Ts}
Let $\{ d_{n} \}$ be a nondecreasing sequence of integers such that $d_{n} \to \infty$ and $d_{1} = d_{2} = 1$ and $d_{n+1} - d_{n} \in \{ 0, 1 \}$ and $d_{n+1} - d_{n}$ does not take the value $1$ for consecutive $n$.

Let $\{ b_{n} \}$ be a nondecreasing sequence of integers such that $b_{n} \to \infty$ and $b_{1} = 3$ and $b_{n} \leq n + 2$.

Set $a_{n} = 2n^{2}+2$.  Set $c_{1} = 1$ and for $n > 1$,
\begin{align*}
c_{n} &= \left\{ \begin{array}{ll} m_{n-d_{n}} &\quad\quad\text{when $d_{n} = d_{n-1}$} \\
c_{n-1} + b_{n-1} &\quad\quad\text{when $d_{n} = d_{n-1}+1$} \end{array} \right.
\end{align*}
Then $\{a_{n}\}, \{b_{n}\}, \{c_{n}\}$ define a quasi-staircase such that $\sum \frac{a_{n}b_{n}^{2}+a_{n+1}b_{n+1}+c_{n+1}}{h_{n}} < \infty$ and $\sum \frac{1}{a_{n}b_{n}} < \infty$.
\end{lemma}
\begin{proof}
Since $r_{n} = a_{n}b_{n}$, we have $6n^{2} + 6 \leq r_{n} \leq (2n^{2}+2)(n+2)$.  Then $\prod_{j=1}^{n-1}(r_{j}+1) \geq n!$ so $h_{n} \geq \prod_{j=1}^{n-1} (r_{j}+1) \geq n!$ so
$\sum \frac{a_{n}b_{n}^{2}+a_{n+1}b_{n+1}}{h_{n}} \leq \sum \frac{(2n^{2}+2)(n+2)^{2} + (2(n+1)^{2}+2)(n+3)}{n!} < \infty$.

Since $a_{n} + 1 < \frac{3}{2}a_{n} \leq \frac{1}{2}b_{n}a_{n}$, then $m_{n} < h_{n+1}$ for all $n$.  

For $n$ such that $c_{n+1} = c_{n} + b_{n}$, since $d_{n+1} = d_{n} + 1$, also $d_{n} = d_{n-1}$.  So for sufficiently large such $n$,
\[
\frac{c_{n+1}}{h_{n}} = \frac{c_{n} + b_{n}}{h_{n}} = \frac{m_{n-d_{n}} + b_{n}}{h_{n}} < \frac{2h_{n-d_{n}+1}}{h_{n}}
< \frac{2}{\prod_{i=1}^{d_{n}-1} r_{n-d_{n}+i}} < \frac{2}{r_{n-1}} < \frac{2}{2(n-1)^{2}+2}
\]
and for sufficiently large $n$ such that $c_{n+1} = m_{n+1-d_{n}}$,
\[
\frac{c_{n+1}}{h_{n}} = \frac{m_{n+1-d_{n}}}{h_{n}} < \frac{h_{n-d_{n}+2}}{h_{n}} 
< \frac{1}{r_{n-1}} < \frac{1}{2(n-1)^{2}+2}
\]
Therefore, as $\sum \frac{2}{2(n-1)^{2} + 2} < \infty$, it follows that $\sum \frac{c_{n+1}}{h_{n}} < \infty$ and the result follows.
\end{proof}

\begin{lemma}\label{L:assump}
If $f : \mathbb{N} \to \mathbb{N}$ is any function such that $f(q) \to \infty$ then there exists $g : \mathbb{N} \to \mathbb{N}$ which is nondecreasing such that $g(1)=1$ and $g(q) \leq f(q)$ and $g(q+2) - g(q) \leq 1$ for all $q$ and $g(q) \to \infty$.
\end{lemma}
\begin{proof}
Set $f^{*}(q) = \inf_{q^{\prime} \geq q} f(q^{\prime})$.  Then $f^{*}(q) \to \infty$ and $f^{*}(q)$ is nondecreasing and $f^{*}(q) \leq f(q)$ for all $q$.
Set $g(1) = 1 \leq f^{*}(1)$.  For $n \geq 0$, set $g(2n+2) = g(2n+1)$ and for $n \geq 1$ set
\[
g(2n+1) = g(2n) + \left\{ \begin{array}{ll} 1 \quad\quad &\text{when $f^{*}(2n+1) > f^{*}(2n-1)$} \\ 0 \quad\quad &\text{otherwise} \end{array} \right.
\]
Then $g$ is nondecreasing and $g(q+2) - g(q) \leq 1$ for all $q$.  Since $f^{*}$ is integer-valued, if $f^{*}(2n+1)-f^{*}(2n-1) \ne 0$ then $f^{*}(2n+1) - f^{*}(2n-1) \geq 1$.  Then $g(2n+1)-g(2n-1) \leq f^{*}(2n+1) - f^{*}(2n-1)$ so for all $n$ we have
\[
g(2n+1) = g(1) + \sum_{m=1}^{n} (g(2m+1) - g(2m-1))
\leq f^{*}(1) + \sum_{m=1}^{n} (f^{*}(2m+1) - f^{*}(2m-1))
= f^{*}(2n+1)
\]
so, as $g(2n+2) = g(2n+1) \leq f^{*}(2n+1) \leq f^{*}(2n+2)$, we have $g(q) \leq f^{*}(q) \leq f(q)$ for all $q$.  If $g(q) \leq C$ for all $q$ then $f^{*}(2n+1) = f^{*}(2n-1)$ eventually, contradicting that $f^{*}(q) \to \infty$.  Therefore $g(q) \to \infty$.
\end{proof}

\begin{theorem}\label{T:attainq}
Let $f : \mathbb{N} \to \mathbb{N}$ be any function such that $f(q) \to \infty$.  There exists a quasi-staircase transformation with $\sum \frac{a_{n}b_{n}^{2}+a_{n+1}b_{n+1}+c_{n+1}}{h_{n}} < \infty$, $\sum \frac{1}{a_{n}b_{n}} < \infty$, $\frac{b_{n}}{a_{n}} \to 0$ and complexity satisfying $\frac{p(q)}{qf(q)} \to 0$.
\end{theorem}
\begin{proof}
By Lemma \ref{L:assump}, we may assume $f$ is nondecreasing and that $f(n+2) - f(n) \leq 1$ for all $n$.
Then $f(n+1) - f(n) \in \{ 0, 1 \}$ and is never $1$ for two consecutive values.  We may also assume $f(1) = 1$.

Set $d_{1} = d_{2} = 1$ and $d_{n} = \floor*{\sqrt[3]{f(n)}}$ for $n > 2$.  Then $d_{n} \to \infty$ is nondecreasing.  Also $d_{n+1} - d_{n} \in \{ 0, 1 \}$ and is never $1$ for two consecutive values.

Set $b_{n} = 3$ for all $n$ such that $\sqrt[3]{f(n)} < 3$ and $b_{n} = \floor*{\sqrt[3]{f(n)}}$ for $n$ such that $\sqrt[3]{f(n)} \geq 3$.  Then $b_{n} \to \infty$ is nondecreasing and $b_{n} \leq f(n) + 2 \leq n+2$ as $f(n) \leq n$ since $f(1) = 1$ and $f(n+2) - f(n) \leq 1$ imply $f(n) \leq 1 + \frac{n}{2}$.

Take the quasi-staircase transformation from Lemma \ref{P:Ts} with defining sequences $\{ a_{n} \}$ and $\{ c_{n} \}$.
As $a_{n} = 2n^{2} + 2$ and $b_{n} = \max(3,\sqrt[3]{f(n)}) \leq \sqrt[3]{n}$, we have $\frac{b_{n}}{a_{n}} \to 0$.

Since $0 \leq d_{n+1} - d_{n} \leq 1$, the sequence $n - d_{n}$ is nondecreasing and attains every value in $\mathbb{N}$.
For each $q$, let $n_{q}$ be the largest $n$ such that $m_{n - d_{n}} \leq q$.  Then $q < m_{n_{q}+1-d_{n_{q}+1}}$ so $n_{q} + 1 - d_{n_{q}+1} > n_{q} - d_{n_{q}}$ and so $1 > d_{n_{q}+1} - d_{n_{q}}$ meaning that $d_{n_{q}+1} = d_{n_{q}}$.  Therefore $c_{n_{q}+1} = m_{n_{q}+1 - d_{n_{q}+1}} = m_{n_{q}-d_{n_{q}} + 1}$.

So $\rho(q) = n_{q} - d_{n_{q}}$ as $m_{n_{q}-d_{n_{q}}} \leq q < m_{n_{q} + 1 -d_{n_{q} + 1}} = m_{n_{q} - d_{n_{q}} + 1}$ and $\beta(q) \leq n_{q}$ since $q < m_{n_{q} - d_{n_{q}} + 1} = c_{n_{q}+1}$.
By Proposition \ref{P:pbounds}, since $q \geq n_{q}$ and $f$ is nondecreasing to infinity and $n_{q} \to \infty$,
\begin{align*}
\frac{p(q)}{qf(q)} &\leq \frac{2 + \sum_{n=\rho(q)}^{\beta(q)} b_{n}}{f(q)} 
\leq \frac{2 + \sum_{n=n_{q} - d_{n_{q}}}^{n_{q}} b_{n}}{f(q)}
\leq \frac{2 + (d_{n_{q}}+1)b_{n_{q}}}{f(n_{q})} \\
&\leq \frac{2 + (\sqrt[3]{f(n_{q})}+1)\sqrt[3]{f(n_{q})}}{f(n_{q})} 
= \frac{2}{f(n_{q})} + \frac{1}{\sqrt[3]{f(n_{q})}} + \frac{1}{(\sqrt[3]{f(n_{q})})^{2}} \to 0 \qedhere
\end{align*}
\end{proof}

\section{Mixing for quasi-staircase transformations}\label{S:mix}

The goal of this section is to prove the following.

\begin{theorem}\label{C:mixing}
	Let $T$ be a quasi-staircase transformation such that $\sum \frac{a_{n}b_{n}+b_{n+1}+c_{n+1}}{h_{n}} < \infty$ and $\sum \frac{1}{a_{n}b_{n}} < \infty$ and $\frac{a_{n}b_{n}^{2}}{h_{n}} \to 0$ and $\frac{a_{n+1}b_{n+1}}{h_{n}} \to 0$ and $\frac{b_{n}}{a_{n}} \to 0$.  Then $T$ is mixing.
\end{theorem}

Throughout this section, we assume that all transformations are on probability spaces.  Recall that $b_{n} \to \infty$ by definition for quasi-staircase transformations.

%Recall that 
%a sequence $\{ t_{n} \}$ is \textbf{mixing} when for all measurable sets $A$ and $B$,
%$
%\mu(T^{t_{n}}A \cap B) \to \mu(A)\mu(B).
%$

We first introduce some notation.

\begin{notation}
For measurable sets $A$ and $B$, write
\[
\muB(A) = \mu(A \cap B) - \mu(A)\mu(B)
\]
\end{notation}

So $\{ t_{n} \}$ is mixing when $\muB(T^{t_{n}}A) \to 0$ for all measurable $A$ and $B$.  The following is left to the reader:
\begin{lemma}
If $A$ and $A^{\prime}$ are disjoint then
\[
\muB(A \sqcup A^{\prime}) = \muB(A) + \muB(A^{\prime})
\quad\quad\text{and}\quad\quad
|\muB(A)| \leq \mu(A)
\]
and, writing $\chi_{B}(x) = \bbone_{B}(x) - \mu(B)$, for $n \in \mathbb{Z}$,
$
\muB(T^{n}A) = \int_{A} \chi_{B} \circ T^{n}d\mu.
$
\end{lemma}

For a rank-one transformation $T$, a sequence $\{ t_{n} \}$ is \textbf{rank-one uniform mixing} when for every union of levels $B$,
$
\sum_{j=0}^{h_{n}-1} |\muB(T^{t_{n}}I_{n,j})| \to 0.
$
Rank-one uniform mixing for a sequence implies mixing for that sequence \cite{CreutzSilva2004} Proposition 5.6.

\begin{notation}
For $h_{n} \leq j < h_{n} + c_{n}$, let $I_{n,j} = T^{j-h_{n}+1}I_{n,h_{n}-1}$ be the union of the $(j-h_{n})^{th}$ stage of the $c_{n}$ spacer levels added above every subcolumn.
Write $$\tilde{h}_{n} = h_{n} + c_{n}$$
\end{notation}

\begin{lemma}\label{L:mixtrick}
Let $T$ be a rank-one transformation, $B$ a union of levels in some column $C_{N}$ and $n \geq N$. Then for any $0 \leq j < \tilde{h}_{n}$ and $0 \leq i \leq r_{n}$,
\[
\muB(I_{n,j}^{[i]}) = \frac{1}{r_{n}+1}\muB(I_{n,j})
\]
\end{lemma}
\begin{proof}
Since $B$ is a union of levels in $C_{N}$, either $I_{n,j} \subseteq B$ or $I_{n,j} \cap B = \emptyset$.  If $I_{n,j} \subseteq B$ then $\mu(I_{n,j}^{[i]} \cap B) = \mu(I_{n,j}^{[i]}) = \frac{1}{r_{n}+1}\mu(I_{n,j}) = \frac{1}{r_{n+1}}\mu(I_{n,j} \cap B)$ and if $I_{n,j} \cap B = \emptyset$ then $\mu(I_{n,j}^{[i]} \cap B) = 0 = \frac{1}{r_{n}+1}\mu(I_{n,j} \cap B)$.
\end{proof}

\begin{proposition}\label{P:finmeas}
	Let $T$ be a quasi-staircase transformation given by $\{ a_{n} \}$, $\{ b_{n} \}$ and $\{ c_{n} \}$ with height sequence $\{ h_{n} \}$.
	Then $T$ is on a finite measure space if and only if $\sum \frac{c_{n} + b_{n}}{h_{n}} < \infty$.
\end{proposition}
\begin{proof}
	Writing $S_{n}$ for the spacers added above the $n^{th}$ column,
	\[
	\mu(S_{n}) = (c_{n}r_n+\frac{1}{2}r_{n}(b_{n}-1))\mu(I_{n+1}) = \left(c_n\frac{r_n}{r_{n}+1}+\frac{1}{2}\frac{r_{n}(b_{n}-1)}{r_{n}+1}\right)\mu(I_{n}) \leq \frac{c_{n} + b_{n}}{h_{n}}~\mu(C_{n})
	\]
	and therefore
	$
	\mu(C_{n+1}) = \mu(C_{n}) + \mu(S_{n}) \leq \big{(} 1 + \frac{c_{n} + b_{n}}{h_{n}} \big{)} \mu(C_{n})
	$.
	Then
	$
	\mu(C_{n+1}) \leq \prod_{j=1}^{n} \big{(}1 + \frac{c_{j} + b_{j}}{h_{j}}\big{)}\mu(C_{1}),
	$
	so the claim follows from \cite{knopp} p.219 that $\prod_{j=1}^{\infty}(1 + \frac{c_{j}+b_{j}}{h_{j}}) < \infty$ if and only if $\sum_{j=1}^{\infty} \frac{c_{j}+b_{j}}{h_{j}} < \infty$.
%	meaning that
%	$
%	\log(\mu(C_{n+1})) \leq \log(\mu(C_{1})) + \sum_{j=1}^{n} \log(1 + \frac{c_{j} + b_{j}}{h_{j}}).
%	$
%	As $\frac{c_{n} + b_{n}}{h_{n}} \to 0$, since $\log(1+x) \approx x$ for $x \approx 0$,
%	$
%	\lim_{n} \log(\mu(C_{n+1})) \lesssim \log(\mu(C_{1})) + \sum_{j=1}^{\infty} \frac{c_{j} + b_{j}}{h_{j}} < \infty.
%	$
	\end{proof}

\subsection{Mixing along most sequences}

For clarity of exposition, we state the results which follow from now-standard techniques for proving mixing on staircases with explanations of how one could modify the corresponding proofs in the literature to our class of transformations.  Detailed proofs, including all necessary modifications, are deferred to the appendix.

\subsubsection{(Weak) power ergodicity}

\begin{proposition}\label{P:pe}
Let $T$ be a quasi-staircase transformation and $B$ a measurable set.  Then
\[
\max_{1 \leq k \leq n} \int\left|\frac{1}{n}\sum_{j=0}^{n-1}\chi_{B}\circ T^{-jk}\right|~d\mu \to 0
\]
\end{proposition}

The proof of Proposition \ref{P:pe} is essentially identical to the proof of (weak) power ergodicity for staircases except that one must replace the height sequence by $\{ h_{n} + c_{n} \}$ and then observe that for a union of levels $B$ in $C_{n}$, any fixed positive integer $k$, and a level $I$ in $C_{n}$ at least $kb_{n}$ above the base,
\[
\mu(T^{k(h_{n}+c_{n})}I \cap B) \approx \frac{1}{b_{n}}\sum_{i=0}^{b_{n}-1} \mu(T^{-ik}I \cap B) \pm \frac{2k}{r_{n}}\mu(I)
\]
which follows from the standard technique that $\mu(T^{-ik}I^{[i]} \cap B) = \frac{1}{r_{n}}\mu(T^{-ik}I \cap B)$ provided that $I$ is at least $ik$ levels above the base of the tower.  From there, deducing the proposition is identical, modulo the obvious replacement of $r_{n}$ by $b_{n}$ throughout the proof, to the proof of (weak) power ergodicity for elevated staircases in \cite{CPR}.

\subsubsection{Mixing between \texorpdfstring{$a_{n}\tilde{h}_{n}$}{a\_n h\_n} and \texorpdfstring{$\tilde{h}_{n+1}$}{h\_(n+1)}}

\begin{proposition}\label{mixingES} Let $T$ be a quasi-staircase transformation such that $\frac{a_n b_n^2}{h_n}\rightarrow 0$ and $\frac{b_{n}}{a_{n}} \to 0$ and $B$ be a union of levels in some fixed $C_{N}$.   For $n > N$, set
\[
M_{B,n} := \max_{a_{n}\tilde{h}_{n} \leq t < \tilde{h}_{n+1}} \quad \sum_{j=0}^{h_{n}-1} |\muB(T^{t}I_{n,j})|
\]
Then $\lim_{n \to \infty} M_{B,n} = 0$.
\end{proposition}

The proof of Proposition \ref{mixingES} follows from the standard argument for proving mixing on (elevated) staircases, see e.g.~\cite{CPR}.  Consider $T^{a_{n}\tilde{h}_{n}}$ applied to a level $I_{n,j}$ in the $n^{th}$ column.  Since every sublevel is pushed through at least one spacer, $T^{a_{n}\tilde{h}_{n}} I_{n,j}$ consists of $a_{n}$ sublevels in $T^{-i}I_{n,j}$ for each $1 \leq i < b_{n}$ plus $a_{n}$ sublevels which are pushed through the top of the next column.  Since $b_{n} \to \infty$, the convergence of the ergodic average $\frac{1}{b_{n}}\sum_{i=0}^{b_{n}-1}T^{-i}$ implies mixing along this sequence.

For the sequence $\{ k_{n}a_{n}\tilde{h}_{n} \}$ where $1 \leq k_{n} < b_{n}$, the resulting average $\frac{1}{b_{n}}\sum_{i=0}^{b_{n}-1}T^{-ik_{n}}$ converges by (weak) power ergodicity so these sequences are likewise mixing.  The general case of times between $a_{n}\tilde{h}_{n}$ and $\tilde{h}_{n+1}$ then follows from the standard interpolation argument and Blum-Hanson trick combined with the Block Lemma.

 \subsubsection{Mixing between \texorpdfstring{$\tilde{h}_{n}$}{h\_n} and \texorpdfstring{$b_{n}\tilde{h}_{n}$}{b\_n h\_n}}

 \begin{proposition}\label{mixingES2}
 Let $T$ be a quasi-staircase transformation with $\frac{b_{n}^{2}}{h_{n}} \to 0$ and $\frac{b_{n}}{a_{n}} \to 0$ and $B$ be a union of levels in some fixed $C_{N}$.  For $n > N$, set
 \[
\widehat{M}_{B,n} := \max_{\tilde{h}_{n}  \leq t < b_{n}\tilde{h}_{n}} \quad \sum_{j=0}^{h_{n}-1} |\muB(T^{t}I_{n,j})|
 \]
Then $ \lim_{n\to\infty} \widehat{M}_{B,n} = 0$.
 \end{proposition}
 
The proof of Proposition \ref{mixingES2} follows nearly immediately from the fact that the resulting ergodic average for $T^{k_{n}\tilde{h}_{n}}$ is already known to converge by weak power ergodicity.

\subsection{Mixing between \texorpdfstring{$b_{n}\tilde{h}_{n}$}{b\_n h\_n} and \texorpdfstring{$a_{n}\tilde{h}_{n}$}{a\_n h\_n}}
 
 The new techniques introduced here apply to the times not covered by the above results.  For $b_{n}\tilde{h}_{n} \leq t < a_{n}\tilde{h}_{n}$, we write $t$ uniquely as $k_{n}\tilde{h}_{n} + y_{n}$ for $b_{n} \leq k_{n} < a_{n}$ and $|y_{n}| \leq \frac{1}{2}\tilde{h}_{n}$ (taking $y_{n}$ positive in the case when $|y_{n}| = \frac{1}{2}\tilde{h}_{n}$).
 
\subsubsection{Mixing using the previous column (mixing when \texorpdfstring{$|y_{n}| \geq a_{n-1}\tilde{h}_{n-1}$}{|y\_n| > a\_(n-1)\tilde{h}\_(n-1)})}\label{prevcol}
 
The first new idea we introduce is showing mixing by invoking known mixing times for the previous tower ($y_{n}$ is already known to be a mixing time for the previous column in the case $|y_{n}| \geq a_{n-1}\tilde{h}_{n-1}$).

We first explain the argument, eliding many details, then present the detailed proofs.  Since $k_{n} < a_{n}$, application of $T^{k_{n}\tilde{h}_{n}}$ to a level effectively gives $b_{n}$ blocks of $a_{n}$ sublevels, each block having passed through the same number of spacers.  The dominant term in the number of spacers is $ik_{n}$ where $1 \leq i < b_{n}$ where each block corresponds to a single $i$.

Since each block will then also have $T^{y_{n}}$ applied to it, the mixing nature of $y_{n}$ can be used to show that each block is mixed.  The goal now is to prove the following.
 
 \begin{proposition}\label{midy}
Let $T$ be a quasi-staircase transformation such that $\frac{a_{n+1}b_{n+1}+c_{n+1}+a_{n}b_{n}^{2}}{h_{n}} \to 0$.  Let $B$ be a union of levels in some column $C_{N}$.  
For $n > N$, set
\[
\widetilde{M}_{B,n} = \max_{b_{n} \leq k < a_{n}} \quad \max_{a_{n-1}\tilde{h}_{n-1} \leq y \leq \tilde{h}_{n} - a_{n-1}\tilde{h}_{n-1}} \sum_{j=0}^{\tilde{h}_{n-1}-1}|\muB(T^{k\tilde{h}_{n}+y}I_{n-1,j})|
\]
Then $\lim_{n \to \infty} \widetilde{M}_{B,n} = 0$.
\end{proposition}

We remark that one could strengthen Proposition \ref{midy} to also include the times when $|y_{n}| < b_{n}$ but we will not need that here (and it is more natural to include those cases in a later argument).

We first establish that, for the sublevels not pushed through the top of the next column, we can replace $T^{k_{n}\tilde{h}_{n} + y_{n}}$ by an ergodic-like average preserving the value of $y_{n}$.  Essentially the proof is the standard technique that sublevels pushed through the top of the next column are mixed (corresponding to the $\frac{k}{a_{n}}\epsilon$ term below) combined with a careful accounting of the sublevels for which that does not occur.
  
 \begin{lemma}\label{first}
 Let $T$ be a quasi-staircase transformation, $B$ a union of levels in some $C_{N}$, $n > N$, 
 $b_{n} \leq k < a_{n}$ and $0 \leq y < \tilde{h}_{n}$.  Let $\epsilon > 0$ such that $\sup_{t \geq b_{n}} \left(\int \left|\frac{1}{t}\sum_{i=0}^{t-1}\chi_{B}\circ T^{-i}\right|~d\mu + \frac{2}{t}\right) < \epsilon$.  Then
 \begin{align*}
 \sum_{j=a_{n}b_{n}+b_{n+1}+c_{n+1}-c_{n}}^{\tilde{h}_{n}-y} &\left|\muB(T^{k\tilde{h}_{n}+y}I_{n,j}) - \frac{a_{n}-k}{r_{n}+1}\sum_{\ell=0}^{b_{n}-1} \muB(T^{y-k\ell}I_{n,j})\right| < \frac{k}{a_{n}}\epsilon; \quad\quad\text{and} \\
  \sum_{j=a_{n}b_{n}+b_{n+1}+c_{n+1}-c_{n}+\tilde{h}_{n}-y}^{\tilde{h}_{n}} &\left|\muB(T^{k\tilde{h}_{n}+y}I_{n,j}) - \frac{a_{n}-k-1}{r_{n}+1}\sum_{\ell=0}^{b_{n}-1} \muB(T^{y-\tilde{h}_{n}-(k+1)\ell}I_{n,j})\right| < \frac{k+1}{a_{n}}\epsilon
 \end{align*}
  \end{lemma}
 \begin{proof}
For $a_{n}b_{n} + b_{n+1} + c_{n+1} - c_{n} \leq j < \tilde{h}_{n}-y$, by Lemmas \ref{mixtrick2} and \ref{mixtrick3},
 \begin{align*}
  \muB(T^{k\tilde{h}_{n}+y}I_{n,j}) &= \sum_{i=0}^{a_{n}-1}\sum_{\ell=0}^{b_{n}-1} \muB(T^{k\tilde{h}_{n}+y}I_{n,j}^{[\ell a_{n}+i]}) + \muB(T^{k\tilde{h}_{n}+y}I_{n,j}^{[r_{n}]}) \\
 &=  
  \sum_{i=0}^{a_{n}-k-1} \sum_{\ell=0}^{b_{n}-1} \muB(T^{-k\ell}I_{n,j+y}^{[\ell a_{n} + i + k]})
 + \sum_{i=a_{n}-k}^{a_{n}-1} \sum_{\ell=0}^{b_{n}-2} \muB(T^{-k\ell - (i+k-a_{n})}I_{n,j+y}^{[\ell a_{n} + i + k+1]}) \\
 &\quad\quad\quad\quad
 + \sum_{i=0}^{k} \muB(T^{k\tilde{h}_{n}+y}I_{n,j}^{[r_{n}-i]}) 
 \end{align*}
 and since $k\ell \leq a_{n}b_{n}$ and $j+y \geq j \geq a_{n}b_{n}$, using Lemma \ref{L:mixtrick},
 \[
 \sum_{i=0}^{a_{n}-k-1} \sum_{\ell=0}^{b_{n}-1} \muB(T^{-k\ell}I_{n,j+y}^{[\ell a_{n} + i + k]})
 = \frac{1}{r_{n}+1} \sum_{i=0}^{a_{n}-k-1} \sum_{\ell=0}^{b_{n}-1} \muB(T^{-k\ell}I_{n,j+y})
 = \frac{a_{n}-k}{r_{n}+1} \sum_{\ell=0}^{b_{n}-1} \muB(T^{y-k\ell}I_{n,j})
 \]
 Likewise, since $k\ell + (i+k-a_{n}) \leq a_{n}b_{n}$,
 \begin{align*}
&\left| \sum_{i=a_{n}-k}^{a_{n}-1} \sum_{\ell=0}^{b_{n}-2} \muB(T^{-k\ell - (i+k-a_{n})}I_{n,j+y}^{[\ell a_{n} + i + k+1]}) \right|
 = \left| \frac{1}{r_{n}+1}  \sum_{i=a_{n}-k}^{a_{n}-1} \sum_{\ell=0}^{b_{n}-2} \muB(T^{-k\ell - (i+k-a_{n})}I_{n,j+y}) \right| \\
 &\quad\quad\quad\quad = \left|\frac{1}{r_{n}+1}  \sum_{i=0}^{k-1} \sum_{\ell=0}^{b_{n}-2} \muB(T^{y-k\ell-i}I_{n,j})\right|
 \leq \frac{k}{r_{n}+1}\sum_{\ell=0}^{b_{n}-2}\int_{T^{y-k\ell}I_{n,j}} \left| \frac{1}{k}\sum_{i=0}^{k-1} \chi_{B}\circ T^{-i}\right|~d\mu
 \end{align*}
 and therefore
 \[
 \sum_{j=a_{n}b_{n}+b_{n+1}+c_{n+1}-c_{n}}^{\tilde{h}_{n}-y} \left| \sum_{i=a_{n}-k}^{a_{n}-1} \sum_{\ell=0}^{b_{n}-2} \muB(T^{-k\ell - (i+k-a_{n})}I_{n,j+y}^{[\ell a_{n} + i + k+1]}) \right|
 < \frac{k(b_{n}-1)}{r_{n}+1}\int \left| \frac{1}{k}\sum_{i=0}^{k-1} \chi_{B}\circ T^{-i}\right|~d\mu
 \]

 For $0 \leq i \leq k-1$, using that $j \geq c_{n+1}-c_{n}+b_{n+1}+a_{n}b_{n}$ and that $I_{n,j}^{[0]} = I_{n+1,j}$,
 \begin{align*}
 T^{k\tilde{h}_{n}+y}I_{n,j}^{[r_{n}-i]}
 &= T^{k\tilde{h}_{n} + y + h_{n+1} - h_{n} - i(\tilde{h}_{n}+b_{n}-1)}I_{n,j}^{[0]} \\
&= T^{\tilde{h}_{n+1} + (k-i-1)\tilde{h}_{n} + c_{n} - c_{n+1} - i(b_{n}-1) + y}I_{n,j}^{[0]}
 = T^{\tilde{h}_{n+1}} I_{n+1,j+ (k-i-1)\tilde{h}_{n} + c_{n} - c_{n+1} - i(b_{n}-1) + y}
 \end{align*}
 therefore, since
 $|\muB(T^{\tilde{h}_{n+1}}I_{n+1,j^{\prime}})| = |\sum_{t=0}^{b_{n+1}-1}(\sum_{i=0}^{a_{n+1}-2}\muB(T^{-t}I_{n+1,j^{\prime}}^{[ta_{n+1}+i+1]}) + \muB(T^{-t-1}I_{n+1,j^{\prime}}^{[(t+1)a_{n+1}]})) + \muB(T^{\tilde{h}_{n+1}}I_{n+1,j^{\prime}}^{[r_{n+1}]})| \leq \frac{a_{n+1}}{r_{n+1}+1}|\sum_{t=0}^{b_{n+1}-1}\muB(T^{-t}I_{n+1,j^{\prime}})| + \frac{2\mu(I_{n+1,j^{\prime}})}{r_{n+1}+1}$ whenever $j^{\prime} \geq b_{n+1}$,
 \begin{align*}
&\left| \muB(T^{k\tilde{h}_{n}+y}I_{n,j}^{[r_{n}-i]}) \right|
\leq \left|\frac{a_{n+1}}{r_{n+1}+1}\sum_{t=0}^{b_{n+1}-1}\muB\left(T^{-t}I_{n+1,j+ (k-i-1)\tilde{h}_{n} + c_{n} - c_{n+1} - i(b_{n}-1) + y}\right)\right| + \frac{2\mu(I_{n+1})}{r_{n+1}+1} \\
&\quad\quad\quad\quad = \left|\frac{a_{n+1}}{r_{n+1}+1}\sum_{t=0}^{b_{n+1}-1}\muB\left(T^{-t}I_{n,j + c_{n} - c_{n+1} - i(b_{n}-1) + y}^{[k-i-1]}\right)\right| + \frac{2\mu(I_{n+1})}{r_{n+1}+1} \\
&\quad\quad\quad\quad  = \left|\frac{a_{n+1}}{r_{n+1}+1}\frac{1}{r_{n}+1}\sum_{t=0}^{b_{n+1}-1}\muB\left(T^{-t}I_{n,j + c_{n} - c_{n+1} - i(b_{n}-1) + y}\right)\right| + \frac{2\mu(I_{n+1})}{r_{n+1}+1} \\
&\quad\quad\quad\quad  \leq \frac{a_{n+1}b_{n+1}}{(r_{n+1}+1)(r_{n}+1)} \int_{T^{y+c_{n}-c_{n+1}-i(b_{n}-1)}I_{n,j}} \left|\frac{1}{b_{n+1}}\sum_{t=0}^{b_{n+1}-1}\chi_{B} \circ T^{-t}\right|~d\mu + \frac{2\mu(I_{n+1,j})}{r_{n+1}+1}
 \end{align*}
 and so
 \begin{align*}
 \sum_{j=a_{n}b_{n}+b_{n+1}+c_{n+1}-c_{n}}^{\tilde{h}_{n}-y} \sum_{i=0}^{k} &\left| \muB(T^{k\tilde{h}_{n}+y}I_{n,j}^{[r_{n}-i]}) \right| \\
&\leq \frac{k}{r_{n}+1} \int \left|\frac{1}{b_{n+1}}\sum_{t=0}^{b_{n+1}-1}\chi_{B} \circ T^{-t}\right|~d\mu + \frac{1}{r_{n}+1} + \frac{2}{(r_{n+1}+1)(r_{n}+1)}
 \end{align*}
 Therefore, since $\sup_{t \geq b_{n}} \left(\int \left|\frac{1}{t}\sum_{i=0}^{t-1}\chi_{B}\circ T^{-i}\right|~d\mu + \frac{2}{t}\right) < \epsilon$,
 \[
 \sum_{j=a_{n}b_{n}+b_{n+1}+c_{n+1}-c_{n}}^{\tilde{h}_{n}-y} \left|\muB(T^{k\tilde{h}_{n}+y}I_{n,j}) - \frac{a_{n}-k}{r_{n}+1}\sum_{\ell=0}^{b_{n}-1}\muB(T^{y-k\ell}I_{n,j})\right| 
 \leq \frac{k b_{n}}{r_{n}+1} \epsilon < \frac{k}{a_{n}}\epsilon
 \]
 
 For $a_{n}b_{n} + c_{n+1} - c_{n} + \tilde{h}_{n} - y \leq j < \tilde{h}_{n}$,
 \[
 T^{k\tilde{h}_{n}+y}I_{n,j} = T^{(k+1)\tilde{h}_{n} + 0}I_{n,j-(\tilde{h}_{n}-y)}
 \]
 and since $a_{n}b_{n} + b_{n+1} + c_{n+1} - c_{n} \leq j - (\tilde{h}_{n} - y) < \tilde{h}_{n} - 0$, the claim follows from the above replacing $k$ by $k+1$, $j$ by $j - (\tilde{h}_{n}-y)$ and $y$ by $0$.
 \end{proof}

We are now able to prove Proposition \ref{midy} using the previous lemma and the mixing nature of $\{ y_{n} \}$.  Unfortunately, the proof cannot be written as directly as one might hope: one cannot just use that $T^{y_{n}}$ is mixing on levels in $C_{n-1}$ directly.  Instead, the proof is essentially the same as the proof that $T^{y_{n}}$ is mixing but applied to the blocks of sublevels.  The techniques are standard for staircase mixing but care must be taken to keep track of the sublevels so there is a significant amount of bookkeeping.

\begin{proof}[Proof of Proposition \ref{midy}]
Let $\epsilon > 0$ such that $\sup_{t \geq b_{n}} \left(\int\left|\frac{1}{t}\sum_{i=0}^{t-1}\chi_{B}\circ T^{-i}\right|~d\mu + \frac{2}{t}\right) < \epsilon$.
Write $y = xa_{n-1}\tilde{h}_{n-1} + z\tilde{h}_{n-1} + w$ for $1 \leq x \leq b_{n}$ and $0 \leq z < a_{n-1}$ and $0 \leq w < \tilde{h}_{n-1}$.
Observe that if $0 \leq i < (b_{n-1} - x)a_{n-1}$ then $I_{n-1,j}^{[i]}$ is a level in $C_{n}$ below $I_{n,\tilde{h}_{n}-y}$ and that if $(b_{n-1} - x)a_{n-1} < i \leq r_{n-1}$ then $I_{n-1,j}^{[i]}$ is a level in $C_{n}$ above $I_{n,\tilde{h}_{n}-y}$.  Then by Lemma \ref{first}, as $\frac{2k+1}{a_{n}}\epsilon \leq \frac{3k}{a_{n}}\epsilon$,
\begin{align}
&\sum_{j=0}^{\tilde{h}_{n-1}-1}\Big{|}\muB(T^{k\tilde{h}_{n}+y}I_{n-1,j}) - \frac{a_{n}-k}{r_{n}+1}\sum_{\ell=0}^{b_{n}-1}\sum_{i=0}^{(b_{n-1}-x)a_{n-1}-1}\muB(T^{y-k\ell}I_{n-1,j}^{[i]}) \tag{$\dagger$} \\
& - \frac{a_{n}-k-1}{r_{n}+1}\sum_{\ell=0}^{b_{n}-1}\sum_{i=(b_{n-1}-x+1)a_{n-1}}^{r_{n-1}}\muB(T^{y-\tilde{h}_{n}-(k+1)\ell}I_{n-1,j}^{[i]})\Big{|}
< \frac{3k}{a_{n}}\epsilon + \frac{a_{n}}{r_{n}+1} + \frac{4(a_{n}b_{n} + b_{n+1} + c_{n+1})}{\tilde{h}_{n}} \nonumber
\end{align}
Now observe that, via Lemma \ref{mixingESlemma}, writing $k^{\prime} = xa_{n-1}+z$,
\begin{align*}
&\sum_{j=0}^{\tilde{h}_{n-1}-1}\left|\frac{1}{b_{n}}\sum_{\ell=0}^{b_{n}-1}\sum_{i=0}^{(b_{n-1}-x)a_{n-1}-1}\muB(T^{y-k\ell}I_{n-1,j}^{[i]})\right|
\leq \frac{1}{b_{n}}\sum_{\ell=0}^{b_{n}-1}\sum_{j=0}^{\tilde{h}_{n-1}-1}\left|\sum_{i=0}^{(b_{n-1}-x)a_{n-1}-1}\muB(T^{y-k\ell}I_{n-1,j}^{[i]})\right| \\
&\quad\quad \leq \frac{c_{n-1}}{\tilde{h}_{n-1}} 
+ \sum_{j=0}^{\tilde{h}_{n-1}-1}\left(\left|\sum_{i=0}^{(b_{n-1}-x)a_{n-1}-1}\muB(T^{k^{\prime}\tilde{h}_{n-1}}I_{n-1,j}^{[i]})\right| + \left|\sum_{i=0}^{(b_{n-1}-x)a_{n-1}-1}\muB(T^{(k^{\prime}+1)\tilde{h}_{n-1}}I_{n-1,j}^{[i]})\right| \right)
\end{align*}
which are precisely the sums $(\star)$ in the proof Proposition \ref{mixingES} (since $x \geq 1$ so $k^{\prime} \geq a_{n-1}$).  Therefore
\[
\sum_{j=0}^{\tilde{h}_{n-1}-1}\left|\frac{a_{n}-k}{r_{n}+1}\sum_{\ell=0}^{b_{n}-1}\sum_{i=0}^{(b_{n-1}-x)a_{n-1}-1}\muB(T^{y-k\ell}I_{n-1,j}^{[i]})\right|  \tag{$\ddagger$}
\to 0
\]
Now observe that for $0 \leq i < a_{n-1}$ and $0 \leq q < b_{n-1}$,
\[
I_{n-1,j}^{[qa_{n-1}+i]} = T^{qa_{n-1}\tilde{h}_{n-1} + \frac{1}{2}a_{n-1}q(q-1) + i\tilde{h}_{n-1} + iq}I_{n-1,j}^{[0]}
\]
so for $0 \leq i < a_{n-1}-1$, as $(b_{n-1}-x+q)(b_{n-1}-x+q-1) - q(q-1) = (b_{n-1}-x)(b_{n-1}-x-1+2q)$,
\[
I_{n-1,j}^{[(b_{n-1}-x+q)a_{n-1}+i+1]}
= T^{(b_{n-1}-x)a_{n-1}\tilde{h}_{n-1} + \frac{1}{2}a_{n-1}(b_{n-1}-x)(b_{n-1}-x-1+2q) + \tilde{h}_{n-1} + q + (i+1)(b_{n-1}-x)} I_{n-1,j}^{[qa_{n-1}+i]}
\]
Set $Q = Q_{q} = -c_{n} + c_{n-1} - (k+1)\ell + \frac{1}{2}a_{n-1}(b_{n-1}-x)(b_{n-1}-x-1+2q) + q - \frac{1}{2}a_{n-1}b_{n-1}(b_{n-1}-1) + b_{n-1}-x$ and note that $|Q| \leq c_{n} + a_{n}b_{n} + 2a_{n-1}b_{n-1}^{2}$.  Then, since $b_{n-1}a_{n-1}\tilde{h}_{n-1} + \tilde{h}_{n-1} - \tilde{h}_{n} = -c_{n} + c_{n-1} - \frac{1}{2}a_{n-1}b_{n-1}(b_{n-1}-1)$,
\[
T^{y-(k+1)\ell-\tilde{h}_{n}}I_{n-1,j}^{[(b_{n-1}-x+q)a_{n-1}+i]}
= T^{z\tilde{h}_{n-1} + w + i(b_{n-1}-x) + Q}I_{n-1,j}^{[qa_{n-1}+i]}
\]
Consider $j$ such that $0 \leq j + Q - a_{n-1}b_{n-1} < \tilde{h}_{n-1} - w - a_{n-1}b_{n-1}$.  If $z + i \geq a_{n-1}$,
\begin{align*}
T^{y-(k+1)\ell-\tilde{h}_{n}}I_{n-1,j}^{[(b_{n-1}-x+q)a_{n-1}+i]}
&= T^{z\tilde{h}_{n-1}}I_{n-1,j+Q+w+i(b_{n-1}-x)}^{[qa_{n-1}+i]}
= I_{n-1,j+Q+w+i(b_{n-1}-x)-zq-(z+i-a_{n-1})}^{[qa_{n-1}+i+z]} \\
&= T^{i(b_{n-1}-x-1)}I_{n-1,j+Q+w-zq-(z-a_{n-1})}^{[qa_{n-1}+i+z]}
\end{align*}
and therefore
\[
\muB(T^{y-(k+1)\ell-\tilde{h}_{n}}I_{n-1,j}^{[qa_{n-1}+i]})
= \frac{1}{r_{n-1}}\muB\left(T^{i(b_{n-1}-x-1)}I_{n-1,j+Q+w-zq-(z-a_{n-1})}\right)
\]
Similarly, if $z + i < a_{n-1}$,
\begin{align*}
T^{y-(k+1)\ell-\tilde{h}_{n}}I_{n-1,j}^{[(b_{n-1}-x+q)a_{n-1}+i]}
&= T^{z\tilde{h}_{n-1}}I_{n-1,j+Q+w+i(b_{n-1}-x)}^{[qa_{n-1}+i]}
= I_{n-1,j+Q+w+i(b_{n-1}-x)-zq}^{[qa_{n-1}+i+z]} \\
&= T^{i(b_{n-1}-x)}I_{n-1,j+Q+w-zq}^{[qa_{n-1}+i+z]}
\end{align*}
so
\[
\muB(T^{y-(k+1)\ell-\tilde{h}_{n}}I_{n-1,j}^{[qa_{n-1}+i]})
= \frac{1}{r_{n-1}}\muB\left(T^{i(b_{n-1}-x)}I_{n-1,j+Q+w-zq-(z-a_{n-1})}\right)
\]
Therefore, as $\frac{x}{r_{n-1}} \leq \frac{b_{n-1}}{r_{n-1}} < \frac{1}{a_{n-1}}$,
\begin{align*}
&\sum_{j=a_{n-1}b_{n-1}-Q}^{\tilde{h}_{n-1}-w-a_{n-1}b_{n-1}}\left|\sum_{i=(b_{n-1}-x+1)a_{n-1}}^{r_{n-1}} \muB(T^{y-(k+1)\ell-\tilde{h}_{n}}I_{n-1,j}^{[i]})\right| \\
&\quad\quad\leq \sum_{j=a_{n-1}b_{n-1}-Q}^{\tilde{h}_{n-1}-w-a_{n-1}b_{n-1}} \left|\sum_{q=b_{n-1}-x+1}^{b_{n-1}-1} \sum_{i=0}^{a_{n-1}-2} \muB(T^{y-(k+1)\ell-\tilde{h}_{n}}I_{n-1,j}^{[qa_{n-1}+i]})\right| + \frac{x+1}{r_{n-1}} \\
&\quad\quad\leq \frac{1}{r_{n-1}}\sum_{q=0}^{x-1}\int\left|\sum_{i=0}^{a_{n-1}-z-1}\chi_{B}\circ T^{i(b_{n-1}-x-1)}\right|~d\mu
+ \frac{1}{r_{n-1}}\sum_{q=0}^{x-1}\int\left|\sum_{i=0}^{z-1}\chi_{B}\circ T^{i(b_{n-1}-x)}\right|~d\mu + \frac{x+1}{r_{n-1}} \\
&\quad\quad\leq \int\left|\frac{1}{a_{n-1}}\sum_{i=0}^{a_{n-1}-z-1}\chi_{B}\circ T^{i(b_{n-1}-x-1)}\right|~d\mu
+ \int\left|\frac{1}{a_{n-1}}\sum_{i=0}^{z-1}\chi_{B}\circ T^{i(b_{n-1}-x)}\right|~d\mu + \frac{x+1}{r_{n-1}}
\end{align*}

Now consider $j$ such that $\tilde{h}_{n-1} - w + a_{n-1}b_{n-1} - Q \leq j < \tilde{h}_{n-1} - a_{n-1}b_{n-1}$.  Then
\[
T^{y-(k+1)\ell-\tilde{h}_{n}}I_{n-1,j}^{[qa_{n-1}+i]} = T^{(z+1)\tilde{h}_{n-1}}I_{n-1,j+Q+w+i(b_{n-1}-x)-\tilde{h}_{n-1}}
\]
so similar reasoning as above shows that
\begin{align*}
&\sum_{j=\tilde{h}_{n-1} - w + a_{n-1}b_{n-1} - Q}^{\tilde{h}_{n-1}-a_{n-1}b_{n-1}}\left|\sum_{i=(b_{n-1}-x+1)a_{n-1}}^{r_{n-1}} \muB(T^{y-(k+1)\ell-\tilde{h}_{n}}I_{n-1,j}^{[i]})\right| \\
&\quad\quad\leq \int\left|\frac{1}{a_{n-1}}\sum_{i=0}^{a_{n-1}-z-1}\chi_{B}\circ T^{i(b_{n-1}-x-1)}\right|~d\mu
+ \int\left|\frac{1}{a_{n-1}}\sum_{i=0}^{z-1}\chi_{B}\circ T^{i(b_{n-1}-x)}\right|~d\mu + \frac{x+1}{r_{n-1}}
\end{align*}

Note that $y < \tilde{h}_{n} - a_{n-1}\tilde{h}_{n-1} = (b_{n-1}-1)a_{n-1}\tilde{h}_{n-1} + h_{n-1} + \frac{1}{2}a_{n-1}b_{n-1}(b_{n-1}-1) + c_{n} < (b_{n-1}-1)a_{n-1}\tilde{h}_{n-1} + 2\tilde{h}_{n-1}$.  Therefore $x \leq b_{n-1} - 1$ and if $x = b_{n-1}-1$ then $z \leq 1$.  When $x \leq b_{n-1} - 1$, both $b_{n-1}-x \geq 1$ and $b_{n-1}-x-1 \geq 1$ so both integrals tend to zero by Proposition \ref{P:pe}.  When $x = b_{n-1}-1$, the first integral tends to zero by Proposition \ref{P:pe} and the second is bounded by $\frac{z}{a_{n-1}} \to 0$.

Since $\frac{|Q|}{\tilde{h}_{n-1}} \leq \frac{c_{n} + a_{n}b_{n} + a_{n-1}b_{n-1}^{2}}{\tilde{h}_{n-1}} \to 0$, then
\[
\sum_{j=0}^{\tilde{h}_{n-1}}\left|\sum_{i=(b_{n-1}-x+1)a_{n-1}}^{r_{n-1}} \muB(T^{y-(k+1)\ell-\tilde{h}_{n}}I_{n-1,j}^{[i]})\right| \to 0
\] 
Therefore equation $(\dagger)$ gives that $\sum_{j=0}^{\tilde{h}_{n-1}-1}|\muB(T^{k\tilde{h}_{n}+y}I_{n-1,j})| \to 0$ as both the above quantity and that in $(\ddagger)$ tend to $0$.  Since this holds uniformly over $k$ and $y$ in the specified range, the claim follows.
\end{proof}

\subsubsection{Mixing for the remaining times}

The times not covered by the previous cases, namely those of the form $k_{n}\tilde{h}_{n} + y_{n}$ where $b_{n} \leq k_{n} < a_{n}$ and $|y_{n}| < a_{n-1}\tilde{h}_{n-1}$, require more care and the introduction of many new ideas; however, the majority of the argument is to split the sequence into cases where, in each case, a (generalization of a) standard mixing proof technique can be applied.

The final goal of this section is to prove the following proposition which ensures mixing on all remaining times.

 \begin{proposition}\label{fourth}
 Let $T$ be a quasi-staircase transformation such that $\sum \frac{a_{n}b_{n}+b_{n+1}+c_{n+1}}{h_{n}} < \infty$ and $\sum \frac{1}{a_{n}b_{n}} < \infty$ and $\frac{a_{n}b_{n}^{2}}{h_{n}} \to 0$ and $\frac{a_{n+1}b_{n+1}}{h_{n}} \to 0$.  Let $B$ be a union of levels in some fixed $C_{N}$.  Then
 \[
 \lim_{n \to \infty} \quad \max_{b_{n} \leq k < a_{n}} \quad \max_{|q| < a_{n-1}\tilde{h}_{n-1}} \left| \muB(T^{k\tilde{h}_{n} + q}B) \right| = 0.
 \]
 \end{proposition}
 
 For ease of exposition, we denote the measure of the levels which may be `safely' ignored as follows.

\begin{notation}
Define $\tau_{n} = \frac{4(a_{n}b_{n} + b_{n+1} + c_{n+1})}{\tilde{h}_{n}}$.
\end{notation}

Our first lemma effectively generalizes the known technique of splitting the tower into the rightmost subcolumns, which will be mixed due to pushing through the top of the next tower, and those on the left which will need to be handled differently.

\begin{lemma}\label{bigsmally}
Let $T$ be a quasi-staircase transformation, $B$ a union of levels in some $C_{N}$, $\epsilon > 0$ such that $\sup_{t\geq b_{N}}\left(\int \left|\frac{1}{t}\sum_{i=0}^{t-1}\chi_{B}\circ T^{-i}\right|~d\mu + \frac{2}{t}\right) < \frac{\epsilon}{3}$, $n > N$, $b_{n} \leq k < a_{n}$ and $0 \leq |y| < a_{n-1}\tilde{h}_{n-1}$.  Then
\[
\left|\muB(T^{k\tilde{h}_{n}+y}B) - \frac{a_{n}-k}{r_{n}+1}\sum_{\ell=0}^{b_{n}-1}\muB(T^{y-k\ell}B)\right| \leq \frac{k}{a_{n}}\epsilon + \tau_{n} + \left(1 - \frac{k}{a_{n}}\right)\frac{1}{b_{n}}\sum_{\ell=0}^{b_{n}-1}\frac{2|y - k\ell|}{\tilde{h}_{n}}
\]
\end{lemma}
\begin{proof}
Consider first when $y \geq 0$.  Write $\beta = \{ a_{n}b_{n} + b_{n+1} + c_{n+1} \leq j < \tilde{h}_{n} - y : I_{n,j} \subseteq B \}$ and $\beta^{\prime} = \{ a_{n}b_{n} + b_{n+1} + c_{n+1} + \tilde{h}_{n} - y \leq j < \tilde{h}_{n} : I_{n,j} \subseteq B \}$ .  By Lemma \ref{first},
\begin{align*}
&\left|\sum_{j\in\beta\cup\beta^{\prime}}\muB(T^{k\tilde{h}_{n}+y}I_{n,j}) - \sum_{\ell=0}^{b_{n}-1}\left(\frac{a_{n}-k}{r_{n}+1}  \sum_{j\in\beta} \muB(T^{y-k\ell}I_{n,j}) - \frac{a_{n}-k-1}{r_{n}+1}\sum_{j\in\beta^{\prime}}\muB(T^{y-\tilde{h}_{n}-(k+1)\ell}I_{n,j})\right)\right|
\end{align*}
is bounded by $\frac{k}{a_{n}}\frac{\epsilon}{3} + \frac{k+1}{a_{n}}\frac{\epsilon}{3} \leq \frac{k\epsilon}{a_{n}}$ and therefore
\begin{align*}
&\left|\muB(T^{k\tilde{h}_{n}+y}B) - \frac{a_{n}-k}{r_{n}+1}\sum_{\ell=0}^{b_{n}-1}\sum_{j\in\beta \cup \beta^{\prime}}\muB(T^{y-k\ell}I_{n,j})\right| \\
&\quad\quad \leq \frac{k}{a_{n}}\epsilon + \frac{\tau_{n}}{2} + \frac{a_{n}-k}{r_{n}+1}\sum_{\ell=0}^{b_{n}-1}\sum_{j \in \beta^{\prime}}\left|\muB(T^{y-k\ell}I_{n,j}) - \frac{a_{n}-k-1}{a_{n}-k}\muB(T^{y-\tilde{h}_{n}-(k+1)\ell}I_{n,j})\right| \\
&\quad\quad \leq \frac{k}{a_{n}}\epsilon + \frac{\tau_{n}}{2} + \frac{a_{n}-k}{r_{n}+1}b_{n}|\beta^{\prime}|\mu(I_{n})\frac{2a_{n}-2k-1}{a_{n}-k}
< \frac{k}{a_{n}}\epsilon + \frac{\tau_{n}}{2} + \left(1 - \frac{k}{a_{n}}\right)\frac{2|\beta^{\prime}|}{\tilde{h}_{n}}
\end{align*}
so the claim follows for $y \geq 0$ as $|\beta^{\prime}| = y-a_{n}b_{n}-c_{n+1} \leq |y - k \ell|$ for all $0 \leq \ell < b_{n}$ (and if $y < a_{n}b_{n}+b_{n+1}+c_{n+1}$ then $\beta^{\prime} = \emptyset$) and since $|\muB(T^{k\tilde{h}_{n}+y}B) - \sum_{j\in\beta\cup\beta^{\prime}}\muB(T^{k\tilde{h}_{n}+y}I_{n,j})| \leq \frac{\tau_{n}}{2}$.

Now consider when $y < 0$.  Then $k\tilde{h}_{n} + y = (k-1)\tilde{h}_{n} + (\tilde{h}_{n}+y)$ so, following the same reasoning as above and swapping the roles of $\beta^{\prime}$ and $\beta$,
\begin{align*}
&\left|\muB(T^{k\tilde{h}_{n}+y}B) - \frac{a_{n}-(k-1)-1}{r_{n}+1}\sum_{\ell=0}^{b_{n}-1}\sum_{j\in\beta \cup \beta^{\prime}}\muB(T^{(y+\tilde{h}_{n})-(k-1+1)\ell}I_{n,j})\right| \\
&\quad\quad\quad\quad < \frac{k-1+1}{a_{n}}\epsilon + \frac{\tau_{n}}{2} + \left(1 - \frac{k-1+1}{a_{n}}\right)\frac{2|\beta|}{\tilde{h}_{n}}
\leq \frac{k}{a_{n}}\epsilon + \tau_{n} + \left(1 - \frac{k}{a_{n}}\right) \frac{2|\beta|}{\tilde{h}_{n}}
\end{align*}
Since in this case $|\beta| \leq |y - k\ell|$ for all $0 \leq \ell < b_{n}$, then
$
|\beta| \leq \frac{1}{b_{n}}\sum_{\ell=0}^{b_{n}-1} |y - k \ell|
$
so the claim follows.
\end{proof}

Our next lemma is a generalization of the now-standard Block Lemma \cite{adams1998smorodinsky}.  It is a `weighted' version of that lemma, though we emphasize that it is not the case that the weights sum to $1$ but rather that each term in the average may be weighted by a value between $0$ and $1$ (so e.g.~the case when most of the weights are $0$ follows trivially).

\begin{lemma}\label{deltaBL}
Let $\epsilon > 0$ and $q,k,p,Q,L \in \mathbb{N}$ and for all $0 \leq \ell < L$, let $0 \leq \delta_{\ell} \leq 1$.  If $\frac{pQ}{L} < \epsilon$ and $\frac{1}{Q} < \epsilon$ and $|\muB(T^{kpt}B)| < \epsilon$ for all $1 \leq t < Q$ then
\[
\left|\frac{1}{L}\sum_{\ell=0}^{L-1}\delta_{\ell}\muB(T^{q-k\ell}B)\right|
< (2\epsilon)^{1/2} + \epsilon
\]
\end{lemma}
\begin{proof}
Using that $T$ is measure-preserving and the Cauchy-Schwarz inequality,
\begin{align*}
&\left|\frac{1}{L}\sum_{\ell=0}^{L-1}\delta_{\ell}\muB(T^{q-k\ell}B)\right|
= \left|\frac{1}{L}\sum_{\ell=0}^{L-1}\int_{B} \delta_{\ell}\chi_{B}\circ T^{q-k\ell}~d\mu\right|
\leq \int \left|\frac{1}{L}\sum_{\ell=0}^{L-1}\delta_{\ell}\chi_{B}\circ T^{q-k\ell}\right|~d\mu \\
&\quad \leq \frac{pQ\floor{\frac{L}{pQ}}}{L} \frac{1}{\floor{\frac{L}{pQ}}} \sum_{j=0}^{\floor{\frac{L}{pQ}}-1} \frac{1}{p}\sum_{i=0}^{p-1} \int \left|\frac{1}{Q}\sum_{t=0}^{Q-1} \delta_{jpQ + i + pt} \chi_{B} \circ T^{-kpt} \right| \circ T^{q - kjpQ - ki}~d\mu + \frac{pQ}{L} \\
&\quad < \frac{1}{\floor{\frac{L}{pQ}}} \sum_{j=0}^{\floor{\frac{L}{pQ}}-1} \frac{1}{p}\sum_{i=0}^{p-1} \int \left|\frac{1}{Q}\sum_{t=0}^{Q-1} \delta_{jpQ + i + pt} \chi_{B} \circ T^{-kpt} \right|~d\mu + \epsilon \\
&\quad \leq \frac{1}{\floor{\frac{L}{pQ}}} \sum_{j=0}^{\floor{\frac{L}{pQ}}-1} \frac{1}{p}\sum_{i=0}^{p-1} \left( \int \left|\frac{1}{Q}\sum_{t=0}^{Q-1} \delta_{jpQ + i + pt} \chi_{B} \circ T^{-kpt} \right|^{2}~d\mu \right)^{1/2} + \epsilon \\
&\quad = \frac{1}{\floor{\frac{L}{pQ}}} \sum_{j=0}^{\floor{\frac{L}{pQ}}-1} \frac{1}{p}\sum_{i=0}^{p-1} \left( \frac{1}{Q^{2}}\sum_{t,u=0}^{Q-1} \delta_{jpQ + i + pt}\delta_{jpQ + i + pu} \muB(T^{kp(t-u)}B) \right)^{1/2} + \epsilon \\
&\quad = \frac{1}{\floor{\frac{L}{pQ}}} \sum_{j=0}^{\floor{\frac{L}{pQ}}-1} \frac{1}{p}\sum_{i=0}^{p-1} \left( 
\frac{1}{Q^{2}}\sum_{t=0}^{Q-1} \delta_{jpQ + i + pt}^{2} \muB(B) + 
\frac{1}{Q^{2}}\sum_{t\neq u} \delta_{jpQ + i + pt}\delta_{jpQ + i + pu} \muB(T^{kp(t-u)}B) \right)^{1/2} + \epsilon \\
&\quad < \frac{1}{\floor{\frac{L}{pQ}}} \sum_{j=0}^{\floor{\frac{L}{pQ}}-1} \frac{1}{p}\sum_{i=0}^{p-1} \left( 
\frac{1}{Q} + 
\frac{1}{Q^{2}}\sum_{t\neq u} \delta_{jpQ + i + pt}\delta_{jpQ + i + pu} \epsilon \right)^{1/2} + \epsilon
\leq \left(\frac{1}{Q} + \frac{1}{Q^{2}}\sum_{t\ne u}\epsilon\right)^{1/2} + \epsilon \qedhere
\end{align*}
\end{proof}

Using the `weighted' Block Lemma, we likewise generalize the Blum-Hanson trick (see e.g.~\cite{adams1998smorodinsky}) to a `weighted' version.

\begin{proposition}\label{deltape}
Let $T$ be a quasi-staircase transformation such that $\frac{b_{n}^{2}}{h_{n}} \to 0$, $\frac{a_{n}b_{n}}{\tilde{h}_{n}} \to 0$ and $\frac{b_{n}}{a_{n}} \to 0$.
Let $B$ be a union of levels in some column $C_{N_{0}}$.  Then
\[
\lim_{N\to\infty} \quad \max_{0 \leq \delta_{\ell} \leq 1} \quad \max_{1 \leq k \leq N} \int \left| \frac{1}{N}\sum_{\ell=0}^{N-1} \delta_{\ell}\chi_{B}\circ T^{-\ell k}\right|~d\mu = 0
\]
\end{proposition}
\begin{proof}
Fix $\epsilon > 0$.  Let $m$ such that $b_{m} \geq 2\ceiling{\epsilon^{-1}}$, $\frac{4(r_{m}+1)\ceiling{\epsilon^{-1}}^{2}}{\tilde{h}_{m}} < \epsilon$ and $\sup_{n \geq m} \widehat{M}_{B,n} < \epsilon$ (using Proposition \ref{mixingES2}).  Take any $N$ such that $\frac{\tilde{h}_{m}\ceiling{\epsilon^{-1}}}{N} < \epsilon$.
Let $k$ and $\delta_{\ell}$ attain the maximum for $N$.

Consider first the case when $k \geq \tilde{h}_{m}$.  Let $n \geq m$ such that $\tilde{h}_{n} \leq k < \tilde{h}_{n+1}$.  Let $p$ such that $(p-1)k < \tilde{h}_{n+1} \leq pk$ so that $pk < \tilde{h}_{n+1} + k < 2\tilde{h}_{n+1}$.  Then for every $1 \leq q < \ceiling{\epsilon^{-1}}$, $\tilde{h}_{n+1} \leq qpk < \ceiling{\epsilon^{-1}}2\tilde{h}_{n+1} \leq b_{n}\tilde{h}_{n}$ meaning that $|\muB(T^{qpk}B)| \leq \widehat{M}_{B,n} < \epsilon$.  Now
\[
\frac{p\ceiling{\epsilon^{-1}}}{N} = \frac{pk\ceiling{\epsilon^{-1}}}{Nk}
< \frac{2\tilde{h}_{n+1}\ceiling{\epsilon^{-1}}}{N\tilde{h}_{n}}
< \frac{4(r_{n}+1)\ceiling{\epsilon^{-1}}}{N}
\leq \frac{4(r_{n}+1)\ceiling{\epsilon^{-1}}}{k}
\leq \frac{4(r_{n}+1)\ceiling{\epsilon^{-1}}}{\tilde{h}_{n}} < \epsilon
\]
so Lemma \ref{deltaBL} implies that $\int \left| \frac{1}{N}\sum_{\ell=0}^{N-1} \delta_{\ell}\chi_{B}\circ T^{-\ell k}\right|~d\mu < (2\epsilon)^{1/2} + \epsilon$.

Consider now when $k < \tilde{h}_{m}$.  Let $p$ such that $(p-1)k < \tilde{h}_{m} \leq pk$ so that $pk < 2\tilde{h}_{m}$ and $p \leq \tilde{h}_{m}$.  Then $\tilde{h}_{m} \leq qpk < \ceiling{\epsilon^{-1}}2\tilde{h}_{m} \leq b_{m}\tilde{h}_{m}$ for $1 \leq q < \ceiling{\epsilon^{-1}}$ so $|\muB(T^{qpk}B)| \leq \widehat{M}_{B,m} < \epsilon$.  Since $\frac{p\ceiling{\epsilon^{-1}}}{N} < \frac{\tilde{h}_{m}\ceiling{\epsilon^{-1}}}{N} < \epsilon$, Lemma \ref{deltaBL} again implies that $\int \left| \frac{1}{N}\sum_{\ell=0}^{N-1} \delta_{\ell}\chi_{B}\circ T^{-\ell k}\right|~d\mu < (2\epsilon)^{1/2} + \epsilon$.
\end{proof}

\newcommand{\gam}[4]{\gamma_{#4}^{#1,#2,#3}}

Since we will have need to work with times not all ranging within the same interval $[h_{n},h_{n+1})$ (which is all one needs to consider when working with staircases), we introduce the following notation to denote which column a time `naturally wants to be thought of as applying to'.

\begin{notation}
 For $t \in \mathbb{Z}$, write $\alpha(t)$ for the unique positive integer such that $\tilde{h}_{\alpha(t)} \leq |t| < \tilde{h}_{\alpha(t)+1}$.
 \end{notation}
 
 The next lemma is the main ingredient in the proof of mixing for our remaining times.  Given a `weighted' average of times which form an arithmetic progression (meaning the times are $q - \ell k$ for some fixed $q$ and $k$ as $\ell$ ranges from $0$ to $L-1$), the lemma states that either the average is already mixed by the previously established facts or that the sequence of times which cannot be guaranteed to already be mixing times has a very specific structure closely resembling that of an arithmetic progression, along with control on the size of the gaps in the potentially `bad' times $\mathcal{L}$.  This structure, along with the control on the size of the gaps, will allow us to deduce mixing by dropping to previous towers in a suitable manner.
 
 The proof breaks into cases.  The bulk of them, where the average is shown to be mixed, follow from a combination of the already established (weak) power ergodicity (similar to how mixing is deduced for staircases) and the generalization of the Blum-Hanson trick combined with the generalized Block Lemma.  Due to the length of the proof, we include expository text throughout, italicized to distinguish it.

\begin{lemma}\label{thirdA}
Let $\epsilon > 0$ and set $\epsilon_{0} = (2\ceiling{\epsilon^{-1}}^{\ceiling{\epsilon^{-1}}+1})^{-1}$.  Let $L,k,q \in \mathbb{Z}$ with $L \geq \epsilon_{0}^{-1}$ and for each $0 \leq \ell < L$, let $0 \leq \delta_{\ell} \leq 1$.  

Let $\alpha_{0} = \max \{ \alpha(q - \ell k) : 0 \leq \ell < L \}$.
Assume that $\max(M_{B,\alpha_{0}}, M_{B,\alpha_{0}-1}, \widehat{M}_{B,\alpha_{0}}, \widetilde{M}_{B,\alpha_{0}}) < \epsilon$ and $b_{\alpha_{0}-1} > 4 \epsilon^{-1}\epsilon_{0}^{-1}$.

Write
\quad\quad
$k = z\tilde{h}_{\alpha_{0}} + y$ \quad for $|y| \leq \frac{1}{2}\tilde{h}_{\alpha_{0}}$ \quad\quad and \quad\quad 
$ q = x\tilde{h}_{\alpha_{0}} + r$ \quad for $|r| \leq \frac{1}{2}\tilde{h}_{\alpha_{0}}$.

Let $k_{\ell},y_{\ell} \in \mathbb{Z}$ such that \quad $\displaystyle q - \ell k = k_{\ell}\tilde{h}_{\alpha_{0}} + y_{\ell}$ \quad with $\displaystyle |y_{\ell}| \leq \frac{1}{2}\tilde{h}_{\alpha_{0}}$ and define
\[
\mathcal{L} = \left\{ 0 \leq \ell < L : \left(k_{\ell} = 0~\text{or}~b_{\alpha_{0}} \leq |k_{\ell}| < a_{\alpha_{0}}\right)~\text{and}~|y_{\ell}| < a_{\alpha_{0}-1}\tilde{h}_{\alpha_{0}-1} \right\}
\]
Then at least one of the following holds:
\begin{itemize}
\item $\displaystyle \left|\frac{1}{L}\sum_{\ell=0}^{L-1}\delta_{\ell} \muB(T^{q - \ell k}B) \right| + \frac{1}{L}\sum_{\ell=0}^{L-1}(1 - \delta_{\ell})\epsilon < 4\epsilon^{1/2}$; or 
\item there exists $p  \in \mathbb{Z}$, $t > 0$ and $0 \leq \ell_{0} < L^{\prime} \leq L$ such that
\[
\mathcal{L} \subseteq \{ \ell_{0} + i t : 0 \leq i < L^{\prime} \} \quad\quad\quad\text{and}\quad\quad\quad |ity - ip\tilde{h}_{\alpha_{0}}| < \frac{1}{3}\tilde{h}_{\alpha_{0}} \quad \text{for all $0 \leq i < L^{\prime}$}
\]
\end{itemize}
\end{lemma}

\begin{proof}
First observe that if $0 \leq \ell < L$ and $\ell \notin \mathcal{L}$ then one of the following must hold:
\begin{itemize}
\item $a_{\alpha_{0}}\tilde{h}_{\alpha_{0}} \leq |q - \ell k| < \tilde{h}_{\alpha_{0}+1}$ (using that $\alpha_{0}$ is maximal);
\item $b_{\alpha_{0}}\tilde{h}_{\alpha_{0}} \leq |q - \ell k| < a_{\alpha_{0}}\tilde{h}_{\alpha_{0}}$ and $|y_{\ell}| \geq a_{\alpha_{0}-1}\tilde{h}_{\alpha_{0}-1}$;
\item $\tilde{h}_{\alpha_{0}} \leq |q - \ell k| < b_{\alpha_{0}}\tilde{h}_{\alpha_{0}}$; or
\item $k = 0$ (so $y_{\ell} = q - \ell k = q$) and $a_{\alpha_{0}-1}\tilde{h}_{\alpha_{0}-1} \leq |q - \ell k| < \tilde{h}_{\alpha_{0}}$.
\end{itemize}

\textit{Each of the above cases falls under one of the hypotheses on times which are already known to be $\epsilon$-mixed, thus $\mathcal{L}$ represents the `potentially bad times'.  As such, the proof need only focus on those values in $\mathcal{L}$ since the rest are already mixed.  We now make this concrete.}\\

Since $\muB(T^{-t}B) = \muB(T^{t}B)$,
\begin{itemize}
\item $a_{\alpha_{0}}\tilde{h}_{\alpha_{0}} \leq |q - \ell k| < \tilde{h}_{\alpha_{0}+1}$ implies $\left|\muB(T^{q - \ell k}B)\right| \leq M_{B,\alpha_{0}}$;
\item $b_{\alpha_{0}}\tilde{h}_{\alpha_{0}} \leq |q - \ell k| < a_{\alpha_{0}}\tilde{h}_{\alpha_{0}}$ and $|y_{\ell}| \geq a_{\alpha_{0}-1}\tilde{h}_{\alpha_{0}-1}$ implies $\left|\muB(T^{q - \ell k}B)\right| \leq \widetilde{M}_{B,\alpha_{0}}$;
\item $\tilde{h}_{\alpha_{0}} \leq |q - \ell k| < b_{\alpha_{0}}\tilde{h}_{\alpha_{0}}$ implies $\left|\muB(T^{q - \ell k}B)\right| \leq \widehat{M}_{B,\alpha_{0}}$; and
\item $a_{\alpha_{0}-1}\tilde{h}_{\alpha_{0}-1} \leq |q - \ell k| < \tilde{h}_{\alpha_{0}}$ implies $\left|\muB(T^{q - \ell k}B)\right| \leq M_{B,\alpha_{0}-1}$.
\end{itemize}
Therefore, by hypothesis, for $\ell \notin \mathcal{L}$, it holds that $\left|\muB(T^{q - \ell k}B)\right| < \epsilon$. 

In particular, if $|\mathcal{L}| < \epsilon L$ then 
\[
\left|\frac{1}{L}\sum_{\ell=0}^{L-1}\delta_{\ell} \muB(T^{q - \ell k}B) \right| + \frac{1}{L}\sum_{\ell=0}^{L-1}(1 - \delta_{\ell})\epsilon
< \frac{1}{L}|\mathcal{L}| + \frac{1}{L}\left(L - |\mathcal{L}|\right)\epsilon + \epsilon
< 3\epsilon
\]
so the first possible conclusion will hold.

\textit{The main idea is to consider the fraction $\frac{y}{\tilde{h}_{\alpha_{0}}}$ and break into cases.  We begin by approximating $\frac{y}{\tilde{h}_{\alpha_{0}}}$ with the closest fraction $\frac{p}{t}$ that has denominator less than $b_{\alpha_{0}-1}L$ (and then do the same to approximate $\frac{r}{\tilde{h}_{\alpha_{0}}}$).}

Let $p,t \in \mathbb{Z}$ with $0 < t < b_{\alpha_{0}-1}L$ such that
$\displaystyle
\left|\frac{y}{\tilde{h}_{\alpha_{0}}} - \frac{p}{t}\right| < \frac{1}{Lb_{\alpha_{0}-1}}
$
and either $(p=0, t=1)$ or $p,t$ are relatively prime.  Then let $u \in \mathbb{Z}$ such that
$\displaystyle
\left|\frac{r}{\tilde{h}_{\alpha_{0}}} - \frac{u}{t}\right| \leq \frac{1}{2t}
$.

\textit{In the case when $\frac{p}{t}$ has a large denominator and large numerator, we expect that $T^{\ell y}$ distributes sublevels somewhat evenly across many levels and so we can apply either (weak) power ergodicity or the (generalized) Blum-Hanson trick.}

\textit{In the case when $\frac{p}{t}$ is very small, $\ell |y| < \tilde{h}_{\alpha_{0}}$ so the set of $\ell$ for which $|y_{\ell}| \approx q - \ell y$ is less than $b_{\alpha_{0}-1}\tilde{h}_{\alpha_{0}-1}$ should be a consecutive set of integers (as there is no `wraparound' coming from  $\ell y$).}

\textit{In the case when $\frac{p}{t}$ has a small denominator, we expect $T^{\ell y}$ to send most of the sublevels into the same small set of levels.  This the arithmetic progression like structure conclusion of the lemma.  The idea (in a later lemma) for handling this case is to `drop down' to an earlier column, similar to in Section \ref{prevcol}.}

For $\ell \in \mathcal{L}$, since $q - \ell k = k_{\ell}\tilde{h}_{\alpha_{0}} + y_{\ell}$ for $|y_{\ell}| < a_{\alpha_{0}-1}\tilde{h}_{\alpha_{0}-1}$ and since
\[
q - \ell k = (x\tilde{h}_{\alpha_{0}} + r - \ell(z\tilde{h}_{\alpha_{)}} + y)
= (x - \ell z)\tilde{h}_{\alpha_{0}} + (r - \ell y)
\]
it follows that $r - \ell y$ must be of the form $n_{\ell}\tilde{h}_{\alpha_{)}} + o_{\ell}$ for some $n_{\ell},o_{\ell} \in \mathbb{Z}$ where $|o_{\ell}| < a_{\alpha_{0}-1}\tilde{h}_{\alpha_{0}-1}$.  Then $|r - \ell y - n_{\ell} \tilde{h}_{\alpha_{0}}| < a_{\alpha_{0}-1}\tilde{h}_{\alpha_{0}-1}$ so $\left|\frac{r - \ell y}{\tilde{h}_{\alpha_{0}}} - n_{\ell}\right| < \frac{1}{b_{\alpha_{0}-1}}$.  Therefore
\[
\left|\frac{u - \ell p}{t} - n_{\ell}\right| 
\leq \left| \frac{u - \ell p}{t} - \frac{r - \ell y}{\tilde{h}_{\alpha_{0}}}\right| + \left|\frac{r - \ell y}{\tilde{h}_{\alpha_{0}}} - n_{\ell}\right|
< \frac{1}{2t} + \frac{\ell}{L b_{\alpha_{0}-1}} + \frac{1}{b_{\alpha_{0}-1}} < \frac{1}{2t} + \frac{2}{b_{\alpha_{0}-1}}
\]

\textit{The above calculation shows that $\frac{u - \ell p}{t}$ is always very close to an integer when $\ell \in \mathcal{L}$.  This will play a key role in establishing the claimed structure of $\mathcal{L}$ in the case when $\frac{p}{t}$ has small denominator.}

\textit{We now split into cases.}
%
%\textit{The first case we consider is when $t$ is small (and $p \ne 0$).  In this case, $\mathcal{L}$ should be a subset of an arithmetic progression (with small gaps) since, for $\ell \in \mathcal{L}$, we know $u - \ell p$ is approximately $n_{\ell}t$ so the relative primeness of $t$ and $p$ yield the conclusion.}

\textit{Case:} $\displaystyle t \leq \frac{b_{\alpha_{0}-1}}{4}$ and $|p| > 0$
\begin{center}
\begin{minipage}{0.95\textwidth}
For $\ell \in \mathcal{L}$, \quad
$\displaystyle
\left|\frac{u-\ell p}{t} - n_{\ell}\right| < \frac{2}{b_{\alpha_{0}-1}} + \frac{1}{2t} \leq \frac{1}{2t} + \frac{1}{2t} = \frac{1}{t}
$
\quad so $|u - \ell p - n_{\ell} t| < 1$.  As those are integers, then $u - \ell p = n_{\ell} t$.  
Let $\ell_{0}$ be the minimal element of $\mathcal{L}$.  Then $(\ell - \ell_{0})p = (u - \ell_{0}p) - (u - \ell p) = (n_{\ell_{0}} - n_{\ell})t$.  As $p$ and $t$ are relatively prime, then $t$ must divide $\ell - \ell_{0}$ so
every $\ell \in \mathcal{L}$ is of the form $\ell = \ell_{0} + ti$ for some $0 \leq i < L$.
Also $\displaystyle |ty - p\tilde{h}_{\alpha_{0}}| < \frac{t}{Lb_{\alpha_{0}-1}}\tilde{h}_{\alpha_{0}} \leq \frac{1}{4L}\tilde{h}_{\alpha_{0}}$ meaning  $|ity - ip\tilde{h}_{\alpha_{0}}| \leq \frac{1}{4}\tilde{h}_{\alpha_{0}}$ for $0 \leq i < L$ so $\mathcal{L}$ has the structure of the second conclusion.
 \end{minipage}
 \end{center}

\textit{Case:} $\displaystyle\frac{b_{\alpha_{0}-1}}{4} < t \leq L$ and $|p| > 0$
\begin{center}
\begin{minipage}{0.95\textwidth}
\textit{In this case $y \approx \frac{p}{t}\tilde{h}_{\alpha_{0}} \gtrapprox \frac{p}{b_{\alpha_{0}}}\tilde{h}_{\alpha_{0}} \approx pa_{\alpha_{0}-1}\tilde{h}_{\alpha_{0}-1}$ so we can deduce mixing from having already established that (most) such times are mixing (and that most times in the progression are of that form).}\\

For $\ell \in \mathcal{L}$, $\displaystyle\left|\frac{u-\ell p}{t} - n_{\ell}\right| < \frac{2}{b_{\alpha_{0}-1}} + \frac{1}{2t} < \frac{4}{b_{\alpha_{0}-1}}$ so $|u - \ell p - n_{\ell}t| < \frac{4}{b_{\alpha_{0}-1}}t$ which means that $u - \ell p~(\mathrm{mod}~t) < \frac{4}{b_{\alpha_{0}-1}}t$ or $u - \ell p~(\mathrm{mod}~t)  > \left(1 - \frac{4}{b_{\alpha_{0}-1}}\right)t$ whenever $\ell \in \mathcal{L}$.

Since $p$ and $t$ are relatively prime, the map $\mathbb{Z} \to \mathbb{Z}/t\mathbb{Z}$ given by $z \mapsto zp~(\mathrm{mod}~t)$ is cyclic and onto.  So at most $\frac{8}{b_{\alpha_{0}-1}}t$ choices of $1 \leq \ell \leq t$ can be in $\mathcal{L}$.  Likewise, for any range of $t$ values of $\ell$, at most $\frac{8}{b_{\alpha_{0}-1}}t$ choices of $\ell$ (out of the $t$ possible) can be in $\mathcal{L}$.  Therefore $|\mathcal{L}| \leq \frac{8}{b_{\alpha_{0}-1}}t \left\lceil \frac{L}{t} \right\rceil$.
Then
\[
|\mathcal{L}| \leq \frac{8}{b_{\alpha_{0}-1}}t\left\lceil\frac{L}{t}\right\rceil \leq \frac{8}{b_{\alpha_{0}-1}}(L + t) \leq \frac{8}{b_{\alpha_{0}-1}}(L + L) < 4\epsilon\epsilon_{0}L
\]
Therefore, as $|\muB(T^{q-\ell k}B)| < \epsilon$ for $\ell \notin \mathcal{L}$, we have  $\left|\frac{1}{L}\sum_{\ell=0}^{L-1}\delta_{\ell}\muB(T^{q-\ell k}B)\right| + \frac{1}{L}\sum_{\ell=0}^{L-1}(1 - \delta_{\ell})\epsilon < \frac{1}{L}|\mathcal{L}| + \epsilon + \epsilon < 4\epsilon\epsilon_{0} + 2\epsilon$.
 \end{minipage}
 \end{center}

\textit{Case:} $\displaystyle|p| < \frac{2t}{b_{\alpha_{0}-1}\epsilon_{0}L}$

\begin{center}
\begin{minipage}{0.95\textwidth}
\textit{In this case, $p$ being very small means $y$ is very small so $\{ \ell : 0 \leq \ell < L \}$ already satisfies the claimed structure (and $\mathcal{L}$ is a subset of that).}

For $0 \leq \ell < L$, then
\quad $\displaystyle
|\ell p| < \frac{2t}{b_{\alpha_{0}-1}\epsilon_{0}} < \frac{2t}{4\epsilon^{-1}\epsilon_{0}^{-1}\epsilon_{0}} = 2\epsilon t
$ \quad so 
\[
\left|\frac{\ell y}{\tilde{h}_{\alpha_{0}}}\right| \leq \ell \left|\frac{y}{\tilde{h}_{\alpha_{0}}} - \frac{p}{t}\right| + \ell \left|\frac{p}{t}\right|
 \leq \frac{\ell}{b_{\alpha_{0}-1}L} + \left|\frac{\ell p}{t}\right| <
 \frac{1}{b_{\alpha_{0}-1}} + 2\epsilon < 3\epsilon
 \]
 For all $0 \leq i < L$, then
 $
 |i \cdot 1 \cdot y - i \cdot 0 \cdot \tilde{h}_{\alpha_{0}}|
 = |iy| < 3\epsilon \tilde{h}_{\alpha_{0}} < \frac{1}{3}\tilde{h}_{\alpha_{0}}$
 Therefore the second possible conclusion of the lemma holds with $L^{\prime} = L$ and $\ell_{0} = 0$ and $t = 1$ and $p = 0$.
 \end{minipage}
 \end{center}
 
From here on, we may then assume that \quad\quad
$
L < t \quad\quad\text{and}\quad\quad |p| \geq \frac{2t}{b_{\alpha_{0}-1}\epsilon_{0}L}
$

\textit{Our goal now is to apply the (generalized) Blum-Hanson trick combined with the Block Lemma.  The complication is that, unlike for staircases, there is no single choice of the power to use in the Block Lemma to reach mixing times but there are at most $1/\epsilon$ possible values.}

For $0 \leq j < \left\lceil \epsilon^{-1} \right\rceil$, define \quad\quad
$
p_{0} = p~(\mathrm{mod}~t) \quad\quad\text{and}\quad\quad p_{j+1} = \left\lceil\frac{t}{p_{j}}\right\rceil p_{j}~(\mathrm{mod}~t)
$

Suppose that $\epsilon t \leq p_{j}$ and $p_{j+1} \leq p_{j} - \epsilon t$ for all $0 \leq j < \left\lceil \epsilon^{-1} \right\rceil$.  Since $p_{j} \geq \epsilon t$, $\ceiling{\frac{t}{p_{j}}} \leq \ceiling{\epsilon^{-1}}$.  Then $\ceiling{\frac{t}{p_{0}}} \ceiling{\frac{t}{p_{1}}} \cdots \ceiling{\frac{t}{p_{\ceiling{\epsilon^{-1}}-1}}} \leq \ceiling{\epsilon^{-1}}^{\ceiling{\epsilon^{-1}}} < L < t$.  Therefore $p_{j+1} = \ceiling{\frac{t}{p_{0}}} \ceiling{\frac{t}{p_{1}}} \cdots \ceiling{\frac{t}{p_{j}}}p_{0}$ for $ j < \ceiling{\epsilon^{-1}}$.

Since $p_{j+1} \leq p_{j} - \epsilon t$, then $p_{j+1} \leq p_{j} - \epsilon t \leq p_{j-1} - 2\epsilon t \leq \ldots \leq p - j\epsilon t$.  Then $p_{\ceiling{\epsilon^{-1}}} \leq p - \ceiling{\epsilon^{-1}}\epsilon t \leq p - t$.  But then $\epsilon t \leq p_{0} \leq p_{\ceiling{\epsilon^{-1}}} \leq p - t < 0$ is a contradiction.
So 
there exists $0 \leq m < \ceiling{\epsilon^{-1}}$ such that $0 < p_{m}$ (since $p_{0} \ne 0$ as $p$ and $t$ are relatively prime) and either $p_{m} < \epsilon t$ or $p_{m} > p_{m-1} - \epsilon t$.

\textit{We now establish that in the bulk of the cases remaining, we can apply the generalized Blum-Hanson trick using $p_{m}$ to determine the power used in the Block Lemma.}

\textit{Case:} $\displaystyle \frac{2t}{b_{\alpha_{0}-1}} \leq p_{m} < \epsilon t$
\begin{center}
\begin{minipage}{0.95\textwidth}
Set $g = \ceiling{\frac{t}{p_{m-1}}}\cdots\ceiling{\frac{t}{p_{0}}}$.  We may assume $m$ is the minimal choice such that $p_{m} < \epsilon t$ so $p_{j} \geq \epsilon t$ for all $0 \leq j < m$.  If $\frac{t}{p_{j}} > \epsilon^{-1}$ then $t > \epsilon^{-1}p_{j} \geq \epsilon^{-1}\epsilon t = t$ is a contradiction, so $\frac{t}{p_{j}} \leq \epsilon^{-1}$ for all $0 \leq j < m$.  Then $g \leq \ceiling{\epsilon^{-1}}^{m} \leq \ceiling{\epsilon^{-1}}^{\ceiling{\epsilon^{-1}}-1}$.  

For $1 \leq i < \ceiling{(2\epsilon)^{-1}}$, then
$\displaystyle
\frac{2t}{b_{\alpha_{0}-1}} \leq i p_{m} < \ceiling{(2\epsilon)^{-1}}\epsilon t < \frac{1}{2}t + \epsilon t
$ so $igp~(\mathrm{mod}~t) \in \left[\frac{2}{b_{\alpha_{0}-1}}t, (\frac{1}{2} + \epsilon)t\right)$.  Then there exists $n_{i},w_{i} \in \mathbb{Z}$ so that $igp = nt + w$ and $w \in \left[\frac{2}{b_{\alpha_{0}-1}}t,(\frac{1}{2}+\epsilon)t\right)$.

Since
$\displaystyle
\left|\frac{i gy}{\tilde{h}_{\alpha_{0}}} - \frac{i g p}{t}\right| < \frac{i g}{L b_{\alpha_{0}-1}} < \frac{2\ceiling{\epsilon^{-1}}^{\ceiling{\epsilon^{-1}}}}{L b_{\alpha_{0}-1}} < \frac{1}{b_{\alpha_{0}-1}}
$, 
we have $\displaystyle
\left|igy - n_{i}\tilde{h}_{\alpha_{0}} - \frac{w_{i}}{t}\tilde{h}_{\alpha_{0}}\right| < \frac{1}{b_{\alpha_{0}-1}}
$
so
\[
\left|igy - n_{i}\tilde{h}_{\alpha_{0}}\right| > \frac{w_{i}}{t}\tilde{h}_{\alpha_{0}} - \frac{1}{b_{\alpha_{0}-1}}\tilde{h}_{\alpha_{0}} \geq \frac{2}{b_{\alpha_{0}-1}}\tilde{h}_{\alpha_{0}} - \frac{1}{b_{\alpha_{0}-1}}\tilde{h}_{\alpha_{0}}
= \frac{1}{b_{\alpha_{0}-1}}\tilde{h}_{\alpha_{0}} > a_{\alpha_{0}-1}\tilde{h}_{\alpha_{0}-1}
\]
and likewise $|igy - n_{i}\tilde{h}_{\alpha_{0}}| \leq (\frac{1}{2} + \epsilon)\tilde{h}_{\alpha_{0}} + \frac{1}{b_{\alpha_{0}-1}}\tilde{h}_{\alpha_{0}} < \tilde{h}_{\alpha_{0}} - a_{\alpha_{0}-1}\tilde{h}_{\alpha_{0}-1}$.

Since $igk = i g z \tilde{h}_{\alpha_{0}} + i gy$, then $igk$ is of the form $v_{i}\tilde{h}_{\alpha_{0}} + \tilde{y}_{i}$ for $v_{i},\tilde{y}_{i} \in \mathbb{Z}$ where $a_{\alpha_{0}-1}\tilde{h}_{\alpha_{0}-1} \leq \tilde{y}_{i} < \tilde{h}_{\alpha_{0}} - a_{\alpha_{0}-1}\tilde{h}_{\alpha_{0}-1}$.  Therefore
$|\muB(T^{i g k}B)| \leq \widetilde{M}_{B,\alpha_{0}} < \epsilon$ (or instead $|\muB(T^{i g k}B)| < M_{B,\alpha_{0}-1} < \epsilon$ when $v_{i} \geq a_{\alpha_{0}}$).

Since
 $\displaystyle
g\ceiling{(2\epsilon)^{-1}} < 2\ceiling{\epsilon^{-1}}^{\ceiling{\epsilon^{-1}}} < \ceiling{\epsilon^{-1}}^{-1}L < \epsilon L
$,
 by Lemma \ref{deltaBL} (with $p = g$ and $Q = \lceil(2\epsilon)^{-1}\rceil$),
\[
\left|\frac{1}{L}\sum_{\ell=0}^{L-1}\muB(T^{q-\ell k}B)\right| + \frac{1}{L}\sum_{\ell=0}^{L-1}(1 - \delta_{\ell})\epsilon 
< (2\epsilon)^{1/2} + \epsilon + \epsilon < 4\epsilon^{1/2}
\]
 \end{minipage}
 \end{center}
 
\textit{Case:} $\displaystyle \frac{2t}{b_{\alpha_{0}-1}\epsilon_{0}L} \leq p_{m} < \frac{2t}{b_{\alpha_{0}-1}}$
\begin{center}
\begin{minipage}{0.95\textwidth}
Let $g \in \mathbb{N}$ minimal so that $gp_{m} \geq \frac{2t}{b_{\alpha_{0}-1}}$.  Then $gp_{m} < \frac{4t}{b_{\alpha_{0}-1}}$ (if not then $(g-1)p_{m} > \frac{4t}{b_{\alpha_{0}-1}} - \frac{2t}{b_{\alpha_{0}-1}}$ contradicts that $g$ is minimal).  
For $1 \leq i < \ceiling{\epsilon^{-1}}$, 
\[
\frac{2t}{b_{\alpha_{0}-1}} \leq igp_{m} < 4\ceiling{\epsilon^{-1}}\frac{t}{b_{\alpha_{0}-1}} < \frac{1}{2}t
\quad\quad\text{meaning}\quad\quad
ig\left\lceil\frac{t}{p_{m-1}}\right\rceil \cdots \left\lceil \frac{t}{p_{0}}\right\rceil p~(\mathrm{mod}~t) \in \left[\frac{2t}{b_{\alpha_{0}-1}}, \frac{t}{2}\right)
\]
Set $g_{0} = g\ceiling{\frac{t}{p_{m-1}}}\cdots \ceiling{\frac{t}{p_{0}}}$.  By the same reasoning as in the $\frac{2t}{b_{\alpha_{0}-1}} \leq p_{m} < \epsilon t$ case, $\ceiling{\frac{t}{p_{j}}} \leq \ceiling{\epsilon^{-1}}$ for all $0 \leq j < m$.  Then
\[
g_{0} = g\left\lceil\frac{t}{p_{m-1}}\right\rceil \cdots \left\lceil \frac{t}{p_{0}}\right\rceil < 
\frac{4t}{b_{\alpha_{0}}} \frac{1}{p_{m}}  \ceiling{\epsilon^{-1}}^{m}
\leq \frac{4t}{b_{\alpha_{0}}} \frac{b_{\alpha_{0}}\epsilon_{0}L}{2t}  \ceiling{\epsilon^{-1}}^{\ceiling{\epsilon^{-1}}}
= 2\epsilon_{0}L \ceiling{\epsilon^{-1}}^{\ceiling{\epsilon^{-1}}} < \epsilon L
\]
Since $ig_{0}p~(\mathrm{mod}~t) \in \left[\frac{2t}{b_{\alpha_{0}-1}}, \frac{t}{2}\right)$, Lemma \ref{deltaBL} then gives the first conclusion.
 \end{minipage}
 \end{center}

\textit{Case:} $\displaystyle 0 < p_{m} < \frac{2t}{b_{\alpha_{0}-1}\epsilon_{0}L}$ and $t < \epsilon_{0}b_{\alpha_{0}-1}L$
\begin{center}
\begin{minipage}{0.95\textwidth}
Set $g = \ceiling{\frac{2t}{b_{\alpha_{0}-1}p_{m}}}$ so $g < \frac{2t}{b_{\alpha_{0}-1}}+1 <\frac{2\epsilon_{0}b_{\alpha_{0}-1}L}{b_{\alpha_{0}-1}}+1 = 2\epsilon_{0}L + 1$.  For $1 \leq i < \ceiling{\epsilon^{-1}}$, then $\frac{2t}{b_{\alpha_{0}-1}} \leq igp_{m} < \ceiling{\epsilon^{-1}}(2\epsilon_{0}L + 1)\frac{2t}{b_{\alpha_{0}-1}\epsilon_{0}L} < \frac{t}{2}$.  Since $g\ceiling{\frac{t}{p_{m-1}}}\cdots\ceiling{\frac{t}{p_{0}}} < (2\epsilon_{0}L + 1)\ceiling{\epsilon^{-1}}^{\ceiling{\epsilon^{-1}}-1} < \ceiling{\epsilon^{-1}}^{-2}L$, then, using the same reasoning as in the above cases, Lemma \ref{deltaBL} gives the first conclusion.
 \end{minipage}
 \end{center}
 
 \textit{Unfortunately, there is one remaining case for which no choice of $p$ will allow us to apply the Blum-Hanson trick.  However, in this case, we can establish that the potentially bad times $\mathcal{L}$ themselves have arithmetic-like structure with control over the gaps.}
 
\textit{Case:} $0 < p_{m} < \frac{2t}{b_{\alpha_{0}-1}\epsilon_{0}L}$ and $t \geq \epsilon_{0}b_{\alpha_{0}-1}L$
\begin{center}
\begin{minipage}{0.95\textwidth}
Since $|p| \geq \frac{2t}{b_{\alpha_{0}-1}\epsilon_{0}L}$ and $p_{m} < \frac{2t}{b_{\alpha_{0}-1}\epsilon_{0}L}$, it follows that $m > 0$.  Set $g = \ceiling{\frac{t}{p_{m-1}}}\cdots\ceiling{\frac{t}{p_{0}}} < \ceiling{\epsilon^{-1}}^{\ceiling{\epsilon^{-1}}-1} = \frac{1}{2}\epsilon_{0}^{-1}\ceiling{\epsilon^{-1}}^{-2} < \epsilon^{2}\epsilon_{0}^{-1}$.  Since $gp~(\mathrm{mod}~t) = p_{m}$, then there exists an integer $v > 0$ (since $m > 0$) such that $gp = vt + p_{m}$.  Then $|gp - vt| = p_{m} < \frac{2}{b_{\alpha_{0}-1}\epsilon_{0}L}$.

For $\ell \in \mathcal{L}$, recall that there exists $n_{\ell} \in \mathbb{Z}$ such that $|u - \ell p - n_{\ell}t| < \frac{1}{2} + \frac{2t}{b_{\alpha_{0}-1}}$.  Let $\ell_{0}$ minimal such that $\ell_{0} \in \mathcal{L}$.  Since $|p| \geq \frac{2t}{b_{\alpha_{0}-1}\epsilon_{0}L}$, then $\frac{1}{b_{\alpha_{0}}} \leq \frac{1}{2}\epsilon_{0} \frac{|p|L}{t}$ so
\[
\left|n_{\ell} - n_{\ell_{0}}\right|
\leq \left|\frac{u - \ell p}{t} - \frac{u - \ell_{0}p}{t}\right| + \frac{1}{t} + \frac{4}{b_{\alpha_{0}-1}} < \frac{(\ell - \ell_{0})|p|}{t} + \frac{5}{b_{\alpha_{0}-1}}
< \frac{L|p|}{t} + \frac{5}{2}\epsilon_{0}\frac{|p|L}{t}
= \left(1 + \frac{5\epsilon_{0}}{2}\right)\frac{|p|L}{t}
\]
Then, since $b_{\alpha_{0}-1}\epsilon_{0} > 4\epsilon^{-1}$,
\begin{align*}
\left|(n_{\ell} - n_{\ell_{0}})vt - (n_{\ell} - n_{\ell_{0}})gp\right| &< |n_{\ell} - n_{\ell_{0}}| \frac{2t}{b_{\alpha_{0}-1}\epsilon_{0}L}
< \left(1 + \frac{5\epsilon_{0}}{2}\right) \frac{2}{b_{\alpha_{0}-1}\epsilon_{0}} |p| \\
&= \frac{2}{b_{\alpha_{0}-1}\epsilon_{0}} |p| + \frac{5}{b_{\alpha_{0}-1}}|p|
< \frac{\epsilon}{2}|p| + \frac{5\epsilon\epsilon_{0}}{4}|p|
\end{align*}
and
\begin{align*}
\left|(n_{\ell}-n_{\ell_{0}})vt - v(\ell - \ell_{0})p\right| &= |v| \left|(n_{\ell} - (u - \ell p)) - (n_{\ell_{0}} - (u - \ell_{0}p))\right|
< |v|\left(1 + \frac{4t}{b_{\alpha_{0}-1}}\right) \\ &< \frac{5|v|t}{b_{\alpha_{0}-1}} 
< \frac{5g|p|}{b_{\alpha_{0}-1}} < \frac{5\ceiling{\epsilon^{-1}}^{\ceiling{\epsilon^{-1}}-1}}{b_{\alpha_{0}-1}}|p|
< \frac{5\epsilon^{2}\epsilon_{0}^{-1}}{b_{\alpha_{0}-1}}|p| < \frac{5\epsilon^{3}}{4}|p|
\end{align*}
Therefore 
$\displaystyle
\left|(n_{\ell}-n_{\ell_{0}})gp - v(\ell_{0} - \ell)p\right| <
\left(\frac{\epsilon}{2} + \frac{5\epsilon\epsilon_{0}}{4} + \frac{5\epsilon^{3}}{4}\right)|p|
$  so $|(n_{\ell} - n_{\ell_{0}})g - v(\ell_{0} - \ell)| < 1$.  As these are integers, then $(n_{\ell} - n_{\ell_{0}})g = v(\ell_{0} - \ell)$ for all $\ell \in \mathcal{L}$.  Set $g^{\star} = \frac{g}{gcd(g,v)}$ and $v^{\star} = \frac{v}{gcd(g,v)}$.  Then $g^{\star}$ divides $\ell_{0} - \ell$ for all $\ell \in \mathcal{L}$ so every $\ell \in \mathcal{L}$ is of the form $\ell = \ell_{0} + ig^{\star}$.
Also for $0 \leq i < L$,
\begin{align*}
|ig^{\star}y - iv^{\star}\tilde{h}_{\alpha_{0}}| &= \frac{i\tilde{h}_{\alpha_{0}}}{gcd(g,v)}\left|\frac{gy}{\tilde{h}_{\alpha_{0}}} - v\right|
< \frac{i\tilde{h}_{\alpha_{0}}}{gcd(g,v)}\left(\left|\frac{gp}{t} - v\right| + \frac{g}{b_{\alpha_{0}-1}L}\right)  \\
&< L\tilde{h}_{\alpha_{0}} \left(\frac{2}{b_{\alpha_{0}-1}\epsilon_{0}L} + \frac{\epsilon^{2}\epsilon_{0}^{-1}}{b_{\alpha_{0}-1}L}\right)
= \frac{2 + \epsilon^{2}}{\epsilon_{0}b_{\alpha_{0}-1}} \tilde{h}_{\alpha_{0}}
< \frac{2 + \epsilon^{2}}{\epsilon_{0}4\epsilon_{0}^{-1}\epsilon^{-1}}\tilde{h}_{\alpha_{0}}
< (2\epsilon + \epsilon^{3})\tilde{h}_{\alpha_{0}}
\end{align*}
 so the second possible conclusion holds with $g^{\star}$ for $t$ and $v^{\star}$ for $p$.
 \end{minipage}
 \end{center}
 
 The above cases cover all possibilities where $p_{m} < \epsilon t$.  Therefore we may assume from here on that
 \[
p_{m} > p_{m-1} - \epsilon t
\]
Set
$\displaystyle
p^{\star} = \left\lfloor \frac{t}{p_{m-1}} \right\rfloor p_{m-1} = \left\lceil\frac{t}{p_{m-1}}\right\rceil p_{m-1} - p_{m-1} = p_{m} + t - p_{m-1} > t - \epsilon t
$.

\textit{The remaining cases follow along similar lines as the cases where $p_{n} < \epsilon t$.}

Note that $p_{j} \geq \epsilon t$ for all $j < m$ by the minimality of $m$ so we again have that $\frac{t}{p_{j}} < \epsilon^{-1}$ for all $j < m$.

\textit{Case:} $\displaystyle t - \epsilon t < p^{\star} \leq t - \frac{2t}{b_{\alpha_{0}-1}}$
\begin{center}
\begin{minipage}{0.95\textwidth}
For $1 \leq i < \ceiling{(2\epsilon)^{-1}}$, $\frac{1}{2}t - \epsilon t < ip^{\star}~(\mathrm{mod}~t) \leq t - \frac{2t}{b_{\alpha_{0}-1}}$.  Then $i\floor{\frac{t}{p_{m-1}}}\ceiling{\frac{t}{p_{m-2}}} \cdots \ceiling{\frac{t}{p_{0}}} p~(\mathrm{mod}~t)$ is nonzero and at least $\frac{2t}{b_{\alpha_{0}-1}}$ away from every multiple of $t$.  Since $\ceiling{\epsilon^{-1}}\floor{\frac{t}{p_{m-1}}}\cdots \ceiling{\frac{t}{p_{0}}} < \ceiling{\epsilon^{-1}}^{\ceiling{\epsilon^{-1}}}\epsilon_{0}L < \epsilon L$, as above Lemma \ref{deltaBL} gives the claim.
 \end{minipage}
 \end{center}

\textit{Case:} $\displaystyle t - \frac{2t}{b_{\alpha_{0}-1}} < p^{\star} \leq t - \frac{2t}{b_{\alpha_{0}-1}\epsilon_{0}L}$
\begin{center}
\begin{minipage}{0.95\textwidth}
Let $g \in \mathbb{N}$ minimal such that $gp^{\star}~(\mathrm{mod}~t) \leq t - \frac{2t}{b_{\alpha_{0}-1}}$.  As in the case where $\frac{2t}{b_{\alpha_{0}-1}\epsilon_{0} L} \leq p_{m} < \frac{2t}{b_{\alpha_{0}-1}}$, $g < 2\epsilon_{0}L$ and then similar reasoning as there using Lemma \ref{deltaBL} gives the claim.
 \end{minipage}
 \end{center}

\textit{Case:} $\displaystyle t - \frac{2t}{b_{\alpha_{0}-1}\epsilon_{0}L} < p^{\star} < t$ and $\displaystyle t \leq \epsilon_{0} b_{\alpha_{0}}L$
\begin{center}
\begin{minipage}{0.95\textwidth}
Set $g = \ceiling{\frac{2t}{b_{\alpha_{0}-1}(t - p^{\star})}} < 2\epsilon_{0}L + 1$.  Then $igp^{\star}~(\mathrm{mod}~t) < t - \frac{2t}{b_{\alpha_{0}-1}(t - p^{\star})}(t - p^{\star}) = t - \frac{2t}{b_{\alpha_{0}-1}}$ and $igp^{\star}~(\mathrm{mod}~t) > t - \ceiling{\epsilon^{-1}}(2\epsilon_{0}L + 1)\frac{2}{b_{\alpha_{0}-1}\epsilon_{0}L}$ so again similar reasoning gives the claim.
 \end{minipage}
 \end{center}

\textit{Case:} $\displaystyle t - \frac{2t}{b_{\alpha_{0}-1}\epsilon_{0}L} < p^{\star} < t$ and $\displaystyle t \geq \epsilon_{0}b_{\alpha_{0}-1}L$ and $|p| \geq \frac{2t}{b_{\alpha_{0}-1}\epsilon_{0}L}$
\begin{center}
\begin{minipage}{0.95\textwidth}
Set $g = \floor{\frac{t}{p_{m-1}}}\ceiling{\frac{t}{p_{m-2}}}\cdots\ceiling{\frac{t}{p_{0}}} < \ceiling{\epsilon^{-1}}^{\ceiling{\epsilon^{-1}}-1}$.  Then $p^{\star} = gp~(\mathrm{mod}~t)$.  Here, as above, $m \ne 0$ so there exists $v \ne 0$ such that $vt - \frac{2t}{b_{\alpha_{0}-1}\epsilon_{0}L} < gp < vt$ and the same argument as in the $0 < p_{m} < \frac{2t}{b_{\alpha_{0}-1}\epsilon_{0}}$ case shows that the second possible conclusion holds.
 \end{minipage}
 \end{center}
 
Therefore the claim is proved as all cases have been covered.
\end{proof}

For clarity of exposition, we introduce notation for the amount of the leftmost subcolumns which will not be pushed through the top of the next tower upon application of $T^{q - \ell k}$.  However, we will declare it to be $0$ if $q - \ell k$ is already known to be a mixing time.

\begin{notation}
 For $\ell,q,k \in \mathbb{Z}$ and $n > 0$ such that $|q - \ell k| < \tilde{h}_{n+1}$, let $d$ be the unique integer such that $|(q-\ell k)- d\tilde{h}_{n}| \leq \frac{1}{2}\tilde{h}_{n}$ and
define
\[
\gam{n}{q}{k}{\ell} = \left\{ \begin{array}{ll} \frac{a_{n} - |d|}{a_{n}} &\quad\text{if $(b_{n} \leq |d| < a_{n}$ or $d=0)$ and $|(q - \ell k) - d\tilde{h}_{n}| < a_{n-1}\tilde{h}_{n-1}$} \\ 0 &\quad\text{otherwise} \end{array}\right.
\]
\end{notation}

Our next lemma is an application of Lemma \ref{bigsmally} to precisely the set of potentially bad times $\mathcal{L}$ in the previous lemma.  The critical point here is that when we split $T^{q - k(\ell_{0}+ti)}B$ into a $\frac{1}{b_{\alpha_{0}}}\sum_{j=0}^{b_{\alpha_{0}}}$ average, the resulting new $q$ and $k$ for each $j$ is the same over all $i$.  This relies in a crucial way on the bounded gaps and arithmetic structure established in the preceding lemma.  Since Lemma \ref{bigsmally} does produce an error term, the $\gamma$ terms are necessary bookkeeping devices to ensure the total error is bounded.

\begin{lemma}\label{thirdB}
Let $\epsilon > 0$.  Let $q,k,t,p,\ell_{0},L^{\prime} \in \mathbb{Z}$ with $\ell_{0},t,L^{\prime} > 0$.
Set $\alpha_{0} = \max \{ \alpha(q - k(\ell_{0} + it)) : 0 \leq i < L^{\prime} \}$.

Assume that $b_{\alpha_{0}-1} > 2 \epsilon^{-1}$ and $\sup_{m \geq b_{\alpha_{0}-1}} \left(\int \left| \frac{1}{m}\sum_{i=0}^{m-1}\chi_{B}\circ T^{-i}\right|~d\mu + \frac{2}{m}\right) < \frac{\epsilon}{3}$.

Let $y,z \in \mathbb{Z}$ such that $k = z\tilde{h}_{\alpha_{0}} + y$ with $|y| \leq \frac{1}{2}\tilde{h}_{\alpha_{0}}$.

Assume, for $0 \leq i < L^{\prime}$, that $\left|ity - ip\tilde{h}_{\alpha_{0}}\right| < \frac{1}{3}\tilde{h}_{\alpha_{0}}$.

For each $i$, let $q - k(\ell_{0} + ti) = k_{i}\tilde{h}_{\alpha_{0}} + y_{i}$ with $|y_{i}| \leq \frac{1}{2}\tilde{h}_{\alpha_{0}}$.  Assume that $|k_{0}| < a_{\alpha_{0}}$ and $|y_{0}| < a_{\alpha_{0}-1}\tilde{h}_{\alpha_{0}-1}$.

Then for each $0 \leq j < b_{\alpha_{0}}$ such that $|k_{i}| < a_{\alpha_{0}}$ and $|y_{i}| < a_{\alpha_{0}-1}\tilde{h}_{\alpha_{0}-1}$, there exists $q_{j},h_{j} \in \mathbb{Z}$ such that $\alpha_{j} = \max \{ \alpha(q_{j} - h_{j}i) : 0 \leq i < L^{\prime} \} < \alpha_{0}$ and for every $0 \leq i < L^{\prime}$,
\begin{align*}
&\left|\muB(T^{q - k(\ell_{0} + it)}B) - \frac{r_{\alpha_{0}}}{r_{\alpha_{0}}+1} \gam{\alpha_{0}}{q}{k}{\ell_{0} + it} \frac{1}{b_{\alpha_{0}}} \sum_{j=0}^{b_{\alpha_{0}}-1} \muB(T^{q_{j} - h_{j}i}B)\right| \\
&\quad\quad \leq 
\left(1 - \gam{\alpha_{0}}{q}{k}{\ell_{0} + it} \right) \epsilon + 
\gam{\alpha_{0}}{q}{k}{\ell_{0} + it} \frac{1}{b_{\alpha_{0}}} \sum_{j=0}^{b_{\alpha_{0}}-1}\left(1 - \gam{\alpha_{j}}{q_{j}}{h_{j}}{i}\right)\epsilon + \tau_{\alpha_{0}} + \frac{1}{a_{\alpha_{0}-1}b_{\alpha_{0}-1}} 
\end{align*}
\end{lemma}
\begin{proof}
\textit{We first use the bounded gaps property to write a concrete expression for $k_{i}$ and $y_{i}$, specifcally to write $q - k(\ell_{0} + ti)$ as a multiple of $\tilde{h}_{\alpha_{0}}$ plus a remainder term where the remainder term and multiple term are linear functions of $i$.}

Since $|y_{0} - ity + ip\tilde{h}_{\alpha_{0}}| < a_{\alpha_{0}-1}\tilde{h}_{\alpha_{0}-1} + \frac{1}{3}\tilde{h}_{\alpha_{0}} < \frac{1}{2}\tilde{h}_{\alpha_{0}}$ and since
\begin{align*}
q - (\ell_{0} + it)k &= q - \ell_{0}k - it(z\tilde{h}_{\alpha_{0}} + y)
= k_{0}\tilde{h}_{\alpha_{0}} + y_{0} - itz\tilde{h}_{\alpha_{0}} - ity \\
&= (k_{0} - itz)\tilde{h}_{\alpha_{0}} + y_{0} - ity
= (k_{0} - itz - ip)\tilde{h}_{\alpha_{0}} + y_{0} - ity + ip\tilde{h}_{\alpha_{0}}
\end{align*}
then \quad\quad $k_{i} = k_{0} - itz - ip$ \quad and \quad $y_{i} = y_{0} - ity + ip\tilde{h}_{\alpha_{0}}$.

\textit{From the linear expressions, we can define what the new powers $q_{j}$ and $h_{j}$ are and show that they are independent of the choice of $i$.}

For $0 \leq j < b_{\alpha_{0}}$, set \quad $q_{j} = y_{0} - jk_{0}$ \quad and \quad $h_{j} = ty - p\tilde{h}_{\alpha_{0}} - jt - jp$ \quad and then
\[
y_{i} - k_{i}j = y_{0} - ity + ip\tilde{h}_{\alpha_{0}} - (k_{0} - itz - ip)j
= (y_{0} - k_{0}j) - (ty - p\tilde{h}_{\alpha_{0}} - jtz - jp)i
= q_{j} - h_{j} i
\]
Since $|k_{i}| < a_{\alpha_{0}}$ and $j < b_{\alpha_{0}}$, then $|q_{j} - h_{j}i| = |y_{i} - k_{i}j| < a_{\alpha_{0}-1}\tilde{h}_{\alpha_{0}-1} + a_{\alpha_{0}}b_{\alpha_{0}} < \tilde{h}_{\alpha_{0}}$.  Therefore $\alpha_{j} \leq \alpha_{0} - 1$.

Consider first $i$ such that $k_{i} = 0$.  Then for any $j$,
$
q - (\ell_{0} + it)k = y_{i} = y_{i} - k_{i}j = q_{j} - h_{j}i
$ so
\[
\muB(T^{q - (\ell_{0}+it)} B) = \frac{1}{b_{\alpha_{0}}}\sum_{j=0}^{b_{\alpha_{0}}-1} \muB(T^{q_{j}-h_{j}i}B)
\]
Since $\gam{\alpha_{0}}{q}{k}{\ell_{0}+it} = 1$ (as $k_{i} = 0$),
\begin{align*}
\left|\muB(T^{q - k(\ell_{0} + it)}B) - \frac{r_{\alpha_{0}}}{r_{\alpha_{0}}+1} \gam{\alpha_{0}}{q}{k}{\ell_{0} + it} \frac{1}{b_{\alpha_{0}}} \sum_{j=0}^{b_{\alpha_{0}}-1} \muB(T^{q_{j} - h_{j}i}B)\right|
= 1 - \frac{r_{\alpha_{0}}}{r_{\alpha_{0}}+1} < \frac{1}{a_{\alpha_{0}}b_{\alpha_{0}}} < \frac{1}{a_{\alpha_{0}-1}b_{\alpha_{0}-1}}
\end{align*}
so the claim holds for $i$ such that $k_{i} = 0$.

From here on, assume $k_{i} \ne 0$.

\textit{Next we bound the error term that will come from Lemma \ref{bigsmally}.  For $j$ such that $\alpha(q_{j} - h_{j}i) \leq \alpha_{0} - 2$, this is straightforwardly bounded by $\frac{1}{a_{\alpha_{0}-1}b_{\alpha_{0}-1}}$ which we will later insist be summable.  The $\gamma$ terms appear exactly for $j$ such that $\alpha(q_{j}-h_{j}i) = \alpha_{0}-1$.}

For $i,j$ such that $\alpha(q_{j} - h_{j}i) < \alpha_{0} - 1$, we have
$\displaystyle
\frac{|q_{j} - h_{j}i|}{\tilde{h}_{\alpha_{0}}} < \frac{\tilde{h}_{\alpha_{0}-1}}{\tilde{h}_{\alpha_{0}}} < \frac{1}{a_{\alpha_{0}-1}b_{\alpha_{0}-1}}
$.

For $i,j$ such that $\alpha(q_{j} - h_{j}i) = \alpha_{0} - 1$, let $d_{i,j}$ be the unique integer such that $|q_{j} - h_{j}i - d_{i,j}\tilde{h}_{\alpha_{0}-1}| \leq \frac{1}{2}\tilde{h}_{\alpha_{0}-1}$.  Then $|d_{i,j}| \leq a_{\alpha_{0}-1}$ since $|q_{j}-h_{j}i| < a_{\alpha_{0}-1}\tilde{h}_{\alpha_{0}-1} + a_{\alpha_{0}}b_{\alpha_{0}}$ so $\gam{\alpha_{0}-1}{q_{j}}{h_{j}}{i} = \frac{a_{\alpha_{0}-1} - |d_{i,j}|}{a_{\alpha_{0}-1}}$.  Then
\[
\frac{|q_{j} - h_{j}i|}{\tilde{h}_{\alpha_{0}}} < \frac{(|d_{i,j}| + 1)\tilde{h}_{\alpha_{0}-1}}{\tilde{h}_{\alpha_{0}}}
< \frac{|d_{i,j}| + 1}{a_{\alpha_{0}-1}b_{\alpha_{0}-1}}
< \frac{2|d_{i,j}|}{a_{\alpha_{0}-1}b_{\alpha_{0}-1}} 
< \frac{|d_{i,j}|}{a_{\alpha_{0}-1}} \frac{\epsilon}{2}
=  \left(1 - \gam{\alpha_{0}-1}{q_{j}}{h_{j}}{i}\right) \frac{\epsilon}{2}
\]
Therefore, for every $i$ with $|k_{i}| < a_{\alpha_{0}}$ and $k_{i} \ne 0$,
\[
\frac{1}{b_{\alpha_{0}}}\sum_{j=0}^{b_{\alpha_{0}}-1} \frac{2|y_{i} - k_{i}j|}{\tilde{h}_{\alpha_{0}}}
\leq \frac{1}{b_{\alpha_{0}}}\sum_{j=0}^{b_{\alpha_{0}}-1} \left(1 - \gam{\alpha_{j}}{q_{j}}{h_{j}}{i}\right) \epsilon + \frac{1}{a_{\alpha_{0}-1}b_{\alpha_{0}-1}}
\]
By Lemma \ref{bigsmally},
\begin{align*}
&\left|\muB(T^{q-(\ell_{0}+ti)k}B) - \frac{a_{\alpha_{0}}-|k_{ii}|}{r_{\alpha_{0}}+1}\sum_{j=0}^{b_{\alpha_{0}}-1}\muB(T^{y_{i} - k_{i}j}B)\right| \\
&\quad\quad\quad\quad \leq 
 \frac{a_{\alpha_{0}}-|k_{i}|}{a_{\alpha_{0}}}\frac{1}{b_{\alpha_{0}}}\sum_{j=0}^{b_{\alpha_{0}}-1} \frac{2|y_{i} - k_{i}j|}{\tilde{h}_{\alpha_{0}}} 
 + \frac{|k_{}|}{a_{\alpha_{0}}}\epsilon
+ \tau_{\alpha_{0}}
\end{align*}
Since $\gam{\alpha_{0}}{q}{k}{\ell_{0}+it} = \frac{a_{\alpha_{0}} - |k_{i}|}{a_{\alpha_{0}}}$, then
\begin{align*}
&\left|\muB(T^{q-(\ell_{0}+ti)k}B) - \frac{r_{\alpha_{0}}}{r_{\alpha_{0}}+1}\gam{\alpha_{0}}{q}{k}{\ell_{0}+it}\sum_{j=0}^{b_{\alpha_{0}}-1}\muB(T^{y_{i} - k_{i}j}B)\right| \\
&\quad\quad\quad\quad \leq 
\gam{\alpha_{0}}{q}{k}{\ell_{0}+it}\left(\frac{1}{b_{\alpha_{0}}}\sum_{j=0}^{b_{\alpha_{0}}-1} \left(1 - \gam{\alpha_{j}}{q_{j}}{h_{j}}{i}\right) \epsilon + \frac{1}{a_{\alpha_{0}}b_{\alpha_{0}}}\right)
 + \left(1 - \gam{\alpha_{0}}{q}{k}{\ell_{0}+it}\right) \epsilon
+ \tau_{\alpha_{0}} \qedhere
\end{align*}
\end{proof}

Our final lemma addresses the situation when the potentially bad times themselves form a controlled arithmetic progression.  In essence, it says that such a progression of potentially bad times can themselves be shown to lead to a `weighted' average with (much) smaller gaps.  The reason for the weights becomes apparent in the proof: $1 - \gamma$ of the tower is already known to be mixing and $\gamma$ of the tower might potentially not be.  The proof is really just combining our previous lemmas without new ideas but performing necessary bookkeeping and sum rearranging.

\begin{lemma}\label{third}
Let $\epsilon > 0$ and set $\epsilon_{0} = (2\ceiling{\epsilon^{-1}}^{\ceiling{\epsilon^{-1}}+1})^{-1}$.  Let $L,k,q \in \mathbb{Z}$ with $L \geq \epsilon_{0}^{-1}$ and for each $0 \leq \ell < L$, let $0 \leq \delta_{\ell} \leq 1$.  

Let $\alpha_{0} = \max \{ \alpha(q - \ell k) : 0 \leq \ell < L \}$.
Assume that $\max(M_{B,\alpha_{0}}, M_{B,\alpha_{0}-1}, \widehat{M}_{B,\alpha_{0}}, \widetilde{M}_{B,\alpha_{0}}) < \epsilon$ and $b_{\alpha_{0}-1} > 4 \epsilon^{-1}\epsilon_{0}^{-1}$ and $\sup_{m \geq b_{\alpha_{0}-1}} \left(\int \left| \frac{1}{m}\sum_{i=0}^{m-1}\chi_{B}\circ T^{-i}\right|~d\mu + \frac{2}{m}\right) < \frac{\epsilon}{3}$.

Then either
\[
\left|\frac{1}{L}\sum_{\ell=0}^{L-1}\delta_{\ell} \muB(T^{q - \ell k}B) \right| + \frac{1}{L}\sum_{\ell=0}^{L-1}(1 - \delta_{\ell})\epsilon < 4\epsilon^{1/2}
\]
or there exist integers $0 < t , L^{\prime} < L$ and, for each $0 \leq j < b_{\alpha_{0}}$, there exist $h_{j}, q_{j} \in \mathbb{Z}$ such that $\alpha_{j} = \max\{\alpha(q_{j} - h_{j}\ell) : 0 \leq \ell < L^{\prime} \} < \alpha_{0}$ and 
\begin{align*}
&\left|\frac{1}{L}\sum_{\ell=0}^{L-1}\delta_{\ell} \muB(T^{q - \ell k}B) \right| + \frac{1}{L}\sum_{\ell=0}^{L-1}(1 - \delta_{\ell})\epsilon \\
&\quad\quad < \frac{1}{b_{\alpha_{0}}}\sum_{j=0}^{b_{\alpha_{0}}-1} \left( \left|\frac{1}{L}\sum_{\ell=0}^{L - 1} \bbone_{\ell < L^{\prime}}\delta_{\ell_{0}+\ell t} \gam{\alpha_{0}}{q}{k}{\ell_{0} + t \ell}
\muB(T^{q_{j}-h_{j}\ell}B)\right| + \frac{1}{L}\sum_{\ell=0}^{L - 1} \left(1 - \bbone_{\ell < L^{\prime}}\delta_{\ell_{0}+\ell t} \gam{\alpha_{0}}{q}{k}{\ell_{0} + t \ell}
\right)\epsilon \right) \\
&\quad\quad\quad\quad + \frac{1}{b_{\alpha_{0}}}\sum_{j=0}^{b_{\alpha_{0}}-1}\frac{1}{L}\sum_{\ell=0}^{L-1}\bbone_{\ell < L^{\prime}}\delta_{\ell_{0}+t\ell}\gam{\alpha_{0}}{q}{k}{\ell_{0} + t \ell}\left(1 - \gam{\alpha_{j}}{q_{j}}{h_{j}}{\ell}\right)\epsilon + \tau_{\alpha_{0}} + \frac{1}{a_{\alpha_{0}-1}b_{\alpha_{0}-1}}
\end{align*}
\end{lemma}
\begin{proof}
Write $q - \ell k = k_{\ell}\tilde{h}_{\alpha_{0}} + y_{\ell}$ with $|y_{\ell}| \leq \frac{1}{2}\tilde{h}_{\alpha_{0}}$.  Define
\[
\mathcal{L} = \left\{ 0 \leq \ell < L : \left(k_{\ell} = 0~\text{or}~b_{\alpha_{0}} \leq |k_{\ell}| < a_{\alpha_{0}}\right)~\text{and}~|y_{\ell}| < a_{\alpha_{0}-1}\tilde{h}_{\alpha_{0}-1} \right\}
\]

By Lemma \ref{thirdA}, we may assume from here on that there exists $p  \in \mathbb{Z}$, $t > 0$ and $0 \leq \ell_{0} < L^{\prime} \leq L$ such that $\mathcal{L} \subseteq \{ \ell_{0} + i t : 0 \leq i < L^{\prime} \}$ and $|ity - ip\tilde{h}_{\alpha_{0}}| < \frac{1}{3}\tilde{h}_{\alpha_{0}}$ for all $0 \leq i < L^{\prime}$ as otherwise the claim follows from Lemma \ref{thirdA}.

By Lemma \ref{thirdB} applied to each $\ell_{0} + ti$,
\begin{align*}
&\left|\sum_{\ell_{0}+ti \in \mathcal{L}} \delta_{\ell_{0}+ti}\muB(T^{q - (\ell_{0} + ti)k}B)\right|
< \left|\sum_{\ell_{0}+ti \in \mathcal{L}} \delta_{\ell_{0}+ti} \frac{r_{\alpha_{0}}}{r_{\alpha_{0}}+1} \gam{\alpha_{0}}{q}{k}{\ell_{0} + it} \frac{1}{b_{\alpha_{0}}} \sum_{j=0}^{b_{\alpha_{0}}-1} \muB(T^{q_{j} - h_{j}i}B)\right| \\
&\quad\quad + \sum_{\ell_{0}+ti\in \mathcal{L}} \delta_{\ell_{0}+ti} \left(
\left(1 - \gam{\alpha_{0}}{q}{k}{\ell_{0} + it} \right) \epsilon + 
\gam{\alpha_{0}}{q}{k}{\ell_{0} + it} \frac{1}{b_{\alpha_{0}}} \sum_{j=0}^{b_{\alpha_{0}}-1}\left(1 - \gam{\alpha_{j}}{q_{j}}{h_{j}}{i}\right)\epsilon + \tau_{\alpha_{0}} + \frac{1}{a_{\alpha_{0}-1}b_{\alpha_{0}-1}} \right) \\
&< \left|\sum_{\ell_{0}+ti \in \mathcal{L}} \delta_{\ell_{0}+ti} \gam{\alpha_{0}}{q}{k}{\ell_{0} + it} \frac{1}{b_{\alpha_{0}}} \sum_{j=0}^{b_{\alpha_{0}}-1} \muB(T^{q_{j} - h_{j}i}B)\right| + \sum_{\ell_{0}+ti \in \mathcal{L}} \delta_{\ell_{0}+ti}\left(1 - \gam{\alpha_{0}}{q}{k}{\ell_{0} + it} \right) \epsilon \\
&\quad\quad + \sum_{\ell_{0}+ti \in \mathcal{L}}\delta_{\ell_{0}+ti} \gam{\alpha_{0}}{q}{k}{\ell_{0} + it} \frac{1}{b_{\alpha_{0}}} \sum_{j=0}^{b_{\alpha_{0}}-1}\left(1 - \gam{\alpha_{j}}{q_{j}}{h_{j}}{i}\right)\epsilon + \sum_{\ell_{0}+ti \in \mathcal{L}}\ \delta_{\ell_{0}+ti} \left(\tau_{\alpha_{0}} + \frac{1}{a_{\alpha_{0}-1}b_{\alpha_{0}-1}}\right)
\end{align*}
Then, since $(1 - \delta) + \delta(1 - \gamma) = 1 - \delta\gamma$,
\begin{align*}
&\left|\sum_{\ell_{0}+ti \in \mathcal{L}} \delta_{\ell_{0}+ti}\muB(T^{q - (\ell_{0} + ti)k}B)\right| + \sum_{\ell_{0}+ti \in \mathcal{L}} (1 - \delta_{\ell_{0}+ti})\epsilon \\
&\quad\quad < \frac{1}{b_{\alpha_{0}}}\sum_{j=0}^{b_{\alpha_{0}}-1}\left|\sum_{\ell_{0}+ti \in \mathcal{L}} \delta_{\ell_{0}+ti}\gam{\alpha_{0}}{q}{k}{\ell_{0}+ti}\muB(T^{q_{j}-h_{j}i}B)\right| + \sum_{\ell_{0}+ti \in \mathcal{L}} \left(1 - \delta_{\ell_{0}+ti}\gam{\alpha_{0}}{q}{k}{\ell_{0}+ti}\right)\epsilon \\
&\quad\quad\quad\quad + \sum_{\ell_{0}+ti \in \mathcal{L}}\delta_{\ell_{0}+ti} \gam{\alpha_{0}}{q}{k}{\ell_{0} + it} \frac{1}{b_{\alpha_{0}}} \sum_{j=0}^{b_{\alpha_{0}}-1}\left(1 - \gam{\alpha_{j}}{q_{j}}{h_{j}}{i}\right)\epsilon + \sum_{\ell_{0}+ti \in \mathcal{L}}\ \delta_{\ell_{0}+ti} \left(\tau_{\alpha_{0}} + \frac{1}{a_{\alpha_{0}-1}b_{\alpha_{0}-1}}\right)
\end{align*}
Since $\muB(T^{-t}B) = \muB(T^{t}B)$, then for $\ell \notin \mathcal{L}$, $|\muB(T^{q - \ell k}B)| < \epsilon$ as it is bounded by one of $M_{B,\alpha_{0}}$, $M_{B,\alpha_{0}-1}$, $\widehat{M}_{B,\alpha_{0}}$ or $\widetilde{M}_{B,\alpha_{0}}$.
For $\ell_{0} + ti \notin \mathcal{L}$, since $\gam{\alpha_{0}}{q}{k}{\ell_{0}+ti} = 0$ and $|\muB(T^{q-(\ell_{0}+ti)}B)| < \epsilon$,
then 
\[
\left|\delta_{\ell_{0}+ti}\muB(T^{q - (\ell_{0} + ti)k}B)\right| + (1 - \delta_{\ell_{0}+ti})\epsilon < \epsilon = \left(1 - \delta_{\ell_{0}+ti}\gam{\alpha_{0}}{q}{k}{\ell_{0}+ti}\right)\epsilon
\]
Therefore, again using that $\gam{\alpha_{0}}{q}{k}{\ell_{0}+ti} = 0$ for $\ell_{0} + ti \notin \mathcal{L}$,
\begin{align*}
&\left|\sum_{i=0}^{L^{\prime}-1} \delta_{\ell_{0}+ti}\muB(T^{q - (\ell_{0} + ti)k}B)\right| + \sum_{i=0}^{L^{\prime}-1} (1 - \delta_{\ell_{0}+ti})\epsilon \\
&\quad\quad < \frac{1}{b_{\alpha_{0}}}\sum_{j=0}^{b_{\alpha_{0}}-1}\left|\sum_{i=0}^{L^{\prime}-1} \delta_{\ell_{0}+ti}\gam{\alpha_{0}}{q}{k}{\ell_{0}+ti}\muB(T^{q_{j}-h_{j}i}B)\right| + \sum_{i=0}^{L^{\prime}-1} \left(1 - \delta_{\ell_{0}+ti}\gam{\alpha_{0}}{q}{k}{\ell_{0}+ti}\right)\epsilon \\
&\quad\quad\quad\quad + \sum_{i=0}^{L^{\prime}-1}\delta_{\ell_{0}+ti} \gam{\alpha_{0}}{q}{k}{\ell_{0} + ti} \frac{1}{b_{\alpha_{0}}} \sum_{j=0}^{b_{\alpha_{0}}-1}\left(1 - \gam{\alpha_{j}}{q_{j}}{h_{j}}{i}\right)\epsilon +\sum_{i=0}^{L^{\prime}-1} \delta_{\ell_{0}+ti} \left(\tau_{\alpha_{0}} + \frac{1}{a_{\alpha_{0}-1}b_{\alpha_{0}-1}}\right)
\end{align*}
Rewriting the sums over $0 \leq i < L^{\prime}$ as over $0 \leq \ell < L$ with the indicator function $\bbone_{\ell < L^{\prime}}$,
\begin{align*}
&\left|\frac{1}{L}\sum_{i=0}^{L^{\prime}-1} \delta_{\ell_{0}+ti}\muB(T^{q - (\ell_{0} + ti)k}B)\right| + \frac{1}{L}\sum_{i=0}^{L^{\prime}-1} (1 - \delta_{\ell_{0}+ti})\epsilon \\
&\quad\quad < \frac{1}{b_{\alpha_{0}}}\sum_{j=0}^{b_{\alpha_{0}}-1}\left|\frac{1}{L}\sum_{\ell=0}^{L-1} \bbone_{\ell < L^{\prime}}\delta_{\ell_{0}+t\ell}\gam{\alpha_{0}}{q}{k}{\ell_{0}+t\ell}\muB(T^{q_{j}-h_{j}\ell}B)\right| + \frac{1}{L}\sum_{\ell=0}^{L -1} \bbone_{\ell < L^{\prime}} \left(1 - \delta_{\ell_{0}+t\ell}\gam{\alpha_{0}}{q}{k}{\ell_{0}+t\ell}\right)\epsilon \\
&\quad\quad\quad\quad + \frac{1}{L}\sum_{\ell=0}^{L-1}\bbone_{\ell < L^{\prime}}\delta_{\ell_{0}+t\ell} \gam{\alpha_{0}}{q}{k}{\ell_{0} + t\ell} \frac{1}{b_{\alpha_{0}}} \sum_{j=0}^{b_{\alpha_{0}}-1}\left(1 - \gam{\alpha_{j}}{q_{j}}{h_{j}}{\ell}\right)\epsilon +\tau_{\alpha_{0}} + \frac{1}{a_{\alpha_{0}-1}b_{\alpha_{0}-1}}
\end{align*}
Since $|\muB(T^{\ell}B)| < \epsilon$ for $\ell \ne \ell_{0} + ti$ as then $\ell \notin \mathcal{L}$,
\begin{align*}
&\left|\frac{1}{L}\sum_{\ell \ne \ell_{0}+ti} \delta_{\ell_{0}+ti}\muB(T^{q - (\ell_{0} + ti)k}B)\right| + \frac{1}{L}\sum_{\ell \ne \ell_{0} + ti} (1 - \delta_{\ell_{0}+ti})\epsilon
< \frac{1}{L} \sum_{\ell \ne \ell_{0} + ti} \epsilon \\
&\quad\quad\quad\quad = \frac{L - L^{\prime}}{L} \epsilon = \frac{1}{L}\sum_{\ell=L^{\prime}}^{L-1} \epsilon
= \frac{1}{L}\sum_{\ell=0}^{L-1} \left(1 - \bbone_{\ell < L^{\prime}}\right) \epsilon
= \frac{1}{L} \sum_{\ell=0}^{L-1} \left(1 - \bbone_{\ell < L^{\prime}}\delta_{\ell_{0}+t\ell}\gam{\alpha_{0}}{q}{k}{\ell_{0}+t\ell}\right) \epsilon
\end{align*}
Therefore the claim follows from the triangle inequality by splitting $0 \leq \ell < L$ into $\{ \ell_{0} + ti : 0 \leq i < L^{\prime} \}$ and $\{ 0 \leq \ell < L : \ell \ne \ell_{0} + ti \}$.
\end{proof}
  
 We are now ready to prove mixing for the remaining times.  The method of proof is to apply Lemma \ref{third} iteratively in the sense that that lemma either gives that the average over the $k$ powers is already mixed or that it can be reduced to a convex combination of `weighted' averages of the same form Lemma \ref{third} applies to.  This iterative process must either terminate, in which case mixing is shown, or reach a point where the powers appearing the average are smaller than $L$, in which case mixing follows from the `weighted' (weak) power ergodicity already established (it is crucial here that $L$ not decrease).
 
 \begin{proof}[Proof of Proposition \ref{fourth}]
 Fix $\epsilon > 0$ and set $\epsilon_{0} = (2\ceiling{\epsilon^{-1}}^{\ceiling{\epsilon^{-1}}+1})^{-1}$.  Using Propositions \ref{mixingES}, \ref{mixingES2}, \ref{midy} and \ref{deltape} and that $\sum_{n} \tau_{n} < \infty$ and $\sum_{n} \frac{1}{a_{n}b_{n}} < \infty$, there exists $N$ such that 
 \begin{itemize}
 \item $b_{N} > 4\epsilon^{-1}\epsilon_{0}^{-1}$;
 \item $\sup_{m \geq N-1} M_{B,m} < \epsilon$;
 \item $\sup_{m \geq N} \widehat{M}_{B,m} < \epsilon$;
 \item $\sup_{m \geq N} \widetilde{M}_{B,m} < \epsilon$;
 \item $\sum_{n=N}^{\infty} \tau_{n} < \epsilon$; 
 \item $\sum_{n=N-1}^{\infty} \frac{1}{a_{n}b_{n}} < \epsilon$; and
 \item  $\sup_{m \geq b_{N}-1} \sup_{k \leq m} \left(\int \left|\frac{1}{m}\sum_{j=0}^{m-1}\chi_{B}\circ T^{-jk}\right|~d\mu + \frac{2}{m}\right) < \frac{\epsilon}{3}$.
 \end{itemize}
 
 Take any $n$ such that $b_{n} > \tilde{h}_{N+1}$.  For $b_{n} \leq k < a_{n}$ and $|q| < a_{n-1}\tilde{h}_{n-1}$, by Lemma \ref{bigsmally},
 \[
 \left|\muB(T^{k\tilde{h}_{n} + q}B)\right| < \frac{a_{n}-k}{a_{n}}\left|\frac{1}{b_{n}}\sum_{\ell = 0}^{b_{n}-1} \muB(T^{q - k\ell}B)\right| + \frac{k}{a_{n}}\epsilon + \tau_{n}
 < \left|\frac{1}{b_{n}}\sum_{\ell = 0}^{b_{n}-1} \muB(T^{q - k\ell}B)\right| + 2\epsilon
 \]
 
 Set $L = b_{n}$.
By Lemma \ref{third}, $\left|\frac{1}{L}\sum_{\ell = 0}^{L-1} \muB(T^{q - k\ell}B)\right| < 4\epsilon^{1/2}$ or there exists $q_{\ell^{\prime}}, k_{\ell^{\prime}}^{\prime},L^{\prime},\ell_{0},t$ such that
\begin{align*}
\left|\frac{1}{L}\sum_{\ell = 0}^{L-1} \muB(T^{q - k\ell}B)\right|
&< \frac{1}{b_{\alpha_{0}}}\sum_{\ell^{\prime}=0}^{b_{\alpha_{0}}-1}\left(\left|\frac{1}{L}\sum_{\ell=0}^{L-1}\bbone_{\ell < L^{\prime}}\gam{\alpha_{0}}{q}{k}{\ell_{0}+t\ell}\muB(T^{q_{\ell^{\prime}}-k_{\ell^{\prime}}^{\prime}\ell}B)\right|
+ \frac{1}{L}\sum_{\ell=0}^{L-1}\left(1 - \bbone_{\ell < L^{\prime}}\gam{\alpha_{0}}{q}{k}{\ell_{0}+t\ell}\right)\epsilon\right) \\
&\quad\quad + \frac{1}{b_{\alpha_{0}}}\sum_{\ell^{\prime}=0}^{b_{\alpha_{0}}-1}\frac{1}{L}\sum_{\ell=0}^{L-1}\bbone_{\ell < L^{\prime}}\gam{\alpha_{0}}{q}{k}{\ell_{0}+t\ell}\left(1 - \gam{\alpha_{\ell^{\prime}}}{q_{\ell^{\prime}}}{k_{\ell^{\prime}}^{\prime}}{\ell}\right)\epsilon + \tau_{\alpha_{0}} + \frac{1}{a_{\alpha_{0}-1}b_{\alpha_{0}-1}}
\end{align*}

\textit{The first average above consists of terms which are exactly of the form that Lemma \ref{third} applies to.  By Lemma \ref{third}, some of those terms are bounded by $4\epsilon^{1/2}$ and the rest likewise split into a convex combination of terms of that same form so our aim is to iterate this process.  At each stage of this process, the `error bound' increases by $\tau_{\alpha}$ and $\frac{1}{a_{\alpha}b_{\alpha}}$ (for strictly decreasing $\alpha$) and by a term similar to the `extra $\epsilon$' term above involving two $\gamma$'s.}

Let $\mathcal{L}^{\prime} = \{ 0 \leq \ell^{\prime} < b_{\alpha_{0}} : \alpha_{\ell^{\prime}} > N~\text{and Lemma \ref{third} does not bound the $\ell^{\prime}$ weighted average by $4\epsilon^{1/2}$} \}$.
Then
\begin{align*}
\frac{1}{b_{\alpha_{0}}}\sum_{\ell^{\prime}=0}^{b_{\alpha_{0}}-1}&\left(\left|\frac{1}{L}\sum_{\ell=0}^{L-1}\bbone_{\ell < L^{\prime}}\gam{\alpha_{0}}{q}{k}{\ell_{0}+t\ell}\muB(T^{q_{\ell^{\prime}}-k_{\ell^{\prime}}^{\prime}\ell}B)\right|
+ \frac{1}{L}\sum_{\ell=0}^{L-1}\left(1 - \bbone_{\ell < L^{\prime}}\gam{\alpha_{0}}{q}{k}{\ell_{0}+t\ell}\right)\epsilon\right)
< \left(1 - \frac{|\mathcal{L}^{\prime}|}{b_{\alpha_{0}}}\right)4\epsilon^{1/2} \\
&+ \frac{|\mathcal{L}^{\prime}|}{b_{\alpha_{0}}} \frac{1}{|\mathcal{L}^{\prime}|}\sum_{\ell^{\prime}\in\mathcal{L}^{\prime}}
\left(\left|\frac{1}{L}\sum_{\ell=0}^{L-1}\bbone_{\ell < L^{\prime}}\gam{\alpha_{0}}{q}{k}{\ell_{0}+t\ell}\muB(T^{q_{\ell^{\prime}}-k_{\ell^{\prime}}^{\prime}\ell}B)\right|
+ \frac{1}{L}\sum_{\ell=0}^{L-1}\left(1 - \bbone_{\ell < L^{\prime}}\gam{\alpha_{0}}{q}{k}{\ell_{0}+t\ell}\right)\epsilon\right)
\end{align*}

Therefore, applying Lemma \ref{third} to each $\ell^{\prime}$ weighted average, since $\alpha_{\ell^{\prime}} \leq \alpha_{0} - 1$ (and suppressing the explicit dependence on $\ell^{\prime}$ of $L^{\prime\prime},\ell_{0}^{\prime},t^{\prime}$ for clarity),
\begin{align*}
&\left|\frac{1}{L}\sum_{\ell = 0}^{L-1} \muB(T^{q - k\ell}B)\right|
<  \left(1 - \frac{|\mathcal{L}^{\prime}|}{b_{\alpha_{0}}}\right)4\epsilon^{1/2} + \tau_{\alpha_{0}} + \tau_{\alpha_{0}-1} 
+ \frac{1}{a_{\alpha_{0}-1}b_{\alpha_{0}-1}} + \frac{1}{a_{\alpha_{0}-1}b_{\alpha_{0}-1}}
\\
&\quad\quad + \frac{|\mathcal{L}^{\prime}|}{b_{\alpha_{0}}}\frac{1}{|\mathcal{L}^{\prime}|}\sum_{\ell^{\prime}\in\mathcal{L}^{\prime}}
\frac{1}{b_{\alpha_{\ell^{\prime}}}}\sum_{\ell^{\prime\prime}=0}^{b_{\alpha_{\ell^{\prime}}}-1} \bigg{(} \left| \frac{1}{L}\sum_{\ell=0}^{L-1} \bbone_{\ell < L^{\prime\prime}} \bbone_{\ell_{0}^{\prime}+t^{\prime}\ell < L^{\prime}} \gam{\alpha_{0}}{q}{k}{\ell_{0} + t(\ell_{0}^{\prime}+t^{\prime}\ell)} \gam{\alpha_{\ell^{\prime}}}{q_{\ell^{\prime}}}{k_{\ell^{\prime}}^{\prime}}{\ell_{0}^{\prime}+t^{\prime}\ell}
\muB(T^{q_{\ell^{\prime},\ell^{\prime\prime}} - k_{\ell^{\prime},\ell^{\prime\prime}}\ell}B)\right| \\
&\quad\quad\quad\quad\quad\quad\quad\quad\quad\quad\quad\quad\quad + \frac{1}{L}\sum_{\ell=0}^{L-1} \left( 1 - \bbone_{\ell < L^{\prime\prime}} \bbone_{\ell_{0}^{\prime}+t^{\prime}\ell < L^{\prime}} \gam{\alpha_{0}}{q}{k}{\ell_{0} + t(\ell_{0}^{\prime}+t^{\prime}\ell)} \gam{\alpha_{\ell^{\prime}}}{q_{\ell^{\prime}}}{k_{\ell^{\prime}}^{\prime}}{\ell_{0}^{\prime}+t^{\prime}\ell}\right)\epsilon \bigg{)} \\
&\quad\quad + \frac{|\mathcal{L}^{\prime}|}{b_{\alpha_{0}}}\frac{1}{|\mathcal{L}^{\prime}|}\sum_{\ell^{\prime}\in\mathcal{L}^{\prime}}
\frac{1}{b_{\alpha_{\ell^{\prime}}}}\sum_{\ell^{\prime\prime}=0}^{b_{\alpha_{\ell^{\prime}}}-1} 
\frac{1}{L}\sum_{\ell=0}^{L-1} \bigg{(} \bbone_{\ell < L^{\prime}}\gam{\alpha_{0}}{q}{k}{\ell_{0}+t\ell}\left(1 - \gam{\alpha_{\ell^{\prime}}}{q_{\ell^{\prime}}}{k_{\ell^{\prime}}^{\prime}}{\ell}\right) \\
&\quad\quad\quad\quad\quad\quad\quad\quad\quad\quad\quad\quad\quad + \bbone_{\ell < L^{\prime\prime}}\bbone_{\ell_{0}^{\prime}+t^{\prime}\ell < L^{\prime}}\gam{\alpha_{0}}{q}{k}{\ell_{0} + t(\ell_{0}^{\prime}+t^{\prime}\ell)}\gam{\alpha_{\ell^{\prime}}}{q_{\ell^{\prime}}}{k_{\ell^{\prime}}^{\prime}}{\ell_{0}^{\prime}+t^{\prime}\ell}\left(1 - \gam{\alpha_{\ell^{\prime\prime}}}{q_{\ell^{\prime},\ell^{\prime\prime}}}{k_{\ell^{\prime},\ell^{\prime\prime}}}{\ell}\right) \bigg{)}\epsilon
\end{align*}

\textit{We now consider the `extra $\epsilon$' terms resulting from this process.  At this stage, there are two such terms.}

Now observe that
\begin{align*}
&\frac{1}{L}\sum_{\ell=0}^{L-1} \left( \bbone_{\ell < L^{\prime}}\gam{\alpha_{0}}{q}{k}{\ell_{0}+t\ell}\left(1 - \gam{\alpha_{\ell^{\prime}}}{q_{\ell^{\prime}}}{k_{\ell^{\prime}}^{\prime}}{\ell}\right) + \bbone_{\ell < L^{\prime\prime}}\bbone_{\ell_{0}^{\prime}+t^{\prime}\ell < L^{\prime}}\gam{\alpha_{0}}{q}{k}{\ell_{0} + t(\ell_{0}^{\prime}+t^{\prime}\ell)}\gam{\alpha_{\ell^{\prime}}}{q_{\ell^{\prime}}}{k_{\ell^{\prime}}^{\prime}}{\ell_{0}^{\prime}+t^{\prime}\ell}\left(1 - \gam{\alpha_{\ell^{\prime\prime}}}{q_{\ell^{\prime},\ell^{\prime\prime}}}{k_{\ell^{\prime},\ell^{\prime\prime}}}{\ell}\right) \right) \\
&= \frac{1}{L}\sum_{\substack{\ell \ne \ell_{0}^{\prime}+t^{\prime}i \\ \ell < L^{\prime}}} \gam{\alpha_{0}}{q}{k}{\ell_{0}+t\ell}\left(1 - \gam{\alpha_{\ell^{\prime}}}{q_{\ell^{\prime}}}{k_{\ell^{\prime}}^{\prime}}{\ell}\right)
+ \frac{1}{L}\sum_{\ell=0}^{L^{\prime\prime}-1} \bbone_{\ell_{0}^{\prime}+t^{\prime}\ell < L^{\prime}}\gam{\alpha_{0}}{q}{k}{\ell_{0} + t(\ell_{0}^{\prime}+t^{\prime}\ell)}
\left(1 - \gam{\alpha_{\ell^{\prime}}}{q_{\ell^{\prime}}}{k_{\ell^{\prime}}^{\prime}}{\ell_{0}^{\prime}+t^{\prime}\ell}
\gam{\alpha_{\ell^{\prime\prime}}}{q_{\ell^{\prime},\ell^{\prime\prime}}}{k_{\ell^{\prime},\ell^{\prime\prime}}}{\ell}\right)
\end{align*}
and that the sets of the original $0 \leq \ell < L$ the two sums range over are disjoint.  The left sum in the above is precisely those $\ell$ for which we do not need to continue iteratively applying Lemma \ref{third} (as they are exactly those which that lemma will bound by $4\epsilon^{1/2}$).  The right sum in the above are those $\ell$ for which the iterative process will need to continue.

\textit{The key point is that the terms on the right have the form $\gamma(1 - \gamma^{\prime} \gamma^{\prime\prime})$ rather than just $\gamma(1 - \gamma^{\prime})$.}

Continue iteratively applying Lemma \ref{third} until all terms are bounded by $4\epsilon^{1/2}$ or have $k_{\ell^{\prime},\ell^{\prime\prime},\ldots} < L$, which must occur as $\alpha$ decrements at each application of the lemma and if $k_{\ell,\ell^{\prime},\ldots} \geq L = b_{n} > \tilde{h}_{N+1}$ then $\alpha_{\ell^{\prime},\ell^{\prime\prime},\ldots} \geq N+1$ (and the hypotheses of the lemma hold as long as $\alpha_{\ell^{\prime},\ell^{\prime\prime},\ldots} > N$).

\textit{If the process terminates with one or more terms having $k$ values less than $L$ then the already established weighted weak power ergodicity (Proposition \ref{deltape}) tells us those terms are already bounded by $\epsilon$.}

If $k_{\ell^{\prime},\ell^{\prime\prime},\ldots} < L$ then, by the seventh requirement on the choice of $N$, as $L = b_{n} > \tilde{h}_{N+1} > b_{N}$,
\[
\left|\frac{1}{L}\sum_{\ell=0}^{L-1} \muB(T^{q_{\ell^{\prime},\ell^{\prime\prime},\ldots} - \ell k_{\ell^{\prime},\ell^{\prime\prime},\ldots}}B)\right|
\leq \int \left|\frac{1}{L}\sum_{\ell=0}^{L-1} \chi_{B} \circ T^{-\ell k_{\ell^{\prime},\ell^{\prime\prime},\ldots}} \right|~d\mu
<
 \frac{\epsilon}{3}
\]

The `extra $\epsilon$' terms weighted by $\gamma$ terms will be of the form
\[
\gamma(1 - \gamma^{\prime}) + \gamma \gamma^{\prime}(1 - \gamma^{\prime\prime}) + \cdots + \gamma \gamma^{\prime} \cdots \gamma^{(m-1)} (1 - \gamma^{(m)}) = \gamma(1 - \gamma\gamma^{\prime} \cdots \gamma^{(m-1)}\gamma^{(m)})
\]
where $m$ is the number of applications of Lemma \ref{third} needed for that $\ell$ (this is written carefully in the case of two applications above).

Then $\left|\frac{1}{L}\sum_{\ell=0}^{L-1}\muB(T^{q-k\ell}B)\right|$ is bounded by a convex combination of terms less than $4\epsilon^{1/2}$ plus a sum of $\tau$'s bounded by $\sum_{n=N}^{\infty} \tau_{n} < \epsilon$ plus a sum of $\frac{1}{ab}$ terms bounded by $\sum_{n=N}^{\infty} \frac{1}{a_{n}b_{n}} < \epsilon$ plus an average over $0 \leq \ell < L$ of terms of the form
\[
\gamma^{\alpha_{0}}\left(1 - \gamma^{\alpha_{\ell^{\prime}}}\gamma^{\alpha_{\ell^{\prime\prime}}}\cdots\gamma^{\alpha_{\ell^{(m)}}}\right)\epsilon
\]
which are all bounded by $\epsilon$ as $0 \leq \gamma \leq 1$.
Therefore
\[
\left|\frac{1}{L}\sum_{\ell = 0}^{L-1} \muB(T^{q - k\ell}B)\right| < 4\epsilon^{1/2} + \epsilon + \epsilon + \epsilon
\quad\quad
\text{meaning that}
 \quad\quad
  \left|\muB(T^{k\tilde{h}_{n} + q}B)\right| < 4\epsilon^{1/2} + 5\epsilon \qedhere
\]
\end{proof}
 
 \subsection{Proof of mixing}

\begin{proof}[Proof of Theorem \ref{C:mixing}]
By Propositions \ref{mixingES}, \ref{mixingES2}, \ref{midy} and \ref{fourth}, for any $B$ which is a union of levels in some $C_{N}$,
$
\lim_{n \to \infty} \max_{\tilde{h}_{n} \leq t < \tilde{h}_{n+1}} |\muB(T^{t}B)| = 0$.  As unions of levels generate the measure algebra, $T$ is Renyi mixing hence mixing.
\end{proof}

\section{Non-superlinear word complexity implies partial rigidity}\label{lowerbound}

\begin{theorem}\label{notmixing}
Let $X$ be a subshift with word complexity $p$ such that $\liminf \frac{p(q)}{q} < \infty$.  Then there exists a constant $\delta_{X} > 0$ such that every ergodic probability measure $\mu$ on $X$ is at least $\delta_{X}$-partially rigid.
\end{theorem}

\subsection{Word combinatorics}

\begin{notation}
For $x$ a finite or infinite word and $-\infty \leq i < j \leq \infty$,
\[
x_{[i,j)} = \text{ the subword of $x$ from position $i$ through position $j-1$ }
\]
\end{notation}

\begin{notation}
$
[w] = \{ x \in X : x_{[0,\len{w})} = w \}
$ for finite words $w$.
\end{notation}

\begin{notation}
For a word $v$ and $0 \leq q < \len{v}$, let $v^{\nicefrac{q}{\len{v}}}$ be the suffix of $v$ of length $q$.

Let $v^{n + \nicefrac{q}{\len{v}}} = v^{\nicefrac{q}{\len{v}}}v^{n}$ for $n \in \mathbb{N}$.
\end{notation}

\begin{definition}
Let $w \in \mathcal{L}(X)$ be a word in the language of a subshift.  A word $v \in \mathcal{L}(X)$ is a \textbf{root} of $w$ if $wv \in \mathcal{L}(X)$ and $\len{v} \leq \len{w}$ and $w$ is a suffix of $v^{\infty}$, i.e.~there exists $q = \nicefrac{p}{\len{v}}$ with $p \geq \len{v}$ such that $w = v^{q}$.
The \textbf{minimal root} of $w$ is the shortest $v$ which is a root of $w$.
\end{definition}

Every word has a unique minimal root as it is a root of itself.

\begin{lemma}\label{A1}
If $uw = wv$ and $\len{v} \leq \len{w}$ then $v$ is a root of $w$.
\end{lemma}
\begin{proof}
As $w$ has $v$ as a suffix, $w = w^{\prime}v$.  Then $uw^{\prime}v = uw = wv = w^{\prime}vv$ so $uw^{\prime} = w^{\prime}v$.  If $\len{w^{\prime}} \geq \len{v}$, repeat this process until it terminates at $w = w^{\prime\prime}v^{n}$ with $\len{w^{\prime\prime}} < \len{v}$.  Then $uw^{\prime\prime} = w^{\prime\prime}v$ so $w^{\prime\prime}$ is a suffix of $v$.
\end{proof}

\begin{lemma}\label{B1}
If $uv = vu$ then $u = v_{0}^{t}$ and $v = v_{0}^{s}$ for some word $v_{0}$ and $t,s \in \mathbb{N}$.
\end{lemma}
\begin{proof}
If $\len{u} = \len{v}$ then $uv = vu$ immediately implies $u=v$.  Let
\[
V = \{ (u,v) : uv = vu, \len{v} < \len{u},~\text{there is no word $v_{0}$ with $u = v_{0}^{t}$ and $v = v_{0}^{s}$ for $s,t \in \mathbb{N}$} \}
\]
and suppose $V \ne \emptyset$.
Let $(u,v) \in V$ such that $\len{u}$ is minimal.  As $\len{u} > \len{v}$, $uv = vu$ implies $u = vu^{\prime} = u^{\prime\prime}v$ for some nonempty words $u^{\prime}, u^{\prime\prime}$.  Then $vu^{\prime}v = uv = vu = vu^{\prime\prime}v$ so $u^{\prime} = u^{\prime\prime}$ and $vu^{\prime} = u^{\prime}v$.  If $\len{u^{\prime}} = \len{v}$ then $u^{\prime} = v$ so $u = v^{2}$ contradicting that $(u,v) \in V$.

Consider when $\len{u^{\prime}} < \len{v}$.  Since $\len{u^{\prime}} < \len{u}$ and $\len{v} < \len{u}$, the minimality of $\len{u}$ implies that $(v,u^{\prime}) \notin V$.  Then $v = v_{0}^{n}$ and $u^{\prime} = v_{0}^{m}$ for some word $v_{0}$ and $n,m \in \mathbb{N}$.  So $u = v_{0}^{n+m}$ meaning $(u,v) \notin V$.
When $\len{v} < \len{u^{\prime}}$, we have $(u^{\prime},v) \notin V$ so $u^{\prime} = v_{0}^{n}$ and $u = v_{0}^{n+m}$.  So $V = \emptyset$.
\end{proof}

\begin{lemma}\label{C1}
If $u$ and $v$ are both roots of a word $w$ and $uu$ is a suffix of $w$ and $\len{v} < \len{u}$ then there exists a suffix $v_{0}$ of $v$ such that $u = v_{0}^{n}$ and $v = v_{0}^{m}$ for some $n,m \in \mathbb{N}$.

In particular, if $v$ is the minimal root of $w$ and $u$ is a root of $w$ and $uu$ is a suffix of $w$ then $u$ is a multiple of $v$, i.e.~there exists $n \in \mathbb{N}$ such that $u = v^{n}$.
\end{lemma}
\begin{proof}
Writing $u^{\prime}$ and $v^{\prime}$ for the appropriate suffixes of $u$ and $v$, we have $w = u^{\prime}u^{t} = v^{\prime}v^{q}$ for some $t,q \in \mathbb{N}$.  Then $u = u_{0}v^{a}$ for some proper (possibly empty) suffix $u_{0}$ of $v$ and $1 \leq a \leq q$.  So $u^{\prime}(u_{0}v^{a})^{t} = v^{\prime}v^{q}$ meaning that $u^{\prime}(u_{0}v^{a})^{t-1}u_{0} = v^{\prime}v^{q-a}$.
As $t \geq 2$, $\len{v^{\prime}v^{q-a}} = \len{u^{\prime}(u_{0}v^{a})^{t-1}u_{0}} \geq \len{u_{0}v^{a}u_{0}} \geq \len{vu_{0}}$ so, as $u_{0}$ is a suffix of $v$, then $v^{\prime}v^{q-a}$ has $u_{0}v$ as a suffix.  This means $vu_{0} = u_{0}v$ so Lemma \ref{B1} gives $v_{0}$ such that $v = v_{0}^{n}$ and $u_{0} = v_{0}^{m}$ so $u = v_{0}^{m+an}$.  If $v$ is the minimal root then $v = v_{0}$ since $v_{0}$ is a root of $w$.
\end{proof}

\begin{lemma}\label{D1}
Let $w$ be a word with minimal root $v$.
If $0 \leq i \leq \frac{1}{2}\len{w}$ and $T^{i}[w] \cap [w] \ne \emptyset$ then $i$ is a multiple of $\len{v}$.
\end{lemma}
\begin{proof}
Let $u$ be the prefix of $B$ of length $i$ and $v_{0}$ be the suffix of $B$ of length $i$.  For $x \in T^{i}[w] \cap [w]$, then $x_{[-i,\len{w})} = uw = wv_{0}$.  By Lemma \ref{A1}, then $v_{0}$ is a root of $w$.  As $\len{v_{0}} = i \leq \frac{1}{2}\len{w}$, $w$ has $v_{0}v_{0}$ as a suffix.  By Lemma \ref{C1}, since $v$ is the minimal root then $v_{0}$ is a multiple of $v$.
\end{proof}

\subsection{Language analysis}

\begin{proposition}\label{words}
There exists $C,k > 0$, depending only on $X$, and $\ell_{n} \to \infty$ and, for each $n$, at most $C$ words $B_{n,j}$ so that $X_{0} = \{ x \in X : \text{every finite subword of $x$ is a subword of a concatenation of the $B_{n,j}$} \}$ has measure one.

Let $h_{n,j} = \len{B_{n,j}}$.  Then $\max_{j} h_{n,j} \leq k\ell_{n}$ and $\min_{j} h_{n,j} \to \infty$.  Let
\[
W_{B_{n,j}} = W_{n,j} = \{ x \in X_{0} : \text{$x$ can be written as a concatenation such that $x_{[0,h_{n,j})} = B_{n,j}$} \} \subseteq [B_{n,j}]
\]
There exists $c_{n,j} \leq k\ell_{n}$ such that the sets $T^{i}W_{n,j}$ are disjoint over $0 \leq i < c_{n,j}$.

For $j$ such that $h_{n,j} > \frac{1}{2}\ell_{n}$, $c_{n,j} \geq \frac{1}{2}\ell_{n}$.

For $j$ such that $h_{n,j} \leq \frac{1}{2}\ell_{n}$, $c_{n,j} = h_{n,j}$.  For such $j$, also $W_{n,j} = T^{\ell_{n}}[B_{n,j}^{\ell_{n}/h_{n,j}}B_{n,j}]$ and $B_{n,j}$ is the minimal root of $B_{n,j}^{\ell_{n}/h_{n,j}}B_{n,j}$.

If $x \in T^{h_{n,j}}W_{n,j} \cap W_{n,j^{\prime}}$ for $j \ne j^{\prime}$ and $h_{n,j^{\prime}} \leq \frac{1}{2}\ell_{n}$ then $x_{(-\infty,0)}$ has $B_{n,j^{\prime}}^{\ell_{n}/h_{n,j^{\prime}}}$ as a suffix and does not have $B_{n,j^{\prime}}^{\ell_{n}/h_{n,j^{\prime}}}B_{n,j^{\prime}}$ as a suffix.
\end{proposition}
\begin{proof}
Since $\liminf \frac{p(q)}{q} < \infty$, \cite{boshernitzan} Theorem 2.2 gives a constant $k$ and $\ell_{n} \to \infty$ such that $p(\ell_{n}+1) - p(\ell_{n}) \leq k$ and $p(\ell_{n}) \leq k \ell_{n}$.  
We perform an analysis similar to Ferenczi \cite{ferenczi1996rank} Proposition 4.

Let $G_{q}$ be the Rauzy graphs: the vertices are the words of length $q$ in $\mathcal{L}(X)$ and the directed edges are from words $w$ to $w^{\prime}$ such that $wa = bw^{\prime} \in \mathcal{L}(X)$ for some letters $a$ and $b$ and we label the edge with the letter $a$.  As $\mu$ is ergodic, exactly one strongly connected component has measure one and the rest have measure zero so we may assume $G_{q}$ is strongly connected.

Let $V_{q}^{RS}$ be the set of all vertices with more than one outgoing edge, i.e.~the right-special vertices.  Let $\mathcal{B}_{q}$ be the set of all paths from some $v \in V_{q}^{RS}$ to some $v^{\prime} \in V_{q}^{RS}$ that do not pass through any $v^{\prime\prime} \in V_{q}^{RS}$.  Then every $v \in V_{q}$ is necessarily along such a path.  Given any word $w$ in $\mathcal{L}(X)$, there exists $x \in X$ such that $x_{[0,\len{w})} = w$ so $w$ is the label of the path from the vertex corresponding to $x_{[-q,0)}$ to the vertex corresponding to $x_{[\len{w}-q,\len{w})}$ hence is a subword of some concatenation of labels of paths in $\mathcal{B}_{q}$.

The labels of the paths between right-special vertices are nested: $\mathcal{B}_{q+1}$ is necessarily a concatenation of paths in $\mathcal{B}_{q}$ since words corresponding to elements of $V_{q+1}^{RS}$ necessarily have right-special suffixes.  There are therefore recursion formulas defining $\mathcal{B}_{q+1}$ in terms of $\mathcal{B}_{q}$ though we do not make use of this fact.

Writing $\mathrm{outdeg}(v)$ for the number of outgoing edges of a vertex,
$
\sum_{v \in V_{\ell_{n}}^{RS}} (\mathrm{outdeg}(v) - 1) = p(\ell_{n}+1) - p(\ell_{n}) \leq k
$
meaning that $|V_{\ell_{n}}^{RS}| \leq k$ and therefore $\sum_{v \in V_{\ell_{n}}^{RS}} \mathrm{outdeg}(v) \leq 2k$.  Therefore $|\mathcal{B}_{\ell_{n}}| \leq 2k$.  No path in $\mathcal{B}_{\ell_{n}}$ properly contains a cycle so $\len{B} \leq p(\ell_{n}) \leq k \ell_{n}$ for any label $B$ of a path in $\mathcal{B}_{\ell_{n}}$.

Let $\mathcal{B}_{n}^{g}$ be the set of all concatenations of paths in $\mathcal{B}_{\ell_{n}}$ of total length at least $\frac{3}{2}\ell_{n}$ and at most $k\ell_{n}$ not properly containing any cycles.  As such a path contains no cycle properly, it has at most $|\mathcal{B}_{\ell_{n}}| \leq 2k$ segments from some vertex in $V_{\ell_{n}}^{RS}$ to another, so there are at most $K = \sum_{t=1}^{2k} (2k)^{t}$ such paths.

Let $\mathcal{B}_{n}^{c}$ be the set of all concatenations of paths in $\mathcal{B}_{\ell_{n}}$ of total length less than $\frac{3}{2}\ell_{n}$ which are simple cycles.  Then $|\mathcal{B}_{n}^{c}| \leq K$ as each path has at most $2k$ segments and at most $2k$ choices for each segment.  Every bi-infinite concatenation of paths in $\mathcal{B}_{\ell
_n}$ is necessarily a concatenation of paths in $\mathcal{B}_{n}^{g} \cup \mathcal{B}_{n}^{c}$.

Let $B$ be the label of a path in $\mathcal{B}_{n}^{g}$ and let $v$ be its minimal root.  Suppose that $\len{v} < \frac{1}{2}\ell_{n}$.  Then the vertex at which the path corresponding to $B$ ends is the word $v^{\ell_{n}/\len{v}}$ as it must be a suffix of $B$.
Let $B^{\prime}$ such that $B = B^{\prime}v$.  Then $\len{B^{\prime}} = \len{B} - \len{v} \geq \frac{3}{2}\ell_{n} - \len{v} > \ell_{n}$.  Then the path corresponding to $B$ reaches its final vertex twice as $B^{\prime}$ has suffix $v^{\ell_{n}/\len{v}}$ corresponding to that vertex.  This means the path properly contains a cycle which is a contradiction.  So all labels of paths in $\mathcal{B}_{n}^{g}$ have minimal root of length at least $\frac{1}{2}\ell_{n}$.
By Lemma \ref{D1}, then $T^{i}W_{n,j} \cap W_{n,j} \ne \emptyset$ for $0 < i \leq \frac{1}{2}\len{B}$ only when $i$ is a multiple of $\len{v}$.  Set $c_{n,j} = \min(\len{v},\frac{1}{2}\len{B}) \geq \frac{1}{2}\ell_{n}$ and then $T^{i}W_{n,j}$ are disjoint over $0 \leq i < c_{n,j}$.

Let $B$ be the label of a simple cycle beginning and ending at the word $w$.  Since $B$ is the label of a path beginning at $w$, every appearance of $B$ as a label in $x \in X$ is preceded by $w$, i.e.~$W_{B} \subseteq T^{\ell_{n}}[wB]$.  Since $B$ either has $w$ as a suffix or $B$ is a root of $w$ by Lemma \ref{A1}, $B$ is a root of $wB$.  Let $v$ be the minimal root of $wB$ and write $B = B^{\prime}v$.  Then $wB^{\prime}$ has $v$ as a root and $\len{wB^{\prime}} = \ell_{n} + \len{B^{\prime}}$ so $wB^{\prime}$ has suffix $v^{\ell_{n}/\len{v}}$.  If $B^{\prime}$ is nonempty then the path corresponding to $B$ passes through its final vertex before the path ends, contradicting that it is a simple cycle.  So $B = v$ is the minimal root of $wB$.

Then Lemma \ref{D1} implies that $T^{i}W_{B} \cap W_{B} \ne \emptyset$ for $0 < i \leq \frac{1}{2}\len{wB}$ only when $i$ is a multiple of $\len{B}$.  So if $\len{B} > \frac{1}{2}\ell_{n}$ then set $c_{n,j} = \min(\len{B},\frac{1}{2}\len{wB}) > \frac{1}{2}\ell_{n}$.  If $\len{B} \leq \frac{1}{2}\ell_{n}$, set $c_{n,j} = \len{B}$.  For such $B$, since $W_{B} \subseteq T^{\ell_{n}} [wB]$, we have that every occurrence of $B$ as a label of a path is preceded by $w = B^{\ell_{n}/\len{B}}$.  Moreover, if $x_{[-\ell_{n},\len{B})} = wB$ then $x_{[0,\len{B})}$ is the label of a path beginning at the vertex $w$ and ending at $w$ so $x \in W_{B}$.

For $x \in W_{B}$, if $x_{(-\infty,0)}$ has $B^{\ell_{n}/\len{B}}B$ as a suffix then the path reaches $w$ prior to the final $B$ in that suffix.  As no word $B^{\prime}$ appearing in the concatenation is the label of a path properly containing a cycle, this means the word preceding $x_{[0,\len{B})} = B$ in $x$ must be $B$, i.e.~$x \in T^{\ell_{n}+\len{B}}[B^{\ell_{n}/\len{B}}B]$ so $x \in T^{\len{B}}W_{B} \cap W_{B}$ and $x \notin T^{\len{B^{\prime}}}W_{B^{\prime}} \cap W_{B}$ for every $B^{\prime} \ne B$ as the path for $B^{\prime}$ does not properly contain a cycle.

Let $\mathcal{B}_{n}^{*} = \mathcal{B}_{n}^{g} \cup \mathcal{B}_{n}^{c}$.  Then $|\mathcal{B}_{n}^{*}| \leq 2K = C$ for all $n$ and every word in $\mathcal{L}(X)$ is a subword of some concatenation of labels of paths in $\mathcal{B}_{n}^{*}$.  Let $\mathcal{R}_{n}$ be the set of all labels of paths in $\mathcal{B}_{n}^{*}$.

Let $\mathcal{D}_{M} = \{ B : \len{B} \leq M \text{ and } B \in \mathcal{R}_{n} \text{ infinitely often} \}$.  Then $|D_{M}| < \infty$ as there only finitely many words of length at most $M$ (as non-superlinear complexity implies finite alphabet rank \cite{DDMP2}).
Let $X_{M}$ be the set of $x \in X$ such that for infinitely many $n$, $x$ cannot be written as a concatenation of labels in $\mathcal{B}_{n}^{*}$ without using at least one label in $\mathcal{D}_{M}$.

For $x \in X_{M}$, there exist infinitely many $t$ such that $x$ has $B_{t}^{r_{t}}$ as a subword for some $B_{t} \in \mathcal{D}_{M}$ and $r_{t} \to \infty$ (since the label $B_{t}$ is preceded by the word $B_{t}^{\floor{\ell_{n}/(\len{B_{t}})}}$).  As $|\mathcal{D}_{M}| < \infty$, there exists $B$ such that $B_{t} = B$ infinitely often.  Then $B^{r_{t}}$ is a subword of $x$ for $r_{t} \to \infty$ meaning $x$ is periodic.  Therefore $\bigcup_{M} X_{M} \subseteq \{ \text{periodic words} \}$ so $\mu(\bigcup_{M} X_{M}) = 0$ as $\mu$ is ergodic hence nonatomic and a periodic word of positive measure would be an atom (there are at most countably many periodic words).

Define $\{ B_{n,j} \}$ to be the set of all labels of paths in $\mathcal{B}_{n}^{*}$ which are in $\mathcal{R}_{n} \setminus \bigcup_{M} \mathcal{D}_{M}$.
If $\liminf_{n} \min_{j} \len{B_{n,j}} < \infty$ then $B_{n,j} = B$ for some fixed $B$ infinitely often (as there are only finitely many words of up to some fixed length).  But then $B \in \mathcal{D}_{\len{B}}$, a contradiction, so $\lim_{n} \min_{j} \len{B_{n,j}} = \infty$.
As $X_{0} = X \setminus \bigcup_{M} X_{M}$, we have $\mu(X_{0}) = 1$.
\end{proof}

\subsection{Measure-theoretic analysis}

\begin{definition}
Let $\displaystyle C_{n,j} = \bigcup_{i=0}^{h_{n,j}-1} T^{i}W_{n,j}$.
\end{definition}

\begin{definition}
For $j$ such that $\len{B_{n,j}} \leq \frac{1}{2}\ell_{n}$, let
\begin{align*}
Z_{n,j} &= [B_{n,j}^{\ell_{n}/h_{n,j}}B_{n,j}] \setminus T^{h_{n,j}}[B_{n,j}^{\ell_{n}/h_{n,j}}B_{n,j}] \\
&= \{ x \in X : x_{[0,\ell_{n}+h_{n,j})} = B_{n,j}^{\ell_{n}/h_{n,j}}B_{n,j}~\text{and}~x_{[-h_{n,j},\ell_{n})} \ne B_{n,j}^{\ell_{n}/h_{n,j}}B_{n,j} \}
\end{align*}
\end{definition}

\begin{proposition}
For $j$ such that $\len{B_{n,j}} \leq \frac{1}{2}\ell_{n}$, the sets $T^{ah_{n,j}}Z_{n,j}$ are disjoint over $0 \leq a \leq \Bigfloor{\frac{\ell_{n}}{h_{n,j}}}$.
\end{proposition}
\begin{proof}
For $0 \leq a < b \leq \Bigfloor{\frac{\ell_{n}}{h_{n,j}}}$ and $x \in T^{ah_{n,j}}Z_{n,j} \cap T^{bh_{n,j}}Z_{n,j}$, writing $z = \ell_{n} - \Bigfloor{\frac{\ell_{n}}{h_{n,j}}}h_{n,j}$, we would have  $x_{[z-(a+1)h_{n,j},z-ah_{n,j})} \neq B_{n,j}$ but $x_{[z-bh_{n,j},z)} = B_{n,j}^{b}$ which is impossible.
\end{proof}

\begin{proposition}
For $j$ such that $\len{B_{n,j}} \leq \frac{1}{2}\ell_{n}$, the sets $T^{i}Z_{n,j}$ are disjoint over $0 \leq i < c_{n,j}$.
\end{proposition}
\begin{proof}
Lemma \ref{D1} as
$B_{n,j}$ is the minimal root of $B_{n,j}^{\ell_{n}/h_{n,j}}B_{n,j}$ and $c_{n,j} \leq \frac{1}{2}\ell_{n} < \frac{1}{2}\len{B_{n,j}^{\ell_{n}/h_{n,j}}B_{n,j}}$.
\end{proof}

\begin{definition}
For $j$ such that $\len{B_{n,j}} > \frac{1}{2}\ell_{n}$, let, for $0 \leq i < c_{n,j}$,
\[
I_{n,j,i} = T^{i}W_{n,j}
\]
and for $j$ such that $\len{B_{n,j}} \leq \frac{1}{2}\ell_{n}$, let, for $0 \leq i < c_{n,j}$,
\[
I_{n,j,i} = T^{i} \Big{(}\bigsqcup_{a=0}^{\floor{\frac{\ell_{n}}{h_{n,j}}}} T^{ah_{n,j}}Z_{n,j}\Big{)}
\]
\end{definition}

As $T$ is measure-preserving, $\mu(I_{n,j,i}) = \mu(I_{n,j,0})$ for all $n$, $j$ and $0 \leq i < c_{n,j}$.

\begin{definition}
Let $\displaystyle \tilde{C}_{n,j} = \bigsqcup_{i=0}^{c_{n,j}-1} I_{n,j,i}$.
For $j$ such that $\len{B_{n,j}} \leq \frac{1}{2}\ell_{n}$, let
$\displaystyle
\widehat{C}_{n,j} = \bigsqcup_{i=0}^{h_{n,j}-1} T^{i}W_{n,j}
$.
\end{definition}

\begin{proposition}\label{526}
For $j$ such that $\len{B_{n,j}} > \frac{1}{2}\ell_{n}$, we have
$\mu(\tilde{C}_{n,j}) \geq \frac{1}{2k}\mu(C_{n,j})$.
\end{proposition}
\begin{proof}
$\displaystyle
\mu(C_{n,j}) \leq h_{n,j} \mu(W_{n,j}) = h_{n,j} \mu(I_{n,j,0}) = \frac{h_{n,j}}{c_{n,j}} \mu(\tilde{C}_{n,j})
\leq \frac{k\ell_{n}}{\frac{1}{2}\ell_{n}} \mu(\tilde{C}_{n,j})
= 2k \mu(\tilde{C}_{n,j}) 
$.
\end{proof}

\begin{proposition}\label{hinf}
$\lim_{n} \max_{j} \{ \mu(I_{n,j,0}) \} = 0$.
\end{proposition}
\begin{proof}
For $j$ such that $\len{B_{n,j}} > \frac{1}{2}\ell_{n}$, we have 
$1 \geq \mu(\tilde{C}_{n,j}) = c_{n,j} \mu(I_{n,j,0}) \geq \frac{1}{2}\ell_{n} \mu(I_{n,j,0})$ and $\ell_{n} \to \infty$.  For $j$ such that $\len{B_{n,j}} \leq \frac{1}{2}\ell_{n}$, we have $1 \geq \mu(\tilde{C}_{n,j}) = h_{n,j} \mu(I_{n,j,0})$ and $\min_{j} h_{n,j} \to \infty$.
\end{proof}

\begin{proposition}\label{Bs2}
$T^{h_{n,j}}W_{n,j} \subseteq \bigcup_{j^{\prime}} W_{n,j^{\prime}}$ and $X_{0} = \bigcup_{j} C_{n,j}$.
\end{proposition}
\begin{proof}
Every $x \in X_{0}$ is a concatenation of words of the form $B_{n,j}$ so every occurrence of $B_{n,j}$ is followed immediately by some $B_{n,j^{\prime}}$ and $x_{[0,\infty)} = uB_{1}B_{2}\cdots$ for some $u$ a suffix of some $B_{n,j}$ and $B_{\ell} \in \{ B_{n,j} \}$.
\end{proof}

\begin{proposition}\label{Bs}
Let $E \subseteq W_{n,j}$.  Then there exists $j^{\prime}$ such that $\mu(T^{h_{n,j}}E \cap W_{n,j^{\prime}}) \geq \frac{1}{C}\mu(E)$.
\end{proposition}
\begin{proof}
$T^{h_{n}}E = T^{h_{n}}E \cap T^{h_{n,j}}W_{n,j} \subseteq T^{h_{n}}E \cap \bigcup_{j^{\prime}} W_{n,j^{\prime}}$ and there are at most $C$ choices of $j^{\prime}$.
\end{proof}

\begin{lemma}\label{BC}
$\mu(W_{n,j}) \geq \frac{1}{k\ell_{n}} \mu(\tilde{C}_{n,j})$.
\end{lemma}
\begin{proof}
For $j$ such that $\len{B_{n,j}} \leq \frac{1}{2}\ell_{n}$, by Proposition \ref{words},
$
T^{-\ell_{n}}W_{n,j} = [B_{n,j}^{\ell_{n}/h_{n,j}}B_{n,j}] \supseteq Z_{n,j}
$
so
\begin{align*}
\mu(W_{n,j}) &\geq \mu(Z_{n,j}) = \frac{1}{\floor{\frac{\ell_{n}}{h_{n,j}}}+1} \mu(I_{n,j,0}) \geq \frac{1}{\frac{\ell_{n}}{h_{n,j}}}\frac{1}{h_{n,j}} \mu(\tilde{C}_{n,j})
= \frac{1}{\ell_{n}}\mu(\tilde{C}_{n,j})
\end{align*}
and for $j$ such that $\len{B_{n,j}} > \frac{1}{2}\ell_{n}$, we have
$
\mu(W_{n,j}) = \frac{1}{c_{n,j}}\mu(\tilde{C}_{n,j}) \geq \frac{1}{k\ell_{n}} \mu(\tilde{C}_{n,j})
$ since $c_{n,j} \leq k\ell_{n}$.
\end{proof}

\begin{proposition}\label{520}
If $\mu(T^{h_{n,j}} W_{n,j} \cap W_{n,j^{\prime}}) \geq \delta \mu(W_{n,j^{\prime\prime}})$ for $j \ne j^{\prime}$ then $\mu(\tilde{C}_{n,j^{\prime}}) \geq \frac{1}{2k}\delta\mu(\tilde{C}_{n,j^{\prime\prime}})$.
\end{proposition}
\begin{proof}
For $j^{\prime}$ such that $h_{n,j^{\prime}} < \frac{1}{2}\ell_{n}$, Proposition \ref{words} states that, as $j \ne j^{\prime}$, 
for $x \in T^{h_{n,j}}W_{n,j} \cap W_{n,j^{\prime}}$, the word $x_{(-\infty,0)}$
has $B_{n,j^{\prime}}^{\ell_{n}/h_{n,j^{\prime}}}$ as a suffix but does not have $B_{n,j^{\prime}}^{\ell_{n}/h_{n,j^{\prime}}}B_{n,j^{\prime}}$ as a suffix.  Therefore $T^{-\ell_{n}}(T^{h_{n,j}}W_{n,j} \cap W_{n,j^{\prime}}) \subseteq [B_{n,j^{\prime}}^{\ell_{n}/h_{n,j^{\prime}}}B_{n,j^{\prime}}] \setminus T^{h_{n,j^{\prime}}}[B_{n,j^{\prime}}^{\ell_{n}/h_{n,j^{\prime}}}B_{n,j^{\prime}}] = Z_{n,j^{\prime}}$.
This means that
$
\mu(Z_{n,j^{\prime}}) \geq \mu(T^{h_{n,j}} W_{n,j} \cap W_{n,j^{\prime}}) \geq \delta \mu(W_{n,j^{\prime\prime}})
$
so
\begin{align*}
\mu(\tilde{C}_{n,j^{\prime}}) &= h_{n,j^{\prime}} \mu(I_{n,j^{\prime},0})
= h_{n,j^{\prime}} \Big{(}\Bigfloor{\frac{\ell_{n}}{h_{n,j^{\prime}}}}+1\Big{)} \mu(Z_{n,j^{\prime}})
\geq h_{n,j^{\prime}}\frac{\ell_{n}}{h_{n,j^{\prime}}} \delta \mu(W_{n,j^{\prime\prime}}) \\
&\geq \ell_{n} \delta \frac{1}{c_{n,j^{\prime\prime}}}\mu(\tilde{C}_{n,j^{\prime\prime}})
\geq \ell_{n} \delta \frac{1}{k\ell_{n}}\mu(\tilde{C}_{n,j^{\prime\prime}}) 
= \delta \frac{1}{k} \mu(\tilde{C}_{n,j^{\prime\prime}})
\end{align*}
For $j^{\prime}$ such that $h_{n,j^{\prime}} > \frac{1}{2}\ell_{n}$, using Lemma \ref{BC} and that $\mu(W_{n,j^{\prime}}) \geq \delta\mu(W_{n,j^{\prime\prime}})$,
\[
\mu(\tilde{C}_{n,j^{\prime}}) = c_{n,j^{\prime}}\mu(W_{n,j^{\prime}})
\geq c_{n,j^{\prime}} \delta \mu(W_{n,j^{\prime\prime}})
\geq c_{n,j^{\prime}} \delta \frac{1}{k\ell_{n}} \mu(\tilde{C}_{n,j^{\prime\prime}})
\geq \frac{\ell_{n}}{2}  \delta \frac{1}{k\ell_{n}} \mu(\tilde{C}_{n,j^{\prime\prime}})
= \frac{1}{2k}\delta \mu(\tilde{C}_{n,j^{\prime\prime}}) \qedhere
\]
\end{proof}

\begin{proposition}\label{525}
For $j$ such that $\len{B_{n,j}} \leq \frac{1}{2}\ell_{n}$, we have
$\displaystyle
\mu(T^{h_{n,j}}I_{n,j,0} \cap I_{n,j,0}) \geq \frac{1}{2} \mu(I_{n,j,0})
$.
\end{proposition}
\begin{proof}
\begin{align*}
\mu(T^{h_{n,j}}I_{n,j,0} \cap I_{n,j,0}) 
&\geq \mu(\bigsqcup_{a=1}^{\floor{\frac{\ell_{n}}{h_{n,j}}}} T^{ah_{n,j}} Z_{n,j})
= \Bigfloor{\frac{\ell_{n}}{h_{n,j}}}\mu(Z_{n,j})
= \frac{\Bigfloor{\frac{\ell_{n}}{h_{n,j}}}}{\Bigfloor{\frac{\ell_{n}}{h_{n,j}}}+1} \mu(I_{n,j,0})
\geq \frac{1}{2}\mu(I_{n,j,0}) \qedhere
\end{align*}
\end{proof}

\subsection{Partial rigidity}

We employ ideas similar to Danilenko's \cite{danilenkopr} proof that exact finite rank implies partial rigidity:

\begin{proposition}\label{prig}
If there exists $\delta > 0$ and $j_{n}$ and $t_{n} \to \infty$ with $\mu(\tilde{C}_{n,j_{n}}) \geq \delta$ (or $\mu(\widehat{C}_{n,j_{n}}) \geq \delta$ when applicable) and $\mu(T^{t_{n}}I_{n,j_{n}} \cap I_{n,j_{n}}) \geq \delta \mu(I_{n,j_{n}})$ then $(X,\mu)$ is $\frac{1}{2}\delta^{2}$-partially rigid.
\end{proposition}
\begin{proof}
Let $A = W_{N,J}$ for some fixed $N$ and $J$.  Define
$
\alpha_{n} = \{ 0 \leq i < c_{n,j_{n}} - h_{N,J} : I_{n,j_{n},i} \subseteq A \}
$.

For $j_{n}$ such that $h_{n,j_{n}} > \frac{1}{2}\ell_{n}$, 
if $x \in I_{n,j_{n},i} \cap W_{N,J}$ then $x_{[-i,-i+h_{n,j_{n}})} = B_{n,j_{n}}$ and $x_{[0,h_{N,J})} = B_{N,J}$ meaning that $(B_{n,j_{n}})_{[i,i+h_{N,J})} = B_{N,J}$.  This implies that $T^{i}W_{n,j_{n}} \subseteq W_{N,J}$ provided $i < h_{n,j_{n}} - h_{N,J}$.

For $j_{n}$ such that $h_{n,j_{n}} \leq \frac{1}{2}\ell_{n}$, if $x \in I_{n,j_{n},i} \cap W_{N,J}$ then $x_{[-i,-i+\ell_{n}/h_{n,j_{n}})} = B_{n,j_{n}}^{\ell_{n}/h_{n,j_{n}}}$ and $x_{[0,h_{N,J})} = B_{N,J}$ so $(B_{n,j_{n}}^{\ell_{n}/h_{n,j_{n}}})_{[i,i+h_{N,J})} = B_{N,J}$ which implies $I_{n,j_{n},i} \subseteq W_{N,J}$ provided $i < h_{n,j_{n}} - h_{N,J}$.

Therefore
$
(|\alpha_{n}| + h_{N,J})\mu(I_{n,j_{n},0}) \geq 
\mu(A \cap \tilde{C}_{n,j_{n}}) \geq |\alpha_{n}|\mu(I_{n,j_{n},0}) 
$.  Likewise, if $\len{B_{n,j_{n}}} \leq \frac{1}{2}\ell_{n}$ then
$
(|\alpha_{n}| + h_{N,J})\mu(W_{n,j_{n}}) \geq 
\mu(A \cap \widehat{C}_{n,j_{n}}) \geq |\alpha_{n}|\mu(W_{n,j_{n}}) 
$ using $\alpha_{n} = \{ 0 \leq i < h_{n,j_{n}} - h_{N,J} : T^{i}W_{n,j_{n}} \subseteq A \}$.

For $m < c_{n,j_{n}}$,
$
\mu(T^{m}\tilde{C}_{n,j_{n}} \symdiff \tilde{C}_{n,j_{n}}) \leq 2m \mu(I_{n,j_{n},0})
$, (and likewise $\mu(T^{m}\widehat{C}_{n,j} \symdiff \widehat{C}_{n,j}) \leq 2m\mu(W_{n,j})$ when applicable)
 therefore
\begin{align*}
\int \big{|}\bbone_{\tilde{C}_{n,j_{n}}}\circ T^{-m} - \bbone_{\tilde{C}_{n,j_{n}}}\big{|}^{2}~d\mu
&= 2\mu(\tilde{C}_{n,j_{n}}) - 2\mu(T^{m}\tilde{C}_{n,j_{n}} \cap \tilde{C}_{n,j_{n}}) 
\leq 2m\mu(I_{n,j_{n},0})
\end{align*}
Therefore for $M < c_{n,j_{n}}$,
\begin{align*}
\Big{|}\frac{1}{M}\sum_{m=1}^{M} &\mu(T^{-m}A \cap \tilde{C}_{n,j_{n}}) - \mu(A \cap \tilde{C}_{n,j_{n}})\Big{|}
= \Big{|}\frac{1}{M}\sum_{m=1}^{M} \mu(A \cap T^{m}\tilde{C}_{n,j_{n}}) - \mu(A \cap \tilde{C}_{n,j_{n}})\Big{|} \\
&\leq \frac{1}{M}\sum_{m=1}^{M} \big{|} \mu(A \cap T^{m}\tilde{C}_{n,j_{n}}) - \mu(A \cap \tilde{C}_{n,j_{n}})\big{|}
\leq \frac{1}{M}\sum_{m=1}^{M} \int_{A} \big{|} \bbone_{\tilde{C}_{n,j_{n}}} \circ T^{-m} - \bbone_{\tilde{C}_{n,j_{n }}}\big{|}~d\mu \\
&\leq \frac{1}{M}\sum_{m=1}^{M} \Big{(}\int \big{|} \bbone_{\tilde{C}_{n,j_{n}}} \circ T^{-m} - \bbone_{\tilde{C}_{n,j_{n }}}\big{|}^{2}~d\mu\Big{)}^{1/2}
\leq \frac{1}{M}\sum_{m=1}^{M} \sqrt{2m\mu(I_{n,j_{n},0})}
\leq \sqrt{2M\mu(I_{n,j_{n},0})} 
\end{align*}
The mean ergodic theorem gives $M$ such that
$
\int \big{|}\frac{1}{M}\sum_{m=1}^{M} \bbone_{A}\circ T^{m} - \mu(A) \big{|}^{2}d\mu < (\frac{1}{4}\delta\mu(A))^{2}
$
so
\begin{align*}
\Big{|}\frac{1}{M}\sum_{m=1}^{M} &\mu(T^{-m}A \cap \tilde{C}_{n,j_{n}}) - \mu(A)\mu(\tilde{C}_{n,j_{n}})\Big{|}
= \Big{|} \int_{\tilde{C}_{n,j_{n}}} \frac{1}{M}\sum_{m=1}^{M} \bbone_{A}\circ T^{m} - \mu(A)~d\mu\Big{|} \\
&\leq \int_{\tilde{C}_{n,j_{n}}} \big{|}\frac{1}{M}\sum_{m=1}^{M} \bbone_{A}\circ T^{m} - \mu(A) \big{|}~d\mu 
\leq \Big{(}\int \big{|}\frac{1}{M}\sum_{m=1}^{M} \bbone_{A}\circ T^{m} - \mu(A) \big{|}^{2}~d\mu\Big{)}^{1/2} 
< \frac{1}{4}\delta\mu(A)
\end{align*}
For $n$ large enough that $c_{n,j_{n}} > M$ and $\sqrt{2M\mu(I_{n,j_{n},0})} < \frac{1}{4}\delta\mu(A)$ (Proposition \ref{hinf} states $\mu(I_{n,j_{n},0}) \to 0$) then
$
|\mu(A \cap \tilde{C}_{n,j_{n}}) - \mu(A)\mu(\tilde{C}_{n,j_{n}})| < \frac{1}{2}\delta\mu(A)
$.
Then 
\begin{align*}
\mu(T^{t_{n}}A \cap A) &\geq \mu(T^{t_{n}}(A \cap \tilde{C}_{n,j_{n}}) \cap (A \cap \tilde{C}_{n,j_{n}}))
\geq \sum_{i \in \alpha_{n}} \mu(T^{t_{n}}T^{i}I_{n,j_{n},0} \cap T^{i}I_{n,j_{n},0}) \\
&= |\alpha_{n}| \mu(T^{t_{n}}I_{n,j_{n},0} \cap I_{n,j_{n},0})
\geq |\alpha_{n}| \delta \mu(I_{n,j_{n},0})
\geq \delta (\mu(A \cap \tilde{C}_{n,j_{n}}) - h_{N,J}\mu(I_{n,j_{n},0})) \\
&> \delta\Big{(}\mu(A) \mu(\tilde{C}_{n,j_{n}}) - \frac{1}{2}\delta\mu(A)\Big{)} - \delta h_{N,J} \mu(I_{n,j_{n},0}) \\
&\geq \delta\Big{(}\mu(A)\delta - \frac{1}{2}\delta\mu(A)\Big{)} - \delta h_{N,J} \mu(I_{n,j_{n},0}) 
= \frac{1}{2}\delta^{2}\mu(A) - \delta h_{N,J} \mu(I_{n,j_{n},0})
\end{align*}
with the same applying to $\hat{C}_{n,j_{n}}$ when applicable.
Therefore for fixed $N$ and $J$ and $0 \leq i < h_{N,J}$,
\[
\liminf \mu(T^{t_{n}}T^{i}W_{N,J} \cap T^{i}W_{N,J})
= \liminf \mu(T^{t_{n}}W_{N,J} \cap W_{N,J}) \geq \frac{1}{2}\delta^{2}\mu(W_{N,J}) = \frac{1}{2}\delta^{2}\mu(T^{i}W_{N,J})
\]
and since the sets $T^{i}W_{N,J}$ generate the Borel algebra, $\mu$ is $\frac{1}{2}\delta^{2}$-partially rigid.
\end{proof}

\begin{proof}[Proof of Theorem \ref{notmixing}]
We aim to apply Proposition \ref{prig}.  Set $\delta = \frac{1}{4k^{2}C^{C+1}}$ which depends only on $X$.

There exists $a_{0}$ such that $\mu(C_{n,a_{0}}) \geq \frac{1}{C}$ since $X_{0} = \bigcup_{j} C_{n,j}$.
If $\len{B_{n,a_{0}}} \leq \frac{1}{2}\ell_{n}$ then $\mu(\widehat{C}_{n,a_{0}}) = \mu(C_{n,a_{0}}) \geq \frac{1}{C}$ and Proposition \ref{525} implies
$
\mu(T^{h_{n,a_{0}}}I_{n,a_{0},0} \cap I_{n,a_{0},0}) \geq \frac{1}{2}\mu(I_{n,a_{0},0})
$
so take $t_{n} = h_{n,a_{0}}$ and $j_{n} = a_{0}$.

Now consider when $\len{B_{n,a_{0}}} > \frac{1}{2}\ell_{n}$ so Proposition \ref{526} implies $\mu(\tilde{C}_{n,a_{0}}) \geq \frac{1}{2k}\mu(C_{n,a_{0}}) \geq \frac{1}{2kC}$.

By Proposition \ref{Bs}, there exists $a_{1}$ such that $\mu(T^{h_{n,a_{0}}}W_{n,a_{0}} \cap W_{n,a_{1}}) \geq \frac{1}{C} \mu(W_{n,a_{0}})$.  If $a_{1} = a_{0}$ then $\mu(\tilde{C}_{n,a_{1}}) = \mu(\tilde{C}_{n,a_{0}}) \geq \frac{1}{2kC}$ and if $a_{1} \ne a_{0}$ then
Proposition \ref{520} implies $\mu(\tilde{C}_{n,a_{1}}) \geq \frac{1}{2k}\mu(\tilde{C}_{n,a_{0}}) \geq \frac{1}{4k^{2}C}$.

Proposition \ref{Bs} then says there exists $a_{2}$ such that
\[
\mu(T^{h_{n,a_{1}}}(T^{h_{n,a_{0}}}W_{n,a_{0}} \cap W_{n,a_{1}}) \cap W_{n,a_{2}}) \geq \frac{1}{C} \mu(T^{h_{n,a_{0}}}W_{n,a_{0}} \cap W_{n,a_{1}}) \geq \frac{1}{C^{2}}\mu(W_{n,a_{0}})
\]
and then Proposition \ref{520} gives $\mu(\tilde{C}_{n,a_{2}}) \geq \frac{1}{C^{2}}\frac{1}{2k}\mu(\tilde{C}_{n,a_{0}}) \geq \frac{1}{4k^{2}C^{3}}$.

Repeating this process, we obtain $a_{\ell}$ for $0 \leq \ell \leq C$ such that $\mu(\tilde{C}_{n,a_{\ell}}) \geq \frac{1}{4k^{2}C^{\ell+1}}\geq \frac{1}{4k^{2}C^{C+1}}$ and
\[
\mu(W_{n,a_{C}} \cap \bigcap_{\ell=0}^{C-1} T^{\sum_{z=\ell}^{C-1} h_{n,a_{z}}}W_{n,a_{\ell}}) \geq \frac{1}{C^{C}}\mu(W_{n,a_{0}})
\]

If any of the $a_{\ell}$ are such that $h_{n,a_{\ell}} \leq \frac{1}{2}\ell_{n}$ then Proposition \ref{525} implies
$
\mu(T^{h_{n,a_{\ell}}}I_{n,a_{\ell},0} \cap I_{n,a_{\ell},0}) \geq \frac{1}{2}\mu(I_{n,a_{\ell},0})
$
so take $t_{n} = h_{n,a_{\ell}}$ and $j_{n} = a_{\ell}$.

If $h_{n,a_{\ell}} > \frac{1}{2}\ell_{n}$ for all $0 \leq \ell \leq C$ then, since there are at most $C$ choices of $j$, for some $q < s$ we must have $a_{q} = a_{s}$ so setting $j_{n} = a_{q}$ and $t_{n} = \sum_{z=q}^{s-1} h_{n,a_{z}}$,
\[
\mu(T^{t_{n}}I_{n,j_{n},0} \cap I_{n,j_{n},0}) = \mu(T^{\sum_{z=q}^{s-1} h_{n,a_{z}}}W_{n,a_{q}} \cap W_{n,a_{s}})
\geq \mu(W_{n,a_{C}} \cap \bigcap_{\ell=0}^{C-1} T^{\sum_{z=\ell}^{C-1} h_{n,a_{z}}}W_{n,a_{\ell}})
 \geq \frac{1}{C^{C}}\mu(W_{n,a_{0}})
\]
\begin{flalign*}
\text{As} &&
\mu(W_{n,a_{0}}) &=
\mu(I_{n,a_{0},0}) = \frac{1}{c_{n,a_{0}}} \mu(\tilde{C}_{n,a_{0}}) \geq \frac{1}{h_{n,a_{0}}} \frac{1}{2kC}
\geq \frac{1}{k\ell_{n}} \frac{1}{2kC}
\geq \frac{1}{k\ell_{n}} \frac{1}{2kC} \mu(\tilde{C}_{n,j_{n}}) && \\
&&&= \frac{1}{k\ell_{n}} \frac{1}{2kC} c_{n,j_{n}} \mu(I_{n,j_{n},0})
\geq \frac{1}{k\ell_{n}} \frac{1}{2kC} \frac{\ell_{n}}{2} \mu(I_{n,j_{n},0})
= \frac{1}{4k^{2}C} \mu(I_{n,j_{n},0}) &&
\end{flalign*}
we then have $\mu(T^{t_{n}}I_{n,j_{n},0} \cap I_{n,j_{n},0}) \geq \frac{1}{4k^{2}C^{C+1}} \mu(I_{n,j_{n},0})$.

In all cases, by Proposition \ref{prig}, we have that $(X,\mu,T)$ is $\frac{1}{2}\delta^{2}$-partially rigid.
\end{proof}

\textbf{Acknowledgments}  The author thanks Ronnie Pavlov for introducing him to this question and for numerous enjoyable discussions.  The author thanks the referees for many helpful suggestions for improving the clarity of exposition, in particular for the suggestion to move the parts of the proof of mixing which follow from standard arguments to the appendix in order to focus on the novel techniques and for suggesting adding significant exposition in and around the use of said techniques.

\appendix

\section{Proofs of mixing properties}

We conclude with detailed proofs of mixing for sequences for which the mixing proof is essentially standard arguments for staircase adapted to our notation.

\subsection{Mixing between \texorpdfstring{$a_{n}\tilde{h}_{n}$}{a\_n h\_n} and \texorpdfstring{$\tilde{h}_{n+1}$}{h\_(n+1)}}

\begin{lemma}\label{hn}
Let $T$ be a quasi-staircase transformation.  Then for any $n$ and $0 \leq \ell < b_{n}$ and $k, i\geq 0$ such that $i + k \leq a_{n}$ and any $j \geq k\ell$,
\[
T^{k\tilde{h}_{n}}I_{n,j}^{[\ell a_{n} + i]} = I_{n,j-k\ell}^{[\ell a_{n} + i + k]}
\]
\end{lemma}
\begin{proof}
There are $c_{n} + \floor*{\frac{\ell a_{n} + i}{a_{n}}} = c_{n} + \ell$ spacers above $I_{n,j}^{[\ell a_{n} + i]}$ so $T^{\tilde{h}_{n}}I_{n,j}^{[\ell a_{n} + i]} = I_{n,j-\ell}^{[\ell a_{n} + i + 1]}$.  Since $i+k \leq a_{n}$, there are also $c_{n} + \ell$ spacers above each $I_{n,j-v\ell}^{[\ell a_{n} + i + v]}$ for $1 \leq v < k$ so applying $T^{h_{n}+c_{n}}$ repeated $k$ times, the claim follows.
\end{proof}

\begin{lemma}\label{mix1}
Let $T$ be a quasi-staircase transformation, $k \in \mathbb{N}$, $B$ a union of levels in some $C_{N}$ and $n \geq N$.  If $k < a_{n}$ and $kb_{n} < h_{n}$ then
\[
\sum_{j=0}^{h_{n}-1} |\muB(T^{k\tilde{h}_{n}}I_{n,j})| \leq  \int \left|\frac{1}{b_{n}}\sum_{\ell=0}^{b_{n}-1} \chi_{B} \circ T^{-k\ell}\right|~d\mu + \frac{k+1}{a_{n}} + \frac{kb_{n}}{h_{n}}
\]
\end{lemma}
\begin{proof}
By Lemma \ref{hn} and then Lemma \ref{L:mixtrick}, for $kb_{n} \leq j < h_{n}$,
\begin{align*}
|\muB(&T^{k\tilde{h}_{n}}I_{n,j})|
= \left|\sum_{\ell=0}^{b_{n}-1}\sum_{i=0}^{a_{n}-1} \muB(T^{k\tilde{h}_{n}}I_{n,j}^{[\ell a_{n} + i]}) + \muB(T^{k\tilde{h}_{n}}I_{n,j}^{[r_{n}]})\right| \\
&\leq  \left|\sum_{\ell=0}^{b_{n}-1}\sum_{i=0}^{a_{n}-k-1} \muB(T^{k\tilde{h}_{n}}I_{n,j}^{[\ell a_{n} + i]})\right| + (b_{n}k+1)\mu(I_{n+1}) \\
&=  \left|\sum_{\ell=0}^{b_{n}-1}\sum_{i=0}^{a_{n}-k-1} \muB(I_{n,j-k\ell}^{[\ell a_{n} + i + k]})\right| + (b_{n}k+1)\mu(I_{n+1}) \\
&=  \left| \sum_{\ell=0}^{b_{n}-1}\sum_{i=0}^{a_{n}-k-1} \frac{1}{r_{n}+1}\muB(I_{n,j-k\ell})\right| + \frac{b_{n}k + 1}{r_{n}+1}\mu(I_{n}) \\
&=  \left| \frac{1}{r_{n}+1} \sum_{\ell=0}^{b_{n}-1}\sum_{i=0}^{a_{n}-k-1} \muB(T^{-k\ell}I_{n,j})\right| + \frac{b_{n}k + 1}{r_{n}+1}\mu(I_{n}) \\
&= \frac{a_{n}-k}{r_{n}+1} \left|\sum_{\ell=0}^{b_{n}-1} \muB(T^{-k\ell}I_{n,j})\right| + \frac{b_{n}k + 1}{r_{n}+1}\mu(I_{n}) 
\leq \frac{1}{b_{n}} \left| \sum_{\ell=0}^{b_{n}-1} \muB(T^{-k\ell}I_{n,j})\right| + \frac{k+1}{a_{n}}\mu(I_{n}) \\
&= \left|\frac{1}{b_{n}}\sum_{\ell=0}^{b_{n}-1} \int_{I_{n,j}} \chi_{B} \circ T^{-k\ell}~d\mu\right| + \frac{k+1}{a_{n}}\mu(I_{n}) 
\leq \int_{I_{n,j}} \left|\frac{1}{b_{n}}\sum_{\ell=0}^{b_{n}-1}\chi_{B} \circ T^{-k\ell}\right|~d\mu + \frac{k+1}{a_{n}}\mu(I_{n})
\end{align*}
Therefore
\begin{align*}
\sum_{j=0}^{h_{n}-1} |\muB(T^{k\tilde{h}_{n}}I_{n,j})|
&\leq \sum_{j=kb_{n}}^{h_{n}-1} |\muB(T^{k\tilde{h}_{n}}I_{n,j})| + kb_{n}\mu(I_{n}) \\
&\leq \sum_{j=kb_{n}}^{h_{n}-1}\left(\int_{I_{n,j}} \left|\frac{1}{b_{n}}\sum_{\ell=0}^{b_{n}-1}\chi_{B} \circ T^{-k\ell}\right|~d\mu + \frac{k+1}{a_{n}}\mu(I_{n,j})\right) + kb_{n}\mu(I_{n}) \\
&\leq \int \left|\frac{1}{b_{n}}\sum_{\ell=0}^{b_{n}-1} \chi_{B} \circ T^{-k\ell}\right|~d\mu + \frac{k+1}{a_{n}} + \frac{kb_{n}}{h_{n}} \qedhere
\end{align*}
\end{proof}

\begin{proposition}\label{P:mix1}
Let $T$ be a quasi-staircase transformation and $k \in \mathbb{N}$.  If $T^{k}$ is ergodic then $\{k\tilde{h}_{n}\}$ and $\{kh_{n} \}$ are rank-one uniform mixing.
\end{proposition}
\begin{proof}
Since $\frac{b_{n}}{h_{n}} \to 0$ and $a_{n} \to \infty$ there exists $N$ such that for all $n \geq N$ we have $k < a_{n}$ and $kb_{n} < h_{n}$.   That $\{k\tilde{h}_{n}\}$ is rank-one uniform mixing follows from Lemma \ref{mix1} since $T^{k}$ is ergodic, $b_{n} \to \infty$, $a_{n} \to \infty$ and $\frac{b_{n}}{h_{n}} \to 0$.  Then
\[
\sum_{j=0}^{h_{n}-1}|\muB(T^{kh_{n}}I_{n,j})| \leq \sum_{j=kc_{n}}^{h_{n}}|\muB(T^{kh_{n}}I_{n,j})| + \frac{kc_{n}}{h_{n}} 
= \sum_{j=0}^{h_{n}-kc_{n}}|\muB(T^{k\tilde{h}_{n}}I_{n,j})| + \frac{kc_{n}}{h_{n}} \to 0
\]
as $\frac{c_{n}}{h_{n}} \to 0$, $k$ is fixed and $\{ k\tilde{h}_{n} \}$ is rank-one uniform mixing.
\end{proof}

\begin{lemma}[\cite{CPR} Proposition A.13]\label{P:CCC}
Let $T$ be a rank-one transformation and $\{ c_{n} \}$ a sequence such that $\frac{c_{n}}{h_{n}} \to 0$.  If $k \in \mathbb{N}$ and $\{ q(h_{n}+c_{n}) \}$ is rank-one uniform mixing for each $q \leq k+1$ and $\{ t_{n} \}$ is a sequence such that $h_{n} + c_{n} \leq t_{n} < (q+1)(h_{n}+c_{n})$ for all $n$ then $\{ t_{n} \}$ is mixing.
\end{lemma}

\begin{lemma}[\cite{CPR} Proposition A.16]\label{ppppp}
Let $T$ be a rank-one transformation and $\{ c_{n} \}$ a sequence such that $\frac{c_{n}}{h_{n}} \to 0$.  If $\{ q(h_{n}+c_{n}) \}$ is rank-one uniform mixing for each fixed $q$ and $k_n \to \infty$ is such that $\frac{k_n}{n} \leq 1$ then for any measurable set $B$,
	$
	\int |\frac{1}{n} \sum \limits_{j=0}^{n-1} \chi_{B} \circ T^{-jk_n}|~d\mu \rightarrow 0.
	$
\end{lemma}

\begin{proof}[Proof of Proposition \ref{P:pe}]
As $T$ is ergodic, Proposition \ref{P:mix1} with $k=1$ gives that $\{ \tilde{h}_{n} \}$ is rank-one uniform mixing, hence mixing, so $T$ is totally ergodic.  Then Proposition \ref{P:mix1} gives that for each fixed $k$ the sequence $\{k\tilde{h}_{n}\}$ is rank-one uniform mixing so Lemma \ref{ppppp} gives the claim.
\end{proof}

\begin{proposition}\label{mixingseqs}
Let $T$ be a quasi-staircase transformation, $B$ a measurable set and $Q > 0$.  Then 
\[
\max_{h_{n}+c_{n}\leq t < Q\tilde{h}_{n}} |\muB(T^{t}B)| \to 0
\]
\end{proposition}
\begin{proof}
As in the proof of Proposition \ref{P:pe}, for each fixed $k$ the sequence $\{k\tilde{h}_{n}\}$ is rank-one uniform mixing so Lemma \ref{P:CCC} gives the claim.
\end{proof}

\begin{lemma}\label{mixtrick2}
Let $T$ be a quasi-staircase transformation.  Let $n > 0$ and $0 \leq x < b_{n}$ and $0 \leq q < a_{n}$.

If $0 \leq \ell < b_{n} - x$ and $0 \leq i < a_{n} - q$ and $j \geq \frac{1}{2}a_{n}x(x-1) + qx + ix + \ell(xa_{n}+q)$ then
\[
T^{(xa_{n}+q)\tilde{h}_{n}}I_{n,j}^{[\ell a_{n} + i]} = I_{n,j - \frac{1}{2}a_{n}x(x-1) - qx - ix - \ell(xa_{n}+q)}^{[(\ell+x)a_{n} + i + q]}
\]
\end{lemma}
\begin{proof}
If $x = 0$ then Lemma \ref{hn} applied with $q$ in place of $k$ gives the claim.  So we can write
\[
xa_{n} + q = (a_{n} - i) + (x-1)a_{n} + (q + i)
\]
and assume all three terms on the right are nonnegative.

Using Lemma \ref{hn},
\[
T^{(a_{n}-i)\tilde{h}_{n}}I_{n,j}^{[\ell a_{n}+i]} = I_{n,j-(a_{n}-i)\ell}^{[\ell a_{n} + i + a_{n} - i]} = I_{n,j-(a_{n}-i)\ell}^{[(\ell+1)a_{n}]}
\]
Now observe that, by Lemma \ref{hn} with $0$ as $i$ and $a_{n}$ as $k$, for any $0 \leq v < x$ and any $a_{n}v \leq z < h_{n}$,
\[
T^{a_{n}\tilde{h}_{n}}I_{n,z}^{[v a_{n}]} = I_{n,z-a_{n}v}^{[(v+1)a_{n}]}
\]
so applying that $x-1$ times for $v = \ell+1, \ell+2, \ldots, \ell+x-1$,
\[
T^{(x-1)a_{n}\tilde{h}_{n}}I_{n,j-(a_{n}-i)\ell}^{[(\ell+1)a_{n}]}
= I_{n,j-(a_{n}-i)\ell - (x-1)\ell a_{n} - \frac{1}{2}x(x-1)a_{n}}^{[(\ell+x)a_{n}]}
\]
since $\sum_{v=\ell+1}^{\ell+x-1} v = \frac{1}{2}(\ell+x)(\ell + x-1) - \frac{1}{2}\ell(\ell+1)= (x-1)\ell + \frac{1}{2}x(x-1)$.  
Then applying Lemma \ref{hn} one final time with $q+i$ in place of $k$,
\begin{align*}
T^{(q+i)\tilde{h}_{n}}I_{nj-(a_{n}-i)\ell - (x-1)\ell a_{n} - \frac{1}{2}x(x-1)a_{n}}^{[(\ell+x)a_{n}]}
&= I_{n,j-(a_{n}-i)\ell - (x-1)\ell a_{n} - \frac{1}{2}x(x-1)a_{n} - (x+\ell)(q+i)}^{[(\ell+x)a_{n}+q+i]} \\
&= I_{n,j - x \ell a_{n} - \frac{1}{2}x(x-1)a_{n} - xi - xq - \ell q}^{[(\ell+x)a_{n}+q+i]} \qedhere
\end{align*}
\end{proof}

\begin{lemma}\label{mixtrick3}
Let $T$ be a quasi-staircase transformation.  Let $n > 0$ and $0 \leq x < b_{n}$ and $0 \leq q < a_{n}$.

If $0 \leq \ell < b_{n} - x - 1$ and $a_{n} - q \leq i < a_{n}$ and $j \geq \frac{1}{2}a_{n}x(x+1) + q(x+1) + i(x+1) + \ell(xa_{n}+1)$ then
\[
T^{(xa_{n}+q)\tilde{h}_{n}}I_{n,j}^{[\ell a_{n} + i]} = I_{n,j - \frac{1}{2}a_{n}x(x+1) - (q+i-a_{n})(x+1) - \ell(xa_{n}+q)}^{[(\ell+x)a_{n} + i + q]}
\]
\end{lemma}
\begin{proof}
The same proof as Lemma \ref{mixtrick2} except we write
$
xa_{n} + q = (a_{n} - i) + x a_{n} + (q + i - x)
$.
\end{proof}

\begin{lemma}\label{mixingESlemma}
Let $T$ be a quasi-staircase transformation.  
Let $B$ be a union of levels $C_N$.  For $n \geq N$ and $k_{n}\tilde{h}_{n} \leq t_{n} < (k_{n}+1)\tilde{h}_{n}$,
\[
\sum \limits_{j=0}^{h_n-1} |\muB(T^{t_n}I_{n,j})| \leq  \sum \limits_{x=0}^{h_n-1} |\muB(T^{k_n\tilde{h}_n}I_{n,x})|  \nonumber + c_{n} \mu(I_{n}) +\sum \limits_{x=0}^{h_n-1} |\muB(T^{(k_n+1)\tilde{h}_n}I_{n,x})|
\]
\end{lemma}
\begin{proof}
Write $t_{n} = k_{n}\tilde{h}_{n} + z_{n}$ for $0 \leq z_{n} < \tilde{h}_{n}$.  Then
	\begin{align}
	\sum \limits_{j=0}^{h_n-1} |\muB(T^{t_n}I_{n,j})| \nonumber 
	&\leq \sum_{j=0}^{h_{n}-z_{n}-1} |\muB(T^{t_n}I_{n,j})| + c_{n}\mu(I_{n}) + \sum_{j=h_{n}-z_{n}+c_{n}}^{h_{n}-1} |\muB(T^{t_n}I_{n,j})|	\\
	&\leq \sum \limits_{j=0}^{h_n-z_n-1} |\muB(T^{k_n\tilde{h}_n}I_{n,j+z_n}) |   + c_{n}\mu(I_{n}) \nonumber 
+\sum \limits_{j=\tilde{h}_n-z_n}^{h_n-1} |\muB(T^{(k_n+1)\tilde{h}_n}I_{n,j+z_n-\tilde{h}_n})| \nonumber \\
	&\leq \sum \limits_{x=0}^{h_n-1} |\muB(T^{k_n\tilde{h}_n}I_{n,x})|  \nonumber + c_{n} \mu(I_{n}) +\sum \limits_{x=0}^{h_n-1} |\muB(T^{(k_n+1)\tilde{h}_n}I_{n,x})|  \nonumber  \qedhere
	\end{align}
\end{proof}

\begin{proof}[Proof of Proposition \ref{mixingES}]
Let $t_{n}$ attain the maximum in $M_{B,n}$.
If $t_{n} \geq (r_{n}-1)\tilde{h}_{n}$ then $h_{n+1} + c_{n+1} - t_{n} \leq c_{n+1} + 2h_{n} + c_{n} + \frac{1}{2}a_{n}b_{n}(b_{n}-1)$ so
\begin{align*}
\sum_{j=0}^{h_{n+1}-1} |\muB(T^{t_{n}}I_{n+1,j})|
&\leq \sum_{j=h_{n+1}+c_{n+1}-t_{n}}^{h_{n+1}-1} |\muB(T^{t_{n}}I_{n+1,j})| + (h_{n+1}+c_{n+1}-t_{n})\mu(I_{n+1}) \\
&\leq \sum_{j=0}^{t_{n}-c_{n+1}-1} |\muB(T^{\tilde{h}_{n+1}}I_{n+1,j})| + \frac{c_{n+1} + 2h_{n} + c_{n} \frac{1}{2}a_{n}b_{n}(b_{n}-1)}{h_{n+1}} \to 0
\end{align*}
since $\{ \tilde{h}_{n+1} \}$ is rank-one uniform mixing.

So we may assume $t_{n} < (r_{n}-1)\tilde{h}_{n}$ and therefore write $t_{n} = k_{n}\tilde{h}_{n} + z_{n}$ for $a_{n} \leq k_{n} < r_{n}-1$ and $0 \leq z_{n} < \tilde{h}_{n}$.
By Lemma \ref{mixingESlemma},
\[
	\sum \limits_{j=0}^{h_n-1} |\muB(T^{t_n}I_{n,j})|  
	\leq \sum \limits_{x=0}^{h_n-1} |\muB(T^{k_n\tilde{h}_n}I_{n,x})|  + c_{n} \mu(I_{n}) +\sum \limits_{x=0}^{h_n-1} |\muB(T^{(k_n+1)\tilde{h}_n}I_{n,x})|   
\]
	
We will show the sum on the left tends to zero; the same argument with $k_{n}+1$ in place of $k_{n}$ gives the same for the right sum.  As $c_{n}\mu(I_{n}) \to 0$, this will complete the proof.

Write $k_{n} = x_{n}a_{n} + q_{n}$ for $0 \leq q_{n} < a_{n}$ and $1 \leq x_{n} < b_{n}$.  
Observe that
\begin{align*}
\sum_{j=0}^{h_{n}-1} |\muB(T^{k_{n}\tilde{h}_{n}}I_{n,j})|
&\leq \sum_{j=0}^{h_{n}-1}  \left|\sum_{\ell=0}^{b_{n}-x_{n}-2}\sum_{i=0}^{a_{n}-1} \muB(T^{k_{n}\tilde{h}_{n}}I_{n,j}^{[\ell a_{n} + i]})\right|
+ 2a_{n}h_{n} \mu(I_{n+1})  \tag{$\star$} \\
&\quad\quad + \sum_{j=0}^{h_{n}-1} \left|\sum_{\ell=b_{n}-x_{n}+1}^{b_{n}-1}\sum_{i=0}^{a_{n}-1} \muB(T^{k_{n}\tilde{h}_{n}}I_{n,j}^{[\ell a_{n} + i]})\right| + \frac{1}{r_{n}+1}  \tag{$\star\star$}
\end{align*}

We handle the sum ($\star\star$) first and return to the sum in $(\star)$ shortly.
	
For $0 \leq \ell < b_{n}$ and $0 \leq i < a_{n}$, we have that
\[
I_{n,0}^{[\ell a_{n} + i]} = T^{(\ell a_{n} + i)\tilde{h}_{n}}I_{n,\frac{1}{2}\ell(\ell-1)a_{n}+i\ell}^{[0]}
\]
since $\frac{1}{2}\ell(\ell-1)a_{n} + i\ell \leq a_{n}b_{n}^{2} + a_{n}b_{n} < h_{n}$ (as $\frac{a_{n}b_{n}^{2}}{h_{n}} \to 0$).

For $b_{n} - x_{n} + 1 \leq \ell < b_{n}$ and $0 \leq i < a_{n}$, since $x + \ell \geq b_{n} + 1$,
\begin{align*}
k_{n}\tilde{h}_{n} + (\ell a_{n} + i)\tilde{h}_{n}
&= (x_{n}a_{n} + q_{n} + \ell a_{n} + i)(h_{n} + c_{n}) \\
&\geq (b_{n}a_{n} + a_{n})\tilde{h}_{n} \\
&= (b_{n}a_{n}+1)h_{n} + b_{n}a_{n}c_{n} + (a_{n}-1)h_{n} + a_{n}c_{n}
\geq h_{n+1}
\end{align*}
since $\frac{1}{2}a_{n}b_{n}(b_{n}-1) \leq h_{n}$.  Also,
\begin{align*}
k_{n}\tilde{h}_{n} + (\ell a_{n} + i)\tilde{h}_{n} + \frac{1}{2}\ell(\ell-1)a_{n} + i\ell 
&= ((x_{n}+\ell)a_{n} + q_{n} + i)(h_{n} + c_{n}) + \frac{1}{2}\ell(\ell-1)a_{n} + i\ell \\
&\leq 2b_{n}a_{n}(h_{n} + c_{n}) + \frac{1}{2}b_{n}(b_{n}-1)a_{n} + a_{n}b_{n} < 2h_{n+1}
\end{align*}

Since a sublevel in $I_{n}$ is a level in $I_{n+1}$ and $\{ h_{n+1} \}$ is rank-one uniform mixing (Proposition \ref{P:mix1}),
\begin{align*}
\sum_{j=0}^{h_{n}-1} \sum_{\ell=b_{n}-x_{n}+1}^{b_{n}-1} \sum_{i=0}^{a_{n}-1} | \muB(T^{k_{n}\tilde{h}_{n}}I_{n,j}^{[\ell a_{n}+i]}) |
&\leq \sum_{y=0}^{h_{n+1}-1} | \muB(T^{h_{n+1}}I_{n+1,y}) | \to 0
\end{align*}

As $2a_{n}h_{n}\mu(I_{n+1}) \leq \frac{2a_{n}h_{n}}{h_{n+1}} \leq \frac{2}{b_{n}} \to 0$ and $r_{n} \to \infty$, it remains only to show that the sum in  ($\star$) tends to zero.
Observe that
\begin{align*}
\sum_{\ell=0}^{b_{n}-x_{n}-2}\sum_{i=0}^{a_{n}-1} \muB(T^{k_{n}\tilde{h}_{n}}I_{n,j}^{[\ell a_{n} + i]}) 
 &= \sum_{\ell=0}^{b_{n}-x_{n}-2}\sum_{i=0}^{a_{n}-q_{n}-1} \muB(T^{k_{n}\tilde{h}_{n}}I_{n,j}^{[\ell a_{n} + i]}) \tag{$\dagger$} \\
&\quad\quad + \sum_{\ell=0}^{b_{n}-x_{n}-2}\sum_{i=a_{n}-q_{n}}^{a_{n}-1} \muB(T^{k_{n}\tilde{h}_{n}}I_{n,j}^{[\ell a_{n} + i]}) \tag{$\ddagger$}
\end{align*}

First, we address $(\dagger)$:
set $y_{n} = \frac{1}{2}a_{n}x_{n}(x_{n}-1) + q_{n}x_{n}$.  For $i < a_{n} - q_{n}$ and $\ell < b_{n} - x_{n} - 1$, we have $y_{n} + ix_{n} + \ell k_{n} \leq 3a_{n}b_{n}^{2}$ so
for $j \geq 3a_{n}b_{n}^{2}$, by Lemma \ref{mixtrick2} and Lemma \ref{L:mixtrick},
\begin{align*}
\sum_{\ell=0}^{b_{n}-x_{n}-2} &\sum_{i=0}^{a_{n}-q_{n}-1} \muB(T^{k_{n}\tilde{h}_{n}}I_{n,j}^{[\ell a_{n} + i]}) 
= \sum_{\ell=0}^{b_{n}-x_{n}-2}\sum_{i=0}^{a_{n}-q_{n}-1} \muB(I_{n,j - y_{n} - ix_{n} - \ell k_{n}}^{[(\ell + x_{n}) a_{n} + i + q_{n}]}) 
 \\
&= \frac{1}{r_{n}+1}\sum_{\ell=0}^{b_{n}-x_{n}-2}\sum_{i=0}^{a_{n}-q_{n}-1} \muB(I_{n,j - y_{n} - ix_{n} - \ell k_{n}}) 
=  \frac{1}{r_{n}+1}\sum_{\ell=0}^{b_{n}-x_{n}-2}\sum_{i=0}^{a_{n}-q_{n}-1} \muB(T^{-\ell k_{n} - ix_{n} - y_{n}} I_{n,j})
\end{align*}
Then, summing over all $3a_{n}b_{n}^{2} \leq j < h_{n}$,
\begin{align*}
\sum_{j=3a_{n}b_{n}^{2}}^{h_{n}-1} &\left| \sum_{\ell=0}^{b_{n}-x_{n}-2}\sum_{i=0}^{a_{n}-q_{n}-1} \muB(T^{k_{n}\tilde{h}_{n}}I_{n,j}^{[\ell a_{n} + i]}) \right| \\
&= \sum_{j=3a_{n}b_{n}^{2}}^{h_{n}-1} \left| \frac{1}{r_{n}+1}\sum_{\ell=0}^{b_{n}-x_{n}-2}\sum_{i=0}^{a_{n}-q_{n}-1} \muB(T^{-\ell k_{n} - ix_{n} - y_{n}} I_{n,j})\right| \\
&\leq \frac{1}{r_{n}+1} \sum_{j=0}^{h_{n}-1} \sum_{\ell=0}^{b_{n}-x_{n}-2} \left| \sum_{i=0}^{a_{n}-q_{n}-1} \muB(T^{-\ell k_{n} - ix_{n} - y_{n}} I_{n,j}) \right| \\
&\leq \frac{1}{r_{n}+1} \sum_{\ell=0}^{b_{n}-x_{n}-2} \int \left|\sum_{i=0}^{a_{n}-q_{n}-1} \chi_{B} \circ T^{-\ell k_{n} - i x_{n} - y_{n}}\right|~d\mu \\
&= \frac{(b_{n}-x_{n}-2)(a_{n}-q_{n})}{r_{n}+1} \int \left| \frac{1}{a_{n}-q_{n}}\sum_{i=0}^{a_{n}-q_{n}-1} \chi_{B} \circ T^{-ix_{n}}\right|~d\mu \\
&\leq \min\left(\frac{a_{n}-q_{n}}{a_{n}}, \int \left| \frac{1}{a_{n}-q_{n}}\sum_{i=0}^{a_{n}-q_{n}-1} \chi_{B} \circ T^{-ix_{n}}\right|~d\mu\right)
\end{align*}
since $\frac{(b_{n}-2)}{r_{n}+1} < \frac{1}{a_{n}}$ and $\int |\chi_{B}|~d\mu \leq 1$.  For a subsequence along which $x_{n} \leq a_{n}-q_{n}$, Proposition \ref{P:pe} implies the integral tends to zero.  For $n$ such that $a_{n} - q_{n} < x_{n} < b_{n}$, the quantity on the left is bounded by $\frac{b_{n}}{a_{n}} \to 0$.

For $(\ddagger)$: set $y_{n}^{\prime} = \frac{1}{2}a_{n}x_{n}(x_{n}+1) + (q_{n} - a_{n})(x_{n}+1)$.
By Lemma \ref{mixtrick3} and Lemma \ref{L:mixtrick}, for $j \geq 3a_{n}b_{n}^{2}$,
\begin{align*}
&\sum_{\ell=0}^{b_{n}-x_{n}-2}\sum_{i=a_{n}-q_{n}}^{a_{n}-1} \muB(T^{k_{n}\tilde{h}_{n}}I_{n,j}^{[\ell a_{n} + i]})
= \sum_{\ell=0}^{b_{n}-x_{n}-2}\sum_{i=a_{n}-q_{n}}^{a_{n}-1} \muB(I_{n,j - y_{n}^{\prime} - i(x_{n}+1) - \ell k_{n}}^{[(\ell+x_{n}) a_{n} + i + q_{n}]})  \\
&= \frac{1}{r_{n}+1} \sum_{\ell=0}^{b_{n}-x_{n}-2}\sum_{i=a_{n}-q_{n}}^{a_{n}-1} \muB(I_{n,j - y_{n}^{\prime} - i(x_{n}+1) - \ell k_{n}}) 
= \frac{1}{r_{n}+1} \sum_{\ell=0}^{b_{n}-x_{n}-2}\sum_{i=a_{n}-q_{n}}^{a_{n}-1} \muB(T^{-\ell k_{n}- i(x_{n}+1) - y_{n}^{\prime}}I_{n,j})
\end{align*}

Similar to the sum ($\dagger$), then
\begin{align*}
&\sum_{j=3a_{n}b_{n}^{2}}^{h_{n}-1}  \left|\sum_{\ell=0}^{b_{n}-x_{n}-2}\sum_{i=a_{n}-q_{n}}^{a_{n}-1} \muB(T^{k_{n}\tilde{h}_{n}}I_{n,j}^{[\ell a_{n} + i]}) \right| \\
&= \sum_{j=3a_{n}b_{n}^{2}}^{h_{n}-1} \left|  \frac{1}{r_{n}+1} \sum_{\ell=0}^{b_{n}-x_{n}-2}\sum_{i=a_{n}-q_{n}}^{a_{n}-1} \muB(T^{-\ell k_{n}- i(x_{n}+1) - y_{n}^{\prime}}I_{n,j})  \right| \\
&\leq \frac{(b_{n}-x_{n}-2)q_{n}}{r_{n}+1} \int \left| \frac{1}{q_{n}} \sum_{i=a_{n}-q_{n}}^{a_{n}-1} \chi_{B} \circ T^{-i(x_{n}+1)} \right|~d\mu \\
&= \frac{(b_{n}-x_{n}-2)q_{n}}{r_{n}+1} \int \left| \frac{1}{q_{n}} \sum_{i^{\prime}=0}^{q_{n}-1} \chi_{B} \circ T^{-i^{\prime}(x_{n}+1)} \right|~d\mu
\leq \min\Big{(} \frac{q_{n}}{a_{n}},  \int \left| \frac{1}{q_{n}} \sum_{i^{\prime}=0}^{q_{n}-1} \chi_{B} \circ T^{-i^{\prime}(x_{n}+1)} \right|~d\mu \Big{)}
\end{align*}
and along any subsequence where $x_{n}+1 \leq q_{n}$, this tends to zero by Proposition \ref{P:pe}, and for $q_{n} \leq x_{n}+ 1 < b_{n} + 1$, the quantity on the left is bounded by $\frac{b_{n}}{a_{n}} \to 0$, completing the proof.
 \end{proof}
 
  \subsection{Mixing between \texorpdfstring{$\tilde{h}_{n}$}{h\_n} and \texorpdfstring{$b_{n}\tilde{h}_{n}$}{b\_n h\_n}}
 
 \begin{proof}[Proof of Proposition \ref{mixingES2}]
Let $t_{n}$ attain the maximum in $\widehat{M}_{B,n}$.
 By Lemma \ref{mixingESlemma}, writing $t_{n} = k_{n}\tilde{h}_{n} + z_{n}$ for $1 \leq k_{n} < b_{n}$ and $0 \leq z_{n} < \tilde{h}_{n}$,
 \begin{align}
	\sum \limits_{j=0}^{h_n-1} &|\muB(T^{t_n}I_{n,j})| \leq \sum \limits_{x=0}^{h_n-1} |\muB(T^{k_n\tilde{h}_n}I_{n,x})|  \nonumber + c_{n} \mu(I_{n}) +\sum \limits_{x=0}^{h_n-1} |\muB(T^{(k_n+1)\tilde{h}_n}I_{n,x})|  \nonumber 
	\end{align}
By Lemma \ref{mix1},
\begin{align*}
\sum_{j=0}^{h_{n}-1} |\muB(T^{k_{n}\tilde{h}_{n}}I_{n,j})| 
&\leq \int \left|\frac{1}{b_{n}}\sum_{\ell=0}^{b_{n}-1} \chi_{B} \circ T^{-k_{n}\ell}\right|~d\mu + \frac{k_{n}+1}{a_{n}} + \frac{k_{n}b_{n}}{h_{n}} \to 0
\end{align*}
since $k_{n} < b_{n}$ so Proposition \ref{P:pe} implies the integral tends to zero.  Similar reasoning for $k_{n} + 1 \leq b_{n}$ then completes the proof.
 \end{proof}

\vspace{-1em}
\dbibliography{ComplexityLinear}


\newcommand{\etalchar}[1]{$^{#1}$}
\providecommand{\bysame}{\leavevmode\hbox to3em{\hrulefill}\thinspace}
\providecommand{\MR}{\relax\ifhmode\unskip\space\fi MR }
\providecommand{\MRhref}[2]{%
  \href{http://www.ams.org/mathscinet-getitem?mr=#1}{#2}
}
\providecommand{\href}[2]{#2}
\begin{thebibliography}{DDMP21}

\bibitem[Ada98]{adams1998smorodinsky}
Terrence~M. Adams, \emph{Smorodinsky's conjecture on rank-one mixing}, Proc.
  Amer. Math. Soc. \textbf{126} (1998), no.~3, 739--744. \MR{1443143}

\bibitem[AFP17]{adamsferenczipeterson17}
Terrence Adams, S\'{e}bastien Ferenczi, and Karl Petersen, \emph{Constructive
  symbolic presentations of rank one measure-preserving systems}, Colloq. Math.
  \textbf{150} (2017), no.~2, 243--255. \MR{3719459}

\bibitem[Bos85]{boshernitzan}
Michael Boshernitzan, \emph{A unique ergodicity of minimal symbolic flows with
  linear block growth}, J. Analyse Math. \textbf{44} (1984/85), 77--96.
  \MR{801288}

\bibitem[Cas97]{cassaigne}
Julien Cassaigne, \emph{Complexit\'{e} et facteurs sp\'{e}ciaux}, Journ\'{e}es
  Montoises (Mons, 1994), no.~1, vol.~4, Bull. Belg. Math. Soc. Simon Stevin,
  1997, pp.~67--88. \MR{1440670}

\bibitem[CFPZ19]{cassaigneetal}
Julien Cassaigne, Anna~E. Frid, Svetlana Puzynina, and Luca~Q. Zamboni, \emph{A
  characterization of words of linear complexity}, Proc. Amer. Math. Soc.
  \textbf{147} (2019), no.~7, 3103--3115. \MR{3973910}

\bibitem[CK19]{CK2}
Van Cyr and Bryna Kra, \emph{Counting generic measures for a subshift of linear
  growth}, J. Eur. Math. Soc. (JEMS) \textbf{21} (2019), no.~2, 355--380.
  \MR{3896204}

\bibitem[CK20a]{CK3}
\bysame, \emph{The automorphism group of a shift of slow growth is amenable},
  Ergodic Theory Dynam. Systems \textbf{40} (2020), no.~7, 1788--1804.
  \MR{4108905}

\bibitem[CK20b]{CK4}
\bysame, \emph{Realizing ergodic properties in zero entropy subshifts}, Israel
  J. Math. \textbf{240} (2020), no.~1, 119--148. \MR{4193129}

\bibitem[CPR23]{CPR}
Darren Creutz, Ronnie Pavlov, and Shaun Rodock, \emph{Measure-theoretically
  mixing subshifts with low complexity}, Ergodic Theory and Dynamical Systems
  \textbf{43} (2023), no.~7, 2293–2316.

\bibitem[Cre21]{Creutz2021}
Darren Creutz, \emph{Mixing on stochastic staircase tansformations}, Studia
  Math. \textbf{257} (2021), no.~2, 121--153. \MR{4194950}

\bibitem[CS04]{CreutzSilva2004}
Darren Creutz and Cesar~E. Silva, \emph{Mixing on a class of rank-one
  transformations}, Ergodic Theory Dynam. Systems \textbf{24} (2004), no.~2,
  407--440. \MR{2054050}

\bibitem[CS10]{CreutzSilva2010}
\bysame, \emph{Mixing on rank-one transformations}, Studia Math. \textbf{199}
  (2010), no.~1, 43--72. \MR{2652597}

\bibitem[Dan16]{danilenkopr}
Alexandre~I. Danilenko, \emph{Actions of finite rank: weak rational ergodicity
  and partial rigidity}, Ergodic Theory Dynam. Systems \textbf{36} (2016),
  no.~7, 2138--2171. \MR{3568975}

\bibitem[DDMP16]{DDMP}
Sebasti\'{a}n Donoso, Fabien Durand, Alejandro Maass, and Samuel Petite,
  \emph{On automorphism groups of low complexity subshifts}, Ergodic Theory
  Dynam. Systems \textbf{36} (2016), no.~1, 64--95. \MR{3436754}

\bibitem[DDMP21]{DDMP2}
\bysame, \emph{Interplay between finite topological rank minimal {C}antor
  systems, {$\mathcal{S}$}-adic subshifts and their complexity}, Trans. Amer.
  Math. Soc. \textbf{374} (2021), no.~5, 3453--3489. \MR{4237953}

\bibitem[dJ77]{deljunco}
Andr\'{e}s del Junco, \emph{A transformation with simple spectrum which is not
  rank one}, Canadian J. Math. \textbf{29} (1977), no.~3, 655--663. \MR{466489}

\bibitem[DOP21]{DOP}
Andrew Dykstra, Nicholas Ormes, and Ronnie Pavlov, \emph{Subsystems of
  transitive subshifts with linear complexity}, Ergodic Theory and Dynamical
  Systems (2021), 1–27.

\bibitem[Fer95]{ferenczichacon}
S\'{e}bastien Ferenczi, \emph{Les transformations de {Chacon} : combinatoire,
  structure g\'eom\'etrique, lien avec les syst\`emes de complexit\'e $2n+1$},
  Bulletin de la Soci\'et\'e Math\'ematique de France \textbf{123} (1995),
  no.~2, 271--292 (fr).

\bibitem[Fer96]{ferenczi1996rank}
\bysame, \emph{Rank and symbolic complexity}, Ergodic Theory Dynam. Systems
  \textbf{16} (1996), no.~4, 663--682. \MR{1406427}

\bibitem[FGH{\etalchar{+}}23]{fghsw21}
M.~Foreman, S.~Gao, A.~Hill, C.E. Silva, and B.~Weiss, \emph{Rank one
  transformations, odometers and finite factors}, Israel J. Math \textbf{255}
  (2023), 231--249.

\bibitem[Kal84]{kalikow}
Steven~Arthur Kalikow, \emph{Twofold mixing implies threefold mixing for rank
  one transformations}, Ergodic Theory Dynam. Systems \textbf{4} (1984), no.~2,
  237--259. \MR{766104}

\bibitem[Kno54]{knopp}
Konrad Knopp, \emph{Theory and application of infinite series}, Blackie and Son
  Limited, 1954.

\bibitem[Ler12]{leroy2}
Julien Leroy, \emph{Some improvements of the {$S$}-adic conjecture}, Adv. in
  Appl. Math. \textbf{48} (2012), no.~1, 79--98. \MR{2845508}

\bibitem[PS22]{PS2}
Ronnie Pavlov and Scott Schmeiding, \emph{Local finiteness and automorphism
  groups of low complexity subshifts}, Ergodic Theory and Dynamical Systems
  \textbf{43} (2022), 1980--2001.

\bibitem[PS23]{PS}
\bysame, \emph{On the structure of generic subshifts}, Nonlinearity \textbf{36}
  (2023), no.~9, 4904--4953.

\bibitem[Ryz93]{Ry93}
V.~V. Ryzhikov, \emph{Joinings and multiple mixing of the actions of finite
  rank}, Funktsional. Anal. i Prilozhen. \textbf{27} (1993), no.~2, 63--78, 96.
  \MR{1251168}

\bibitem[Sil08]{silva2008invitation}
C.~E. Silva, \emph{Invitation to ergodic theory}, Student Mathematical Library,
  vol.~42, American Mathematical Society, Providence, RI, 2008. \MR{2371216}

\end{thebibliography}
\end{document}